\documentclass{article}

\usepackage{ amssymb, amsmath, amsthm}

\usepackage{graphicx, psfrag, setspace, subfigure, gensymb}
\usepackage{tabularx}
\usepackage{xltabular}
\usepackage{thm-restate}
\usepackage{amsfonts}
\usepackage{bbm}
\usepackage{fix-cm}
\usepackage{comment}
\usepackage{enumitem}
\usepackage[mathcal]{eucal}
\usepackage{tikz,tikz-cd}
\usepackage{pgfplots}
\usepackage{booktabs}
\usetikzlibrary{pgfplots.fillbetween}
 \usetikzlibrary{tikzmark, arrows.meta}
\usetikzlibrary{matrix,arrows,decorations.pathmorphing,decorations.pathreplacing}
\tikzset{commutative diagrams/diagrams={baseline=-2.5pt},commutative diagrams/arrow style=tikz}
\let\amsamp=&
\tikzcdset{scale cd/.style={every label/.append style={scale=#1},
    cells={nodes={scale=#1}}}}

\usepackage{xcolor}
\usepackage[
  colorlinks=true,
  citecolor=blue!50!black,
  linkcolor=purple!50!black  
]{hyperref}

\usepackage{float}
\usepackage{setspace} 
\usepackage{enumitem}
\usepackage[
  backend=biber,
  style=alphabetic
]{biblatex}
\usepackage{float}
\usepackage{caption}
\usepackage{mathtools}
\usepackage{listings}
\usepackage{quiver}
\usetikzlibrary{nfold}
\usetikzlibrary{braids,positioning}
\usetikzlibrary{decorations.markings}
\usepackage{svg}

\usepackage{color}

\definecolor{mygreen}{rgb}{0,0.6,0}
\definecolor{mygray}{rgb}{0.5,0.5,0.5}
\definecolor{mymauve}{rgb}{0.58,0,0.82}

\graphicspath{ {./images/} }
\usepackage[export]{adjustbox}

\DeclareFontFamily{U}{mathx}{}
\DeclareFontShape{U}{mathx}{m}{n}{<-> mathx10}{}
\DeclareSymbolFont{mathx}{U}{mathx}{m}{n}
\DeclareMathAccent{\widehat}{0}{mathx}{"70}
\DeclareMathAccent{\widecheck}{0}{mathx}{"71}

\newcommand{\beq}[1]{\begin{equation}\label{#1} }
\newcommand{\eeq}{\end{equation}}

\theoremstyle{plain}
\makeatletter
\@ifundefined{newcounteralias}{}{%
  \renewcommand\thmt@autorefsetup{\@xa\def\csname\thmt@envname autorefname\@xa\endcsname\@xa{\thmt@thmname}}%
}
\newtheorem{prop}{Proposition}[section]
\newtheorem*{prop*}{Proposition}
\newtheorem{thm}{Theorem}
\newtheorem{lem}{Lemma}
\newtheorem*{lemma*}{Lemma}
\newtheorem{cor}[prop]{Corollary}

\theoremstyle{remark}
\newtheorem{rem}[prop]{Remark}

\theoremstyle{definition}
\newtheorem{defn}[prop]{Definition}

\newtheorem{eg}[prop]{Example}

\makeatletter \@addtoreset{equation}{section} \makeatother

\theoremstyle{plain}
\newtheorem{keythm}{Theorem}

\newtheorem{keyconj}[keythm]{Conjecture}

\let\oldmarginpar\marginpar
\renewcommand\marginpar[1]{\-\oldmarginpar[\framebox{\setstretch{\marginparstretch}\begin{minipage}{\marginparwidth}{\raggedleft\scriptsize #1}\end{minipage}}]{\framebox{\setstretch{\marginparstretch}\begin{minipage}{\marginparwidth}{\raggedright\scriptsize #1}\end{minipage}}}}

\let\oldmarginpar\marginpar
\renewcommand\marginpar[1]{\-\oldmarginpar[\raggedleft\footnotesize #1]%
	        {\raggedright\footnotesize #1}}

\newskip\stdskip                      
\begin{document}

\title{On the symplectic cohomology of orbifold divisor complements}
\author{Elliot Gathercole, Bogdan Simeonov}
\date{}

\maketitle

\setcounter{tocdepth}{1}
\hypersetup{
    bookmarksopen=true,
    bookmarksdepth=2
}
\begin{abstract}
    This article studies the symplectic cohomology of Liouville domains whose contact boundaries admit a periodic Reeb flow, with a particular focus on affine varieties which arise as orbifold divisor complements. We describe a formula for the Conley-Zehnder indices which allows for a simple determination of the positive equivariant symplectic cohomology in cases when the inertia orbifold of the compactifying divisor has only even degree cohomology. This is then applied to the study of orbifold compactifications of $\mathbb{C}^n$. Moreover, we analyse the action of multiplication by the principal orbit on V-shaped symplectic cohomology and prove that if the compactification $(\mathcal{X},\mathcal{D})$ is positively monotone, then this class induces a periodic structure on symplectic cohomology, reducing the problem of computing an infinite-dimensional group to a finite-dimensional repeating portion. We combine the index computations and this periodic structure to compute the symplectic cohomology of various Milnor fibers of surface singularities, the rational homology balls $B_{p,q}$ as well as ample divisor complements in orbifold K3 surfaces and Calabi-Yau threefolds, which includes the Milnor fibers of the exceptional unimodal singularities. 
\end{abstract}
\tableofcontents

\section{Introduction}

\subsection{Scope of the paper and main results}
The purpose of this paper is to study the symplectic cohomology (and its variants) of Liouville domains $W$ in the presence of a periodic Reeb flow on the contact boundary $\partial W$. Contact structures of this type are frequently referred to as \textit{Besse} (\cite{Kegel_Lange_2021}): the quotient of $\partial W$ by the Reeb circle flow is then an orbifold. A convenient framework to describe many such symplectic manifolds is by realising them, up to completion, as complements of a quasismooth orbifold divisor $\mathcal{X}\setminus \mathcal{D}$. The contact boundary in this model is represented by the circle Seifert bundle associated to the normal bundle of $\mathcal{D}$ in $\mathcal{X}$. 

\begin{eg}
Consider the Milnor fiber $M_f$ of a quasihomogeneous polynomial $f(z_1,\dots, z_n)$ of degree $d$ with weights $w_1,\dots, w_n$. The Liouville domain $M_f$ is given by $\{f=1\}\cap B_\delta\subset \mathbb{C}^n$, where $B_\delta$ is a ball centered at the origin of a large enough radius $\delta$ and the Liouville form is given by restricting the standard one on $\mathbb{C}^n$. Its completion can be compactified to $\{f=w^d\}\subset \mathbb{P}(w_1,\dots, w_n,1)$ by adding the compactifying divisor $\{f=0\}\subset \mathbb{P}(w_1,\dots, w_n)$. A concrete example is $f=z_1z_2+z_3^{p}$ with weights $(p-1,1,1)$, which defines the $A_{p-1}$ Milnor fiber. The compactifying curve in this case is a teardrop orbifold $\mathbb{P}(1,p-1)$ and the contact boundary (which is a lens space $L(p,p-1)$) is realised as a Seifert fibered $3$-manifold with one exceptional fiber over a $2$-sphere.
\end{eg}

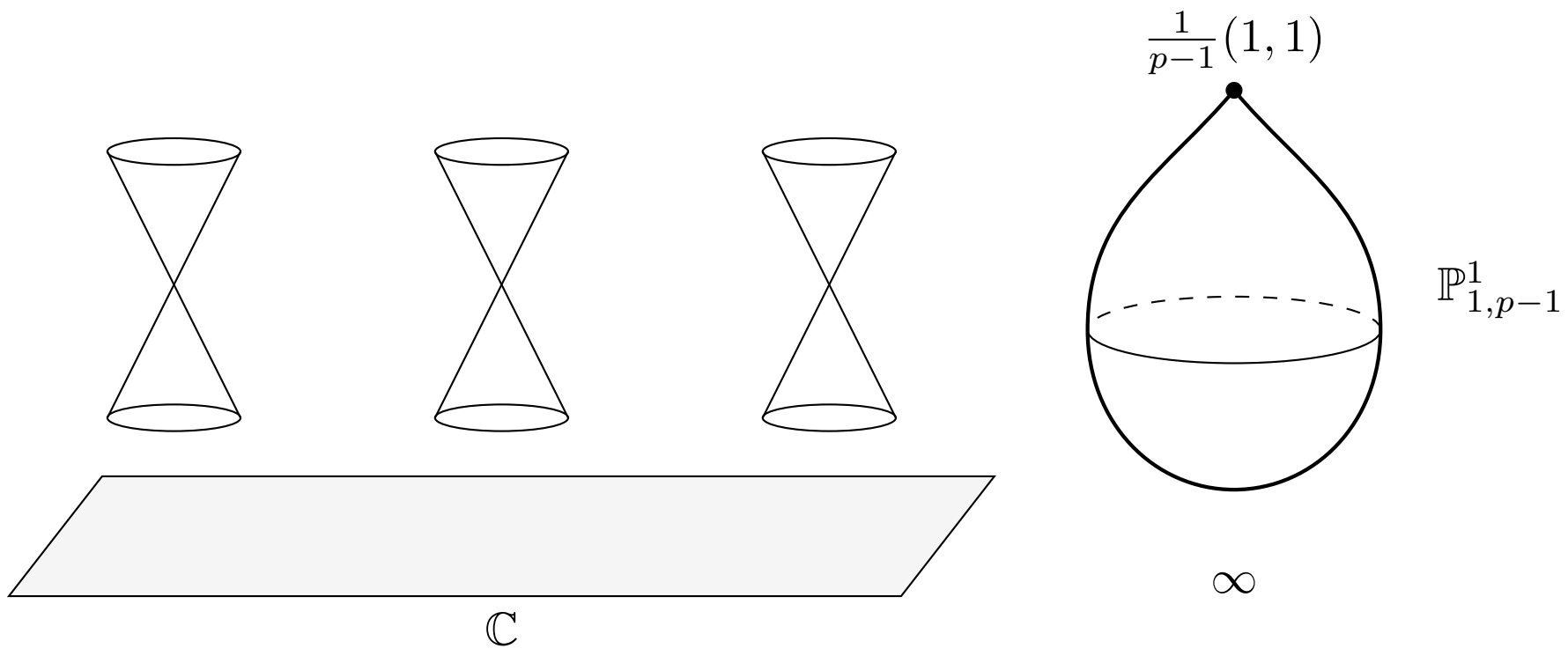
\begin{figure}[H]
\centering
\begin{tikzpicture}
  \draw[fill=gray!8]
    (-6.2,-2.3) -- (0.5,-2.3)
    -- (1.2,-1.4) -- (-5.5,-1.4) -- cycle;
  \node[below] at (-2.5,-2.3) {$\mathbb{C}$};

\begin{scope}[xshift=0,yshift=-70]
\draw (-3,3.5) node (v1) {} -- (-2,1.5);
\draw (-2,3.5) -- (-3,1.5);

\draw  (-2.5,3.5) ellipse (0.5 and 0.1);
\draw  (-2.5,1.5) ellipse (0.5 and 0.1);
\end{scope}

\begin{scope}[xshift=70,yshift=-70]
\draw (-3,3.5) node (v1) {} -- (-2,1.5);
\draw (-2,3.5) -- (-3,1.5);

\draw  (-2.5,3.5) ellipse (0.5 and 0.1);
\draw  (-2.5,1.5) ellipse (0.5 and 0.1);
\end{scope}

\begin{scope}[xshift=-70,yshift=-70]
\draw (-3,3.5) node (v1) {} -- (-2,1.5);
\draw (-2,3.5) -- (-3,1.5);

\draw  (-2.5,3.5) ellipse (0.5 and 0.1);
\draw  (-2.5,1.5) ellipse (0.5 and 0.1);
\end{scope}
  
  \begin{scope}[shift={(3,0)}]
    \draw[thick]
      (0,1.5)
      .. controls (-0.5,0.9) and (-1.1,0.6) .. (-1.1,-0.3)
      .. controls (-1.1,-1) and (-0.6,-1.5) .. (0,-1.5)
      .. controls (0.6,-1.5) and (1.1,-1) .. (1.1,-0.3)
      .. controls (1.1,0.6) and (0.5,0.9) .. (0,1.5);

    \draw[dashed]
      (1.1,-0.3) arc[start angle=0,end angle=180,
                     x radius=1.1,y radius=0.25];
    \draw
      (-1.1,-0.3) arc[start angle=180,end angle=360,
                      x radius=1.1,y radius=0.25];

    \fill (0,1.5) circle (1.8pt);
    \node[above] at (0,1.5) {$\frac{1}{p-1}(1,1)$};

    \node[right] at (1.4,0) {$\mathbb{P}^1_{1,p-1}$};
  \node[below] at (0,-2) {$\infty$};
  \end{scope}
\end{tikzpicture}
\caption{The $A_{p-1}$ Milnor fiber compactified, seen as a pencil with the compactifying divisor at $\infty$.}
\end{figure}

In theory, to compute symplectic cohomology, one needs to enumerate the Reeb orbits, compute their Conley-Zehnder indices and lastly compute the Floer cylinders between them. In practice, this is impossible unless one already has good control over the Reeb dynamics and even then, the computation of the Floer cylinders is extremely challenging. The work of Diogo and Lisi \cite{DiogoLisi2019} describes an algorithm to perform these three steps when $W=X\setminus D$ is a smooth divisor complement, in other words when the Reeb flow is an honest free circle action, under some monotonicity assumptions on the pair $(X,D)$, reducing the enumeration of Floer cylinders to Gromov-Witten theory. In the more general context of simple normal crossings compactifications, there is no decisive algorithm like the one for smooth divisor complements, but the Reeb dynamics at infinity can be explicitly described \cite{mclean2012growth}, and moreover there are powerful computational tools available like the spectral sequence of Ganatra and Pomerleano \cite{GP20}.

The guiding principle of this paper is that realising $W$ as an orbifold divisor complement allows one to work with the smallest, most efficient chain-level model of symplectic cohomology that has the fewest amount of families of Reeb orbits. Moreover, admitting a periodic Reeb flow carries consequences which might not be at all apparent if one works exclusively using simple normal crossings compactifications. One of the consequences is a repeating, periodic structure on the symplectic cohomology:

\begin{restatable}{keythm}{monotonecompactifications} \label{thm:monotonecompactifications}
Let $W$ be a Liouville domain which can be compactified to an orbifold with cyclic quotient singularities $(\mathcal{X},\mathcal{D})$, where $\mathcal{D}$ is smooth as an orbifold and whose isotropy locus is of codimension at least one. Suppose moreover that the compactification is monotone in the sense that $c_1^\mathcal{X}=\lambda \mathcal{D}, \lambda>1$.

Then, the class of the principal orbit $\check{M}_1\in \widecheck{\mathrm{SH}}^{-2(\lambda-1)}(W)$ induces a periodic structure on V-shaped symplectic cohomology
\begin{equation*}
    \widecheck{\mathrm{SH}}^\bullet(W)\simeq \widecheck{\mathrm{SH}}^{\bullet-2(\lambda-1)}(W).
\end{equation*}
Moreover, $\mathrm{SH}^\bullet(W) \simeq \widecheck{\mathrm{SH}}^\bullet(W)$ for all $\bullet <0$, so the same periodic structure persists for the negative degree part of symplectic cohomology.
    
\end{restatable}

Previously, it was noticed by Ritter \cite{ritter_floer_2014} that if a conical symplectic manifold $W$ admits a global Hamiltonian circle action generated by a Hamiltonian of linear slope at infinity, then there is a Seidel element $S_g(1)$ acting on its symplectic cohomology and inducing a periodic structure. Unfortunately, in our Liouville setup, if $\mathrm{SH}^\bullet(W)\neq 0$ then no such structure can exist. This led Uebele in \cite{uebele2019periodicreebflowsproducts} to consider the case when the circle action only exists on the contact boundary. Assuming quite strict index-positivity assumptions on $\partial W$ which imply that $\widecheck{\mathrm{SH}}$ is an invariant of the contact boundary, independent of the filling, the main result of \cite{uebele2019periodicreebflowsproducts} is the existence of a class $S_g(1) \in \widecheck{\mathrm{SH}}^\bullet(W)$ corresponding to the principal orbit which induces a periodic structure on V-shaped symplectic cohomology under multiplication.

The content of the theorem above is an extension of the results \cite{uebele2019periodicreebflowsproducts} which shows that a similar periodic structure exists on a much wider class of $W$ for which the symplectic cohomology depends on the filling. The notable case of the $A_{p-1}$ Milnor fiber is already mentioned in \cite[Section 4.4.2]{uebele2019periodicreebflowsproducts}: its contact boundary is not index-positive, yet the Hochschild cohomology computations of \cite{EtguLekili2017} imply that there is a periodic structure of degree $-2$. This was the motivating example for the ideas of this paper and we provide a new, purely geometric computation of the symplectic cohomology of the $A_{p-1}$ Milnor fiber in Section \ref{sec:Ap-1}.

It is worth noting that, typically, if $W$ admits a positively monotone compactification, then $\widecheck{\mathrm{SH}}(W)$ depends on $W$ and not just $Y:=\partial W$, and moreover it fails the strong index conditions from Uebele's work \cite{uebele2019periodicreebflowsproducts}, so our Theorem \ref{thm:mainthm} is a genuine improvement of the results of \cite{uebele2019periodicreebflowsproducts}.

The second point we want to stress is that the smooth orbifold divisor complement model can be much more efficient at describing the Reeb dynamics and, in particular, the positive equivariant symplectic cohomology.

\begin{eg}
    Consider again the $A_{p-1}$ Milnor fiber. Using the orbifold compactification chain model, there are exactly $p-1$ orbits contributing in degree $1$, namely one orbit $\hat{m}_1$ coming from the space of principal orbits, together with exceptional orbits $\hat{\delta}_{1/p-1}, \dots,\hat{\delta}_{(p-2)/(p-1)}$ of period less than $1$. In fact, $\mathrm{SH}^1\simeq \mathrm{SH}^1_{+,S^1}$ happens to be of rank $p-1$ so this is the most efficient possible chain level model. In contrast, the $A_{p-1}$ Milnor fiber also admits an snc compactification which is given by, for example, resolving the $\frac{1}{p-1}(1,1)$ point, which is isomorphic to $X=\mathrm{Bl}_{p}\mathbb{P}^2$, where the blowups all happen on a line (this is shown e.g. in \cite{evans2011symplectic}). Using the Ganatra-Pomerleano model \cite{GP20} of the compactification $(X, D_1\cup D_{-(p-1)})$ where $D_1,D_{-(p-1)}$ are strict transforms of lines, and using the fact that $-K_X=D_{-(p-1)}+2D_1$, one obtains an infinite set of generators of degree $1$ arising from multiple covers of the simple orbit in the circle normal bundle to the stratum $D_{-(p-1)}$ which is far from optimal.
\end{eg}

In some cases, as the one described above, the model is so small that in fact there is no differential on the first page Morse-Bott spectral sequence for positive equivariant symplectic cohomology. This is the content of our second result: following a detailed computation of Conley-Zehnder indices which involves the data of $c_1^\mathcal{X}$, as well as the age grading of cyclic quotient singularities, in Section \ref{sec:positiveequivariant} we show:
\begin{restatable}{keythm}{degeneration} \label{thm:degeneration}
    Suppose that $H^\bullet(I\mathcal{D}; \mathbb{Q})$ with the standard grading (not the Chen-Ruan grading) is concentrated in even degrees, where $I\mathcal{D}$ is the inertia orbifold of $\mathcal{D}$ thought of as a topological space. Then the Morse-Bott spectral sequence for positive equivariant symplectic cohomology degenerates on the first page. Moreover, if $H_1(W;\mathbb{Z})=0$ and $c_1^{TW}$ is torsion, then \begin{equation}
        \mathrm{SH}^\bullet_{+,S^1}(W)=\bigoplus_{N\geq 0} A[-2NI]=A\oplus A[-2I]\oplus\dots.
    \end{equation}
    where $2I$ is the index of the principal orbit and 
\begin{equation*}
 A:=\bigoplus_{S\in\pi_0I\mathcal{D}} H^\bullet(S;\mathbb{Q})\left[-\left(\mu_{CZ}(S)-(n-1)+\frac{1}{2}\dim_\mathbb{R}S\right)+1\right]   
\end{equation*}
are the contributions of the orbits of period $\leq 1$, with degree shifted down by the transverse Robbin-Salamon index.
\end{restatable}
\begin{rem}
    The last part of Theorem \ref{thm:degeneration} assumed that $c_1^{TW}$ is torsion and that $H_1(W;\mathbb{Z})=0$ in order to simplify the statement: under such conditions, the positive equivariant symplectic cohomology is $\mathbb{Z}$-graded. In general, one could either work with a fractional grading, or restrict attention to the contractible orbits.
\end{rem}

\begin{rem}
    Previous computations of Conley-Zehnder indices, for example for Brieskorn manifolds as in Kwon-van Koert's work \cite{kwon_brieskorn_2016} (see also \cite{Chiang2014Open} ) and Ustilovsky's \cite{Ustilovsky1999}, relied on an embedding of the contact manifold into $\mathbb{C}^n$ and an explicit trivialization of the contact structure, which does not easily generalize to an arbitrary orbifold divisor complement $W$.
\end{rem}

\begin{rem}
    Usually, the Morse-Bott spectral sequence for positive equivariant symplectic cohomology depends on a local system of coefficients. However, by \cite[Lemma 8.7]{kwon_brieskorn_2016}, this local system is trivial in the case that the linearised Reeb flow is complex-linear in some unitary trivialization of the contact structure along each $\gamma$ (or when $\gamma$ is an isolated good orbit). This condition is satisfied in our case, as the Reeb flow rotates the fiber of a prequantization bundle over a symplectic orbifold (in fact, the linearised Reeb flow can be made diagonal and is given by Equation \eqref{eq:linearizedReeb}, exactly as in \cite[Proposition 3.4]{LiZhou2026}).
\end{rem}

We briefly explain how to compute this in practice. For simplicity, let $\mathcal{X}$ and $\mathcal{D}$ have a single isolated cyclic quotient orbifold point of type $\frac{1}{w}(1,b_2,\dots, b_n)$ with $1$ corresponding to the direction normal to $\mathcal{D}$ and suppose $c_1^\mathcal{X}= \lambda \mathcal{D} \in H^2(X;\mathbb{Q})$ and $H_1(W)=0$. In this case, the inertia orbifold $I\mathcal{D}$ is simply \begin{equation*}I\mathcal{D}=\mathcal{D} \amalg \underbrace{[pt/\mu_w]\amalg \dots \amalg [pt/\mu_w]}_{w-1 \text{ copies}}\end{equation*}

The space of simple principal orbits is parametrized by $\mathcal{D}$, whereas each of the $[pt/\mu_w]$ corresponds to an exceptional orbit $\delta_{i/w}$ of period $<1$, where $1\leq i\leq w-1$. The relevant part of the degree shift is \begin{equation*}
    \mu_{CZ}(S)-(n-1)+\frac{1}{2}\dim_\mathbb{R} S=\begin{cases}
        2(\lambda-1), \,S=\mathcal{D}\\
        2\left(\frac{i}{w}(\lambda-1)-\frac{\sum_{j=2}^n\overline{ib_j}}{w}\right), \,S \text{ corresponds to the } i\text{-th copy of }[pt/\mu_w] 
    \end{cases}
\end{equation*}
where $\overline{ib_j}$ denotes the residue modulo $w$.

\begin{eg}
Consider again the compactification of the $A_{p-1}$ Milnor fiber into $\mathcal{X}=\{z_1z_2+z_3^p+w^p=0\}\subset \mathbb{P}(p-1,1,1,1)$. The divisor $\mathcal{D}\simeq \mathbb{P}^1_{1,p-1}$ is a hyperplane section and $c_1^\mathcal{X}=2\mathcal{D}$, showing that the simple principal orbits contribute $H^\bullet(S^2)[-2+1]$ and similarly the $m$-th covers of principal orbits will contribute $H^\bullet(S^2)[-2m+1]$. Moreover, the orbifold point is of type $\frac{1}{p-1}(1,1)$ so that the exceptional orbits of type $\delta_{i/{p-1}}$ with $1\leq i \leq p-2$ contribute $H^\bullet(pt)[0+1]$, and their further iterations $\delta_{m/{p-1}}$ contribute more generally $H^\bullet(pt)[-2\lfloor \frac{m}{p-1}\rfloor+1]$. Adding all of these together, one recovers the fact that \begin{equation*}\mathrm{SH}^\bullet_{+,S^1}(W)=\begin{cases}
    \mathbb{Q}^{p-1}, \bullet =1\\
    \mathbb{Q}^p, \bullet <1 \text { odd}\\
    0, \text { otherwise}
\end{cases}\end{equation*}
which was computed in \cite{kwon_brieskorn_2016}.
    
\end{eg}

An interesting situation where the conditions of Theorem \ref{thm:degeneration} apply is in the case of Kahler orbifold compactifications of $\mathbb{C}^n$. A long-standing conjecture of Brenton and Morrow \cite{brenton1978compactifications} stating that any smooth such compactification with a smooth divisor at infinity must be $(\mathbb{P}^n, \mathbb{P}^{n-1})$ was recently resolved by Li and Zhou \cite{LiZhou2026} using symplectic techniques. At around the same time, a similar result resolving various open problems on homology cell compactifications \cite{Fujita1980} was also established in \cite{Li2026Compactification} using entirely algebraic tools. Both of these results depend crucially on the Kobayashi-Ochiai theorem, which states that a Fano variety of dimension $n$ and Fano index greater than $n$ must be projective space.

However, there is no theorem analogous to the Kobayashi-Ochiai criterion for weighted projective spaces, whose $\mathbb{Q}$-Cartier Fano index can be arbitrarily high. Hence, determining the Fano index alone is not enough to characterize weighted projective space. Motivated by this, we prove something stronger: the vanishing of symplectic cohomology implies the existence of unirulings of high anticanonical degree:

\begin{restatable}{keythm}{orbcompclassification}
\label{thm:orb_com_classification}
    Let $(\mathcal{X},\mathcal{D})$ be a Kähler orbifold compactification of $\mathbb{C}^n$ with $d$ isolated orbifold points of order $w_1,\dots, w_d$. Then:
    \begin{itemize}
        \item $\mathcal{X}$ is a $\mathbb{Q}$-homology projective space, projective, $\mathbb{Q}$-Fano, of Picard rank one of index $r=n+1+\sum_{j=1}^d (w_j-1)$. 
        \item There is a uniruling by rational curves $C\simeq \mathbb{P}(1,w_{\max})$ passing through the point with maximal isotropy order and of orbifold anticanonical degree $n+1$: $$\mathrm{orbdeg}\,C:=d_{g^{-1}}C=-K_\mathcal{X} \cdot C+n-\iota_{age}=n+1$$ where $n-\iota_{age}$ denotes the age grading of the opposite of the twisted sector where the orbifold point of $C$ lands in.
        \item If $d=1$, then $(\mathcal{X},\mathcal{D})=(\mathbb{P}(1^n,k),\mathbb{P}(1^{n-1},k))$ for a $k \in \mathbb{N}$.
    \end{itemize}
\end{restatable}

\begin{rem}
The idea to use vanishing of symplectic cohomology to construct unirulings originated in the work of Anne-Laure Biolley \cite{biolley2003cohomologie}, \cite{biolley2004floer} and similar arguments were then used in \cite[Section 2.1]{AlbersFrauenfelderOancea2017}, \cite[Theorem 5.4]{Zhou2021Symplectic}, \cite{Cho2022Symplectic}, \cite[Theorem 12]{Li_2023_Exact_CY}, \cite{pieloch2024sections} as well as \cite[Section 7]{abbondandolo2024closedcharacteristicsminimalaction}.
\end{rem}

While there is no orbifold counterpart of the Kobayashi-Ochiai theorem \cite{KobayashiOchiai}, it is conjectured in \cite{Li2025Minimal} that there is an orbifold analogue to the Miyaoka-Cho-Shepherd-Barron theorem \cite{Cho2002CharacterizationsOP} (see also \cite{Kebekus2002}) characterizing $n$-dimensional projective space by the property of the existence of a locally unsplit dominating family of rational curves of degree $n+1$.

Motivated by these results, we conjecture:
\begin{keyconj}\label{conj:compactifications}
    Each Kähler orbifold compactification of $\mathbb{C}^n$ with isolated cyclic quotient singularities and a single smooth orbifold divisor must be a weighted projective space $\mathbb{P}(1,w_1,\dots, w_n)$.
\end{keyconj}

Interpreted differently, this conjecture is about classifying the periodic i.e. quasi-regular Sasakian Reeb flows that can arise on a sphere $S^{2n-1}$ equipped with the standard contact structure.

\begin{rem}
    Classifying pseudofree circle actions on odd spheres is a classical question in topology. Fintushel and Stern \cite{Fintushel1987} conjectured that any pseudo-free circle action on $S^5$ admits at most three exceptional orbits (this is now known as the Montgomery-Yang problem, an algebraic variant of which was proposed by Kollár \cite{kollar_is_2008} and resolved recently in \cite{jo2026algebraic}). Of these, there are the linear actions producing weighted projective space, but also many non-linear actions as well. In higher dimensions, the story becomes even more complicated: Montgomery and Yang in  \cite{MontgomeryYang1971} show that there are pseudo-free actions on topological $7$-spheres with arbitrarily many exceptional orbits and this was later extended to higher dimensions by Petrie \cite{petrie1974equivariant}. However, we claim that imposing the constraint on the contact structure to be the standard one limits the number of possible exceptional orbits. It is an interesting question whether there is any non-standard contact structure on $S^{2n-1}$ which admits periodic Reeb flows with arbitrarily many exceptional orbits.
\end{rem}

A possible strategy to prove this conjecture might go as follows: as a first step, using the high Fano index, one could attempt to compute the Hilbert function of $\mathcal{X}$. Then, one needs to show that there is a quasismooth representative $\mathcal{Z}$ in a linear system $|w\mathcal{D}|$ which avoids the orbifold point of order $w$, for some $w$. Proceeding inductively, either restricting attention to this lower-dimensional $\mathcal{Z}$ or by passing to a cyclic cover ramified over $\mathcal{Z}$, will reduce the Fano index and the number of singularities. However, each of these steps is highly nontrivial. We have been able to carry all of them out only in the case $d=1$.

Finally, we mention some computations. Using a slightly enhanced version of Theorem \ref{thm:monotonecompactifications} we obtain in Section \ref{sec:K3} a computation of the symplectic cohomology of the complement of an orbifold rational curve in an orbifold $K3$ surface:
\begin{prop*} [K3 surface compactifications] Let $S$ be a $K3$ surface, and $D$ a star-shaped configuration of smooth rational curves in $S$. We can contract all components of $D$ except the central one to obtain an orbifold rational curve $\mathcal{D}=\mathbb{P}^1_{a_1,\dots, a_d}$ inside an orbifold $K3$ surface $\mathcal{S}$. Assume that $\mathcal{D}$ is ample and moreover that $H_1(\mathcal{S}\setminus \mathcal{D};\mathbb{Z})$ is generated by the class of a principal orbit (for example, if it is zero, this is satisfied). Let 
    \begin{equation*}
        a=\sum_{i=1}^da_i-d.
    \end{equation*}
    Then:
    \begin{equation*}
        \mathrm{SH}^\bullet(W)=\begin{cases}
            \mathbb{Q} &\bullet=0
            \\ 0 &\bullet=1
            \\ \mathbb{Q}^{22}&\bullet=2
            \\ \mathbb{Q}^{a}&\bullet=3
            \\\mathbb{Q}^{a+1}&\bullet\geq 4.
        \end{cases}
    \end{equation*}
\end{prop*}
This, in particular, applies to all of the Milnor fiber of exceptional unimodal singularities defined by an invertible polynomial $w$, in which case the number $a+1$ is equal to $22-\mu(w)$.


\newcommand{\CentralRadius}{1.08}
\newcommand{\SmallRadius}{0.42}
\newcommand{\PictureUnit}{8.5mm}

\tikzset{
  sphere edge/.style={draw=black!85,line width=0.55pt},
  equator/.style={draw=black!44,line width=0.34pt},
  hidden equator/.style={equator,dash pattern=on 1.4pt off 1.3pt},
  neighbourhood/.style={draw=black!40,line width=0.5pt,
    dash pattern=on 3pt off 2.4pt,fill=black!2},
  curve label/.style={font=\small,inner sep=1pt},
  small curve label/.style={font=\scriptsize,inner sep=0pt},
  panel title/.style={font=\small},
  panel caption/.style={font=\footnotesize,align=center},
  orbifold point/.style={fill=black,draw=none},
}

\newcommand{\DrawSphere}[4]{%
  \begin{scope}[shift={#1}]
    \shade[sphere edge,top color=white,bottom color=black!12,
      shading angle=20] (0,0) circle[radius=#2];
    \draw[hidden equator] (#2,0)
      arc[start angle=0,end angle=180,x radius=#2,y radius={0.18*(#2)}];
    \draw[equator] (-#2,0)
      arc[start angle=180,end angle=360,x radius=#2,y radius={0.18*(#2)}];
    \node[#4] at (0,{0.50*(#2)}) {#3};
  \end{scope}%
}

\begin{figure}
\begin{tabular}{@{}c@{\hspace{4mm}}c@{\hspace{4mm}}c@{}}

\begin{tikzpicture}[x=\PictureUnit,y=\PictureUnit,baseline=(baseline)]
  \coordinate (baseline) at (0,0);
  \path[use as bounding box] (-4.85,-5.55) rectangle (4.85,2.8);
  \node[panel title] at (0,2.62) {Before contraction};


  \DrawSphere{(0,0)}{\CentralRadius}{
  }{curve label}

  \foreach \angle/\length in {90/1,0/4,-90/4,180/4}{
    \foreach \j in {1,...,\length}{
      \pgfmathsetmacro{\distance}{\CentralRadius+\SmallRadius+
        2*\SmallRadius*(\j-1)}
      \DrawSphere{(\angle:\distance)}{\SmallRadius}{$-2$}{small curve label}
      \pgfmathsetmacro{\joint}{\distance-\SmallRadius}
      \fill[black!80] (\angle:\joint) circle[radius=0.023];
    }
  }

\end{tikzpicture}
&
\begin{tikzpicture}[x=\PictureUnit,y=\PictureUnit,baseline=(baseline)]
  \coordinate (baseline) at (0,0);
  \path[use as bounding box] (0,-5.55) rectangle (1.8,2.8);
  \draw[-{Latex[length=2.8mm,width=1.7mm]},line width=0.65pt]
    (0,-0.82) -- (1.8,-0.82);
  \node[font=\normalsize] at (0.9,-0.44) {$\pi$};
  \node[font=\scriptsize,align=center] at (0.9,-1.32)
    {contract\\the legs};
\end{tikzpicture}
&
\begin{tikzpicture}[x=\PictureUnit,y=\PictureUnit,baseline=(baseline)]
  \coordinate (baseline) at (0,0);
  \path[use as bounding box] (-2.95,-5.55) rectangle (2.95,2.8);
  \node[panel title] at (0,2.62) {After contraction};

  \begin{scope}[shift={(0,-0.82)}]

    \shade[sphere edge,top color=white,bottom color=black!12,shading angle=20]
      (0,2.05)
      .. controls (0.10,1.65) and (0.35,1.45) .. (0.68,1.33)
      .. controls (1.00,1.20) and (1.20,1.00) .. (1.33,0.68)
      .. controls (1.45,0.35) and (1.65,0.10) .. (2.05,0)
      .. controls (1.65,-0.10) and (1.45,-0.35) .. (1.33,-0.68)
      .. controls (1.20,-1.00) and (1.00,-1.20) .. (0.68,-1.33)
      .. controls (0.35,-1.45) and (0.10,-1.65) .. (0,-2.05)
      .. controls (-0.10,-1.65) and (-0.35,-1.45) .. (-0.68,-1.33)
      .. controls (-1.00,-1.20) and (-1.20,-1.00) .. (-1.33,-0.68)
      .. controls (-1.45,-0.35) and (-1.65,-0.10) .. (-2.05,0)
      .. controls (-1.65,0.10) and (-1.45,0.35) .. (-1.33,0.68)
      .. controls (-1.20,1.00) and (-1.00,1.20) .. (-0.68,1.33)
      .. controls (-0.35,1.45) and (-0.10,1.65) .. (0,2.05)
      -- cycle;

    \draw[hidden equator] (-2.05,0)
      .. controls (-0.95,0.53) and (0.95,0.53) .. (2.05,0);
    \draw[equator] (-2.05,0)
      .. controls (-0.95,-0.53) and (0.95,-0.53) .. (2.05,0);
    \draw[hidden equator] (0,2.05)
      .. controls (-0.62,0.92) and (-0.62,-0.92) .. (0,-2.05);
    \draw[equator] (0,2.05)
      .. controls (0.62,0.92) and (0.62,-0.92) .. (0,-2.05);

    \node[curve label,fill=white,fill opacity=0.85,text opacity=1]
      at (0,0.13) {
      };
    \foreach \angle/\idx in {90/1,0/2,-90/3,180/4}{
      \fill[orbifold point] (\angle:2.05) circle[radius=0.042];
      \node[font=\small] at (\angle:2.37) {$p_{\idx}$};
    }
  \end{scope}

\end{tikzpicture}

\end{tabular}
\caption{Right: the K3 surface $\mathcal{X}_{15}\subset \mathbb{P}(2,3,5,5)$ contains a rational orbifold curve $\mathcal{C}_{3,15}\simeq \mathbb{P}^1_{2,5,5,5}$. The complement $\mathcal{X}\setminus \mathcal{C}$ is isomorphic to the quotient of the Milnor fiber of $x^5z+z^3+w^3$ by a free $\mu_3$ action. Left: describes the preimage of $\mathcal{D}$ in a crepant resolution of $\mathcal{X}$ which is a smooth $K3$ surface.}
\end{figure}
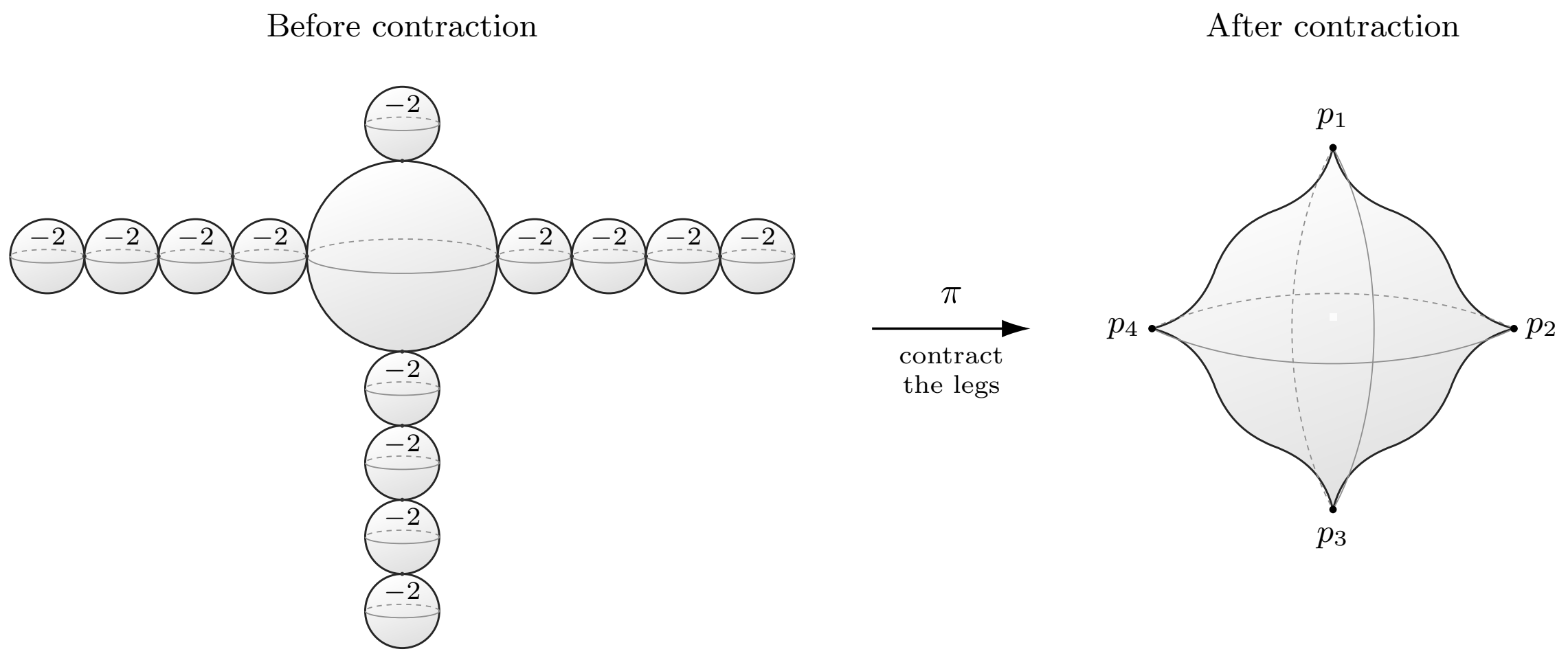

\begin{rem}
    In \cite[Conjecture 1.4]{LU22}, there is a general mirror symmetry conjecture describing the wrapped Fukaya category $\mathcal{W}(\check{V}_{w})$ of the Milnor fiber of an invertible polynomial $w$ as a matrix factorization category $\mathrm{mf}([\mathbb{A}^{n+1}/\Gamma_{w^T}], w^T+x_0x_1\dots x_n)$ using the transpose polynomial $w^T$. In that paper, Lekili and Ueda compute $\mathrm{HH}^\bullet(\mathrm{mf}([\mathbb{A}^{n+1}/\Gamma_{w^T}], w^T))$ in the case when $w^T$ (and hence $w$ as well) defines an exceptional unimodal singularity, and show in particular that it has ranks $(1,1,23)$ in degrees $(0,1,2)$. It is likely (though we do not show this) that deforming $w$ by the term $x_0\dots x_n$ will change the Hochschild cohomology to match with our computation above.
\end{rem}

We also compute in Section \ref{sec:CY3} the symplectic cohomology of orbifold divisor complements in Calabi-Yau threefolds with orbifold singularities.

\begin{prop*}[CY3 compactifications]
Suppose $W$ is a smooth, simply-connected Liouville domain which compactifies, via a simply-connected divisor $\mathcal{D}$, to a three-dimensional orbifold $\mathcal{X}$ with isolated orbifold points and $c_1^\mathcal{X}=0$. Suppose that:\begin{itemize}
        \item $H^4(W)=0$
        \item $H^3\rightarrow \mathrm{SH}^3(W)$ is injective
        \item every Reeb orbit of $\partial W$ (with its canonical contact structure as a Seifert bundle over $\mathcal{D}$) corresponding to a twisted sector of $\mathcal{X}$ of age one has smaller fractional period than the Reeb orbits corresponding to twisted sectors of age two.
    \end{itemize}  Then: $$\mathrm{SH}^\bullet(W)=\begin{cases}
        0, \bullet <0 \text{ or } \bullet =1\\
        \mathbb{Q}, \bullet =0\\
        \mathbb{Q}^{b_2(W)+1+\alpha_1}, \bullet =2\\
        \mathbb{Q}^{b_3(W)+\alpha_1}, \bullet =3\\
        \mathbb{Q}^{b_2(\mathcal{D})+2\alpha_1}, \bullet =4 \text{ or } \bullet \geq 6\\
        \mathbb{Q}^{b_2(\mathcal{D})+2\alpha_1-1}, \bullet =5
    \end{cases}$$
    where $\alpha_1$ is the number of non-identity twisted sectors on $\mathcal{X}$ with age $1$.
\end{prop*}

\begin{eg}
    Take the Calabi-Yau threefold defined by $\mathcal{X}_7\subset \mathbb{P}(1,1,1,1,3)$ relative to a hyperplane $\mathcal{D}_7$. This satisfies all the assumptions, since the complement $W=\mathcal{X}_7\setminus \mathcal{D}_7$ is identified with the Milnor fiber of the same polynomial that defines $\mathcal{D}_7$, which we could take, for example, to be $x^7+y^7+z^7+xu^2$. The Lagrangian $3$-spheres in this Milnor fiber span $H_3(W)$ and ensure, via Viterbo functoriality, that $H^3(W)$ injects into $\mathrm{SH}^3(W)$.
    Using Griffiths residues, one computes the Hodge numbers of $\mathcal{D}_7$ and in particular, $(b_0,b_2,b_4)=(1,43,1)$. Moreover, $\alpha_1=1$ and $b_2(W)=0, b_3(W)=288$ since the Milnor number is $288$.

    Therefore, 
    \begin{table}[H]
\begin{center}
\begin{tabular}{llcccccl}
    \toprule
     \(\mathrm{SH}^{<0}\) & \(\mathrm{SH}^{0}\)  &\(\mathrm{SH}^{1}\)          &\(\mathrm{SH}^2\)          &\(\mathrm{SH}^3\)          & \(\mathrm{SH}^4\) & \(\mathrm{SH}^{5}\) & \(\mathrm{SH}^{\geq 6}\) \\
    \midrule
     0& \(\mathbb{Q} \)                 & \(0\)             &  \(\mathbb{Q} ^2\)          &\(\mathbb{Q} ^{289}\)  & \(\mathbb{Q} ^{45}\) & \(\mathbb{Q} ^{44}\) & \(\mathbb{Q}^{45}\)\\
    \bottomrule
  \end{tabular}
\end{center}
\end{table}
\end{eg}

Finally, we mention also some log-Fano computations which are carried out in Section \ref{sec:logFanocomputations}.
\begin{prop*}
Consider $$f_m=x^2z+y^2w+z^{2m+1}+w^{2m+1}$$ which defines a canonical three-fold singularity.
The Milnor fiber $W_{2m+1}$ of $f_m$ compactifies as $\mathcal{X}_{2m+1}\subset\mathbb{P}(m,m,1,1,1)$ and the compactifying divisor $\mathcal{D}_{2m+1}$ is an orbifold surface whose stabilizer locus is the locus $z=w=0$ and forms a whole $\mathbb{P}^1$ with $\frac{1}{m}(0,1,1)$ singularities. The symplectic cohomology of $W_{2m+1}$ is given by $$\mathrm{SH}^\bullet(W_{2m+1})=\begin{cases}
        \mathbb{Q}^{m+(2m+2)^2}, \bullet =3\\
        \mathbb{Q}^{3m+4},\bullet \leq 1\\
        0, \text{ otherwise }
    \end{cases}$$
\end{prop*}

Moreover, positive equivariant symplectic cohomology can distinguish contact structures, like in the following family of examples:
\begin{eg}
    Consider the cone on the Johnson-Kollár surface (\cite{JK01}, \cite{Gugiatti2023FullEC}) $X_{8k+4}\subset \mathbb{P}(2,2k+1,2k+1,4k+1)$. This is a canonical 3-fold singularity with Milnor number $\mu=3(12k+9)$. Moreover, the link of this singularity is diffeomorphic to $(S^2\times S^3)\#_7$, but our computation of Conley-Zehnder indices shows that the positive equivariant symplectic cohomology depends on $k$ and is lacunary and hence is an invariant of the contact structure of the boundary (\cite{gutt2017positive}).
    \begin{equation}
    \mathrm{SH}_{+,S^1}^\bullet\cong\begin{cases}
        \mathbb{Q}^{6k+1}& \bullet=3,
        \\\mathbb{Q}^{12k+9}& \bullet=1,
        \\\mathbb{Q}^{12k+10}& \bullet=-2r-1,\;r\in\mathbb{Z}_{\geq 0}.
    \end{cases}
    \end{equation}
    Therefore, this gives an infinite family of inequivalent non-index-positive contact structures on $(S^2\times S^3)\#_7$ and their mean Euler characteristic satisfies \begin{equation}
        \chi_m=-\frac{1}{2} \mathrm{rk}\,H^\bullet_{CR}(X_{8k+4})=-(6k+5)
    \end{equation} 
\end{eg}

\subsection{Main ideas and proof strategies}

\subsubsection{Positive equivariant symplectic cohomology}
The main ingredient in the proof of Theorem \ref{thm:degeneration} is a computation of the Conley-Zehnder index in the context of a filling of a prequantization bundle over an orbifold. This is done by using an orbifold disk cap in the normal bundle of $\mathcal{D}$. However, there is a subtlety: orbifold vector bundles over an orbifold disk are not necessarily trivial. One can instead apply the desingularization procedure of Chen-Ruan \cite{chen_new_2004} to get an honest vector bundle over the underlying smooth disk of the orbifold disk. Using this desingularized frame, we compute the index $\mu(\gamma, c_{orb})$ relative to the trivialization of $TW$ induced by the desingularized cap $c_{orb}$.

This is enough, for one, to compute $\mu(\gamma)$ modulo $2$. In the case when $\gamma$ is null-homologous in $W$ and $c_1^{TW}$ is torsion, one can compute $\mu(\gamma)$ using a capping surface inside $W$ and the result is independent of the capping. The difference between $\mu(\gamma)$ and $\mu(\gamma, c_{orb})$ is measured by an orbifold Chern class. 
Using these ideas, in Section \ref{sec:grading_formulae} we show the following:
\begin{restatable*}{prop}{gradingformula}\label{prop:gradingformula}
    Suppose $\gamma$ is a Reeb orbit in $\partial W$, null-homologous in $W$, which belongs to a Morse-Bott family $\mathcal{N}$, with $\mathcal{N}/S^1\subset \partial W/S^1$ the corresponding locus of orbifold points on the divisor $\mathcal{D}$. Let $B$ be the class in $H_2(\mathcal{X};\mathbb{Q})$ of the surface obtained by gluing a normal fiber to a capping surface for $\gamma$ in $W$. Then, the transverse Robbin-Salamon index satisfies 
    \begin{equation}
     \mu_{CZ}(\gamma)=2(c_1^\mathcal{X}B-\iota_{age}-\lfloor B\cdot \mathcal{D}\rfloor)+n-1-\frac{1}{2}\dim_\mathbb{R} \mathcal{N}/S^1.
    \end{equation}
    where $\iota_{age}$ is an age grading contribution depending on the type of quotient singularity that the family of orbits $\gamma$ asymptote to.

    Perturbing the Morse-Bott Reeb dynamics to non-degenerate Hamiltonian dynamics has the following effect: after perturbation, which first breaks the Morse-Bott families of orbits into separate $S^1$ families indexed by critical points $q$ of a Morse function on $\mathcal{D}$, followed by a time-dependent perturbation breaking the $S^1$ families into a minimum and maximum, the Conley-Zehnder index is modified by $\mu_{Morse}+\mu_{S^1}-\frac{1}{2}\dim_\mathbb{R} \mathcal{N}/S^1$ and therefore
    \begin{equation}
        \mu_{CZ}(\gamma^q)=2(c_1^\mathcal{X}B-\iota_{age}-\lfloor B\cdot \mathcal{D}\rfloor)+n+(\mu_{S^1}-1)+(\mu_{Morse}(q)-\dim_\mathbb{R} \mathcal{N}/S^1).
    \end{equation}
    The corresponding grading in symplectic cohomology is given by $n-\mu_{CZ}$.
\end{restatable*}

Here, $\mu_{Morse}$ denotes the Morse index of a critical point of a function $f/S^1$ on $\partial W/S^1$ which comes from an $S^1$-invariant function $f$ on $\partial W$, whereas $\mu_{S^1}$ is $0$ or $1$, depending on whether $\gamma$ corresponds to a maximum or minimum after breaking the circle families of orbits. We will frequently use the notation $\check{\gamma}$ resp. $\hat{\gamma}$ to denote orbits with $\mu_{S^1}=1$ resp. $\mu_{S^1}=0$.

In particular, if the compactification is monotone so that $c_1^\mathcal{X}= \lambda \mathcal{D} \in H^2(\mathcal{X};\mathbb{Q})$, then the formula above does not depend on the spherical class $B$. In such a scenario, the class of the principal orbit associated to a maximum\footnote{or minimum, depending on convention} in its Morse-Bott family of generic orbits has index $2(\lambda-1)$. The exceptional orbits are more interesting:

\begin{eg}

Suppose $\mathcal{X}$ has an isolated $\frac{1}{m}(1,b_2,\dots, b_n)$ orbifold point on $\mathcal{D}$ (with the weight $1$ corresponding to the normal direction) and let $\delta_{i/m}$ denote a contractible Reeb orbit which is the $i$-th cover of the exceptional orbit over that point with $m \nmid i$. Then \begin{equation}
          \mu_{CZ}(\delta_{i/m})=2\left(c_1^\mathcal{X}B-\frac{\overline{i}+\overline{ib_2}+\dots+\overline{ib_n}}{m}-\left\lfloor \frac{i}{m}\right\rfloor\right)+n-1,
      \end{equation}
      where $\overline{a}$ denotes the residue of $a$ mod $m$. When $c_1^\mathcal{X}=\lambda \mathcal{D}$, this is equal to \begin{equation}
      \begin{split}
          \mu_{CZ}(\delta_{i/m})&=2\left(\lambda \frac{i}{m}-\frac{\overline{i}+\overline{ib_2}+\dots+\overline{ib_n}}{m}-\left\lfloor \frac{i}{m}\right\rfloor\right)+n-1
          \\&=\frac{2}{m}\left((\lambda-1)i-\left(\overline{ib_2}+\dots+\overline{ib_n}\right)\right)+n-1.
          \end{split}
      \end{equation}
\end{eg}

Equipped with these index formulas, one is in good shape to compute the positive equivariant symplectic cohomology by using the Morse-Bott spectral sequence. One can think of $\mathrm{SH}_{+,S^1}(W)$ as an unparametrized version of symplectic cohomology and each Morse-Bott family of orbits $\mathcal{N}\subset \partial W$ contributes $H^\bullet(\mathcal{N}/S^1;\mathbb{Q})$ shifted downwards by a suitable degree to the first page of the spectral sequence. In the case at hand, the spaces of orbits $\mathcal{N}/S^1$ are components of the inertia orbifold of $\mathcal{D}$. Theorem \ref{thm:degeneration} follows by showing that every contribution to the first page is actually in odd degree, by using the assumption that $H^\bullet(I\mathcal{D})$ is concentrated in even degrees only and applying the explicit computations of Conley-Zehnder indices.

\subsubsection{Periodicity}

There is the following fundamental obstacle in directly extending the framework of \cite{uebele2019periodicreebflowsproducts} and \cite{ritter_floer_2014} to our setting: since the rotation by the Reeb flow does not necessarily extend to the interior and since the contact boundary is not necessarily index-positive, there is no homology-level Seidel element. In other words, the rotation by the Reeb flow on the chain level does not necessarily commute with the differential.

To ameliorate this, we instead work with a replacement of the Seidel element, to first order: this is the class of the principal orbit. The downside of this is that one necessarily has to work on the chain level to make this work. The steps that underlie Theorem \ref{thm:monotonecompactifications} are to first show that $\check{M}_1$ is indeed closed and hence multiplication by it defines a chain map between suitable Floer complexes. Secondly, one studies the action filtration: by a monotonicity argument, we show that any pair of pants which crosses into the interior must consume a given energy penalty $\hbar$. Passing to the associated graded, we show that multiplication by $\check{M}_1$ has the same effect as rotating by the Reeb flow, and thus induces a filtered quasi-isomorphism. Informally, there is the following equation on the chain level:$$\check{M}_1\cdot - = \mathcal{S}(-)+\text{higher action corrections}$$

\begin{eg}[Motivating example]
Consider $T^*S^2$ and its compactification as a quadric minus a hyperplane section i.e. $\mathbb{P}^1\times \mathbb{P}^1\setminus \Delta$. The symplectic cochains $SC^\bullet(T^*S^2)$ have a chain-level model as described in Diogo's thesis \cite{Diogo2012}:

\begin{table}[H]
\begin{center}
\begin{tabular}{llccccl}
    \toprule
     \dots & \(SC^{-2}\)  &\(SC^{-1}\)          &\(SC^0\)          &\(SC^1\)          & \(SC^2_{}\) &  \\
    \midrule
     \dots& \( \check{M}_1\)  &\(\hat{M}_1\)  & &  & \\
     
  \dots& \( \check{m}_2\)  &\(\hat{m}_2\)  &  \(\check{m}_1\)& \(\hat{m}_1\) & \\
  
     \dots&  &\(\)  &  \(e\) &  &  \(S\) \\
    \bottomrule
  \end{tabular}
\end{center}
\end{table}
    Here $\check{M}_1$ is the principal orbit. Diogo enumerates the various pairs of pants and shows that the chain-level multiplication satisfies $\check{M}_1 \cdot \hat{m}_1=\hat{m}_2-\hat{M}_1$: in other words, the lowest-energy contribution is a rotation of $\hat{m}_1$ (which in fact corresponds to a degree zero Gromov-Witten invariant on $\Delta$), yet there is a  higher order correction $\hat{M}_1$ which corresponds to a Floer cylinder going deep into the interior owing to the fact that $T^*S^2$ is not index-positive. Nevertheless, there is still a $2$-periodic structure on the negative part of symplectic cohomology induced by $[\check{M}_1]$.
\end{eg}

The main result of Part \ref{part:Vshapedperiodicity} of the paper, which is a slightly more general version of Theorem \ref{thm:monotonecompactifications}, is proved in Section \ref{sec:main_thm}:

\begin{restatable}{keythm}{mainthm}
\label{thm:mainthm} Let $W$ be a Liouville domain with torsion $c_1^W$. Suppose that the Reeb flow on the contact boundary of $W$ is periodic and that there is a cocycle $\mathfrak{M}\in\mathrm{CF}^{-2I}(H)$ which agrees with the principal orbit $\check{M}_1^H$ to first order, meaning $\mathfrak{M}-\check{M}_1^H\in\mathrm{CF}^{-2I}_{[-1+\hbar,\infty)}(H)$, where $H$ is a J-shaped Hamiltonian as in Definition \ref{def:principalorbit}.
Then $\widecheck{\mathrm{SH}}^\bullet(W)$ is a module over $\mathbb{Q}[t^{\pm 1}]$, where $t\cdot x=[\mathfrak{M}]\cdot  x$ and there is an isomorphism
\begin{equation*}
    \widecheck{\mathrm{SH}}^\bullet(W)\simeq \widecheck{\mathrm{SH}}^{\bullet-2I}(W).
\end{equation*}
\end{restatable}

\begin{rem} 
    The reason we use $\mathfrak{M}$ instead of $\check{M}_1$ is that sometimes $\check{M}_1$ itself is not closed in our chain model, but after some higher action corrections it does become closed. An example where we describe such an $\mathfrak{M}$ is given in Section \ref{sec:K3}, where we study the $K3$ compactifications of the Milnor fibers of the exceptional unimodal singularities.
\end{rem}
There are various approaches to construct such a class $\mathfrak{M}$. In particular, we show in Section \ref{sec:constructingmathfrakM} that in the case of positively monotone compactifications (this notion of monotonicity is similar to the one in \cite{diogo2018morsebottsplitsymplectichomology}) with isolated cyclic quotient singularities, one does not even need the higher-order corrections to produce a closed representative:
\begin{restatable}{keycor}{logFanoCor}{(Positively monotone compactifications)}
\label{cor:log_Fano}
    Suppose $(\mathcal{X}, \mathcal{D})$ is a monotone compactification with cyclic quotient singularities on $\mathcal{D}$ and such that $ c_1^\mathcal{X}= \lambda \mathcal{D}\in H^2(X; \mathbb{Q})$ with $\lambda>1$. Then there exists a class $\mathfrak{M}$ as in Theorem \ref{thm:mainthm} and moreover $-2I=-2(\lambda-1)$, which implies Theorem \ref{thm:monotonecompactifications}.
\end{restatable}

It is also possible to construct a closed class $\mathfrak{M}$ outside the positively monotone (or 'log Fano') case by studying the Morse-Bott-Floer spectral sequence computing symplectic cohomology. We explain a criterion that allows one to do this in Section \ref{sec:constructingmathfrakM}.

\subsection{Orbifold compactifications of $\mathbb{C}^n$}
We recall some results about compactifications of homology cells. We denote by $X$ and $D$ the underlying varieties of the orbifolds $\mathcal{X}$ and $\mathcal{D}$.
\begin{prop}(\cite[Proposition 2.10]{LiZhou2026}, also \cite[Lemma 2.1 and Theorem 2.2]{Li2026Compactification}, \cite[Theorem 2]{Fujita1980}, \cite{Sommese1976}) Suppose $(\mathcal{X},\mathcal{D})$ is an orbifold compactification of $\mathbb{C}^n$, with $\mathcal{D}$ an irreducible quasismooth divisor. Then:
$\mathcal{X}$ is Fano, projective with $\mathcal{D}$ embedded as a hypersurface. Moreover, $H^\bullet(X;\mathbb{Q})\simeq H^\bullet(\mathbb{P}^n;\mathbb{Q}), H^\bullet(D;\mathbb{Q})\simeq H^\bullet(\mathbb{P}^{n-1};\mathbb{Q})$ and $-K_\mathcal{X}=(r+1)\mathcal{O}(\mathcal{D})$ for some integer $r>0$.
\end{prop}

From this alone, we can determine $H^\bullet(I\mathcal{D};\mathbb{Q})$, as well as the grading shifts by Conley-Zehnder indices that appear in Theorem \ref{thm:degeneration}. Then we compare the well-known computation of $\mathrm{SH}_{+,S^1}^\bullet(\mathbb{C}^n)$ to one using the orbifold divisor complement model laid out in Theorem \ref{thm:degeneration}. Surprisingly, this comparison is enough to determine the $\mathbb{Q}$-Fano index entirely from the dimension and the isotropy orders $w_i$, completing the first part of Theorem \ref{thm:orb_com_classification}.

For the second part, we use the vanishing of symplectic cohomology $\mathrm{SH}^\bullet(\mathbb{C}^n)=0$ and compare it to our chain model depending on a hypothetical orbifold compactification. By a neck-stretching argument, we compare the counts of Floer cylinders which destroy the unit class to counts of $J_W$ holomorphic planes and, in turn, to a Gromov-Witten invariant $\langle [\mathrm{pt}], \mathbf{1}_{1/w_{max}}\rangle^\mathcal{X}_\beta$ which is shown to be non-vanishing. Thankfully, all the curves that we need to account for in this computation are simple and somewhere injective, so no multiply-covered transversality issues arise.

The third part of Theorem \ref{thm:orb_com_classification} relies on some birational geometry. One begins with the following characterization of $\mathbb{P}(1^n,k)$:
\begin{prop}\cite[Lemma 13]{Liu_2017}
Let $X$ be an $n$-dimensional $\mathbb{Q}$-Fano variety containing a prime divisor $Z\simeq \mathbb{P}^{n-1}$ in its smooth locus. If $\rho(X)=1$ and $n\geq 3$, then $X=\mathbb{P}(1^n,k)$ for some positive integer $k$.  
\end{prop}

Recall also the famous Kobayashi-Ochiai theorem:
\begin{thm}(Kobayashi-Ochiai \cite{KobayashiOchiai})
    Suppose $Z$ is a smooth, $(n-1)$-dimensional Fano variety of Fano index $n$. Then $Z\simeq\mathbb{P}^{n-1}$.
\end{thm}
Our strategy then is to produce a divisor $Z\simeq \mathbb{P}^{n-1}$ inside $X$. This requires showing that a certain line bundle is base-point-free. We sketch the proof:
\begin{itemize}
    \item Using the comparison between the two computations of $\mathrm{SH}^\bullet_{+,S^1}$, one can pin down the type of the unique singularity on $X$ as $\frac{1}{k}(1,...,1)$ and the Fano index as $n+k$ with $-K_{\mathcal{X}}=(n+k)\mathcal{L}$ and $\mathcal{L}\simeq \mathcal{O}(\mathcal{D})$.
    \item The next step is to find a smooth (as a variety) representative $Z$ in the linear system $|\mathcal{O}(k\mathcal{D})|$. If a smooth section existed, it would have Fano index $n$, which by Kobayashi-Ochiai's theorem \cite{KobayashiOchiai} guarantees that $Z\simeq \mathbb{P}^{n-1}$. 
    \item To find such a smooth $Z$, by Bertini's theorem, it is enough to show that $|\mathcal{O}(k\mathcal{D})|$ is base point free. One computes the Hilbert polynomial and uses the ruling of $\mathcal{X}$ by curves with $C\cdot \mathcal{D}=\frac{1}{k}$. Equipped with this, one shows that various cohomology groups vanish and a little bit of birational geometry suffices to complete the proof.
\end{itemize}

When the number $d$ of singular points is greater than one, however, it becomes more difficult to control the singularity types: the various $w_i$ influence each other. Determining the Hilbert function is essentially a linear algebra problem, but it becomes increasingly difficult as $d$ increases.

\subsection{Organisation of the paper}

The paper is divided into three parts:

Part \ref{part:positiveeqsh} of the paper focuses on positive equivariant symplectic cohomology. We describe a formula (Proposition \ref{prop:gradingformula}) for the transverse Conley-Zehnder index of Reeb orbits in prequantization bundles over orbifolds relative to a filling $W$. The data involved is purely algebro-geometric: it depends on the first Chern class of $\mathcal{X}$, as well as the age grading of the cyclic quotient singularities on $\mathcal{D}$. This computation allows, by means of a Morse-Bott spectral sequence, for an efficient computation of $\mathrm{SH}^\bullet_{+,S^1}(W)$ whenever, for example, $\mathcal{D}$ has isolated cyclic quotient singularities and $H^\bullet(D;\mathbb{Q})$ is concentrated in even degrees only.

Part \ref{part:Vshapedperiodicity} is devoted to studying (non-equivariant) symplectic cohomology and its V-shaped counterpart $\widecheck{\mathrm{SH}}^\bullet(W)$ (which is known to be isomorphic to Rabinowitz Floer cohomology when using field coefficients, as in \cite{CieliebakFrauenfelderOancea2010}). We first describe various monotonicity lemmas and energy estimates, as well as a suitable discrete action filtration, ensuring that all Floer cylinders and pairs of pants that contribute on the associated graded (the zeroth page of the filtration spectral sequence) have energy less than $\hbar$ and thus do not cross into the interior of $W$. Thus, on the level of the associated graded, every Floer solution is entirely contained in the symplectisation. We then proceed to relate $\check{M}_1\cdot$ to the Reeb rotation $\mathcal{S}$ on the associated graded, eventually showing that $\check{M}_1\cdot$ is a filtered quasi-isomorphism.

The final Part \ref{part:computations} of the paper is dedicated to applying our results to compute the symplectic cohomology of various examples, including the Milnor fibers of the exceptional unimodal singularities which compactify to orbifold K3 surfaces, as well as various rational homology balls including $B_{p,1}$.

\subsection{Acknowledgments}
B.S. would like to thank his supervisor Yankı Lekili, as well as Alessio Corti, Kyler Siegel, Jonny Evans, Nick Sheridan, Calum Spicer, Calum Crossley, Yuhan Sun and Luís Diogo for the opportunity to discuss mathematics and for providing valuable insights that shaped the content of this paper.

E.G. would like to thank Jonny Evans for useful discussions around the subjects of this paper, and the organisers of the workshop on partially wrapped Fukaya categories at Helmsley for providing the authors with a great environment to discuss mathematics.

The ideas underlying this paper were nascent before the advent of LLMs powerful enough to do mathematics well. Nevertheless, the authors used large language models for proofreading at the final stage of writing up this document, as well as for extensive literature searches (in particular, in the field of birational geometry), as well as for TikZ image generation.

B.S. is funded by a London School of Geometry and Number Theory–Imperial College London PhD studentship, which is supported by the Engineering and Physical Sciences Research Council [EP/S021590/1].

E.G. is supported by TÜBİTAK grant no. 225N333.

\section{Background}
\label{sec:Floer_theory}
\subsection{Setup, conventions and notations}
Let $W$ be a Liouville domain, in other words a compact manifold with connected boundary admitting a symplectic form $\omega=d\lambda$ such that $\alpha:=\lambda|_{\partial W}$ is a contact form. We denote by $\widehat{W}$ its completion. The Reeb vector field $R$ on $Y:=\partial W$ is defined by the conditions $\alpha(R)=1, \iota_R d\alpha=0$. For a Hamiltonian function $H:W\rightarrow \mathbb{R}$, the Hamiltonian vector field is defined via the formula $\omega(-,{X_H})=dH$. Under this convention, if this Hamiltonian is autonomous and linear at infinity, then $X_H$ agrees with $h'(\rho)R$ on $\mathbb{R}_{\rho> 0}\times \partial W$. Throughout the paper, we will focus on the case when the Reeb flow is periodic. A particularly nice case is when the Reeb flow generates a pseudo-free circle action, meaning that only an isolated set of orbits have a finite stabilizer. In this case, the quotient $\partial W/S^1$ is an orbifold with isolated orbifold points.

Throughout, we will work over $\mathbb{Q}$. The rational coefficients are needed for when we apply the Morse-Bott spectral sequence for $\mathrm{SH}_{+,S^1}(W)$, or the long exact sequence between $\mathrm{SH}$ and $\widecheck{\mathrm{SH}}$.
\begin{center}
\begingroup
\renewcommand{\arraystretch}{1.2}
\begin{xltabular}{0.95\linewidth}{@{}l X@{}}
\toprule
\textbf{Notation} & \textbf{Meaning} \\
\midrule
$\widecheck{\mathrm{SH}}(W)$
  & Denotes the V-shaped symplectic cohomology, which is isomorphic to Rabinowitz Floer cohomology $RFH(W)$.\\
$\mathcal{A}_H(\gamma)$
  & Denotes the action functional, which in our convention is given by $\int_0^1 H(\gamma(t)) -\int_\gamma \lambda$\\
$E_{geo}(u)$
  & Geometric energy, which is equal to $\frac{1}{2}\int_S |du-X_H\otimes \beta|^2 \mathrm{d}\mathrm{vol}$ \\
$E_{top}(u)$
  & Topological energy given by $\int_S \tilde{u}^* (\omega-d(H\beta))$ where $\tilde{u}$ is the graph of $u$. The difference is $E_{top}(u)-E_{geo}(u)=-\int_S \tilde{u}^*(d_S(H\beta))-||\overline{\partial} u||^2_{L^2}$ where the last term vanishes when $u$ is a Floer solution. Moreover, when $d\beta=0$, this is given by $-\int_S \tilde{u}^* d_SH\wedge \beta$ which is non-negative under the split-monotonicity condition $d_SH\wedge\beta\leq 0$.\\
$(\mathcal{X}, \mathcal{D})$ resp. $(X,D)$
  & Denotes the compactification of $W$ obtained by performing a symplectic cut of $\partial W$ using the periodic Reeb flow. The $\mathcal{X}$ is an orbifold, and $X$ its underlying space.\\
$c_1^{log}$
  & Denotes $c_1^\mathcal{X}-[\mathcal{D}]$\\
$\rho$
  & Symplectisation coordinate \\
$\hbar_{cross}, \hbar, N$
  & Energy bound for curves crossing into the interior of $W$, with $\hbar$ some fixed value $0<\hbar<\min(\hbar_{cross}, \frac{1}{2N})$ with $N$ denoting the lcm of all the denominators of periods of Reeb orbits \\
$a(p)$
  & An increasing function $a:\mathbb{Z}\rightarrow \mathbb{R}$ used to define a discrete action filtration on the Floer cochains.\\ 
$\check{M}_1^H, \mathfrak{M}$ 
  & Denotes the class of the principal orbit in some chain group $\mathrm{CF}(H,J)$, respectively a cocycle whose lowest action part is equal to $\check{M}_1^H$\\
$I\mathcal{D}$
  & Denotes the inertia orbifold of $\mathcal{D}$\\
$\iota_{age}(g)$
  & Denotes the age grading of a twisted sector on $\mathcal{X}$, which is usually left implicit.\\
$\frac{\overline{i}}{w}$
  & Denotes a fraction between $0$ and $1$, where $\overline{i}$ signifies the non-negative reduction of $i$ modulo $w$.\\
$H^\bullet_{CR}(\mathcal{D})$
  & Denotes the Chen-Ruan cohomology of an orbifold, which is as a vector space $H^\bullet(I\mathcal{D})$, but the gradings are shifted by the age grading.\\

$\mathrm{orbdeg}\,C$
  & Denotes the orbifold degree $-K_\mathcal{X}\cdot C+\iota_{age}(g^{-1})=-K_\mathcal{X}\cdot C+n-\iota_{age}(g)$ of a curve with domain $\mathbb{P}(1,w)$ into $\mathcal{X}$ which evaluates at the twisted sector corresponding to $g$\\
$\mathfrak{X}_H$
  & Denotes the set of Hamiltonian 1-orbits for a Hamiltonian $H$\\
$h(\rho)$
  & Denotes a profile function for Hamiltonians linear at infinity, satisfying $h(1)=-\epsilon$ very small.\\
$\mathcal{S}$
  & Denotes the rotation by the Reeb flow on $\partial W$.\\
\bottomrule
\end{xltabular}
\endgroup
\end{center}

\subsection{Floer data}
To define $\mathrm{SH}^\bullet(W)$ and $\widecheck{\mathrm{SH}}^\bullet(W)$, we restrict our attention to Hamiltonians $H\in C^\infty(\widehat{W})$ and almost-complex structures $J$ on $\widehat{W}$ such that
 \begin{equation}
  \label{dissipative}
 \begin{cases}
      H|_{\partial W}<0,
     \\ H-\mu\rho \text{ is constant outside a compact set, for some $\mu>0$}.
     \\ J \text{ is of contact type outside a compact set}.
     \end{cases}
 \end{equation}

We consider regular $\mathbb{R}/\mathbb{Z}$-dependent Floer data $t\mapsto (H_t,J_t)$ such that for all $t\in\mathbb{R}/\mathbb{Z}$, $(H_t,J_t)$ satisfies \eqref{dissipative}. We can use this data to define the Floer complex $\mathrm{CF}^\bullet(H,J)$, the vector space which has as a basis
\begin{equation}
    \mathfrak{X}_H=\{\gamma:\mathbb{R}/\mathbb{Z}\rightarrow\widehat{W}\;|\; \dot{\gamma}=X_H\},
\end{equation}
the set of 1-periodic orbits of $X_H$, and a differential $d$. The coefficient of $\gamma_+$ in $d\gamma_-$ counts rigid elements of $\mathcal{M}(\gamma_+,\gamma_-;H,J)/\mathbb{R}$, where $\mathcal{M}(\gamma_+,\gamma_-;H,J)$ is the moduli space of solutions of Floer's equation
which carries an $\mathbb{R}$-action given by translation in the $s$-direction.

We say $H_-\leq H_+$ whenever $H_{-,t}(x)\leq H_{+,t}(x)$ for all $t\in\mathbb{R}/\mathbb{Z}$ and $x\in\widehat{W}$. A monotone homotopy between $(H_-,J_-)$  and $(H_+,J_+)$ is a regular $\mathbb{R}\times\mathbb{R}/\mathbb{Z}$-dependent Floer datum $(s,t)\mapsto (H_{s,t},J_{s,t})$ such that
\begin{equation}
\begin{cases}
    \partial_sH_{s,t}\leq 0,
    \\ (H_{s,t},J_{s,t})= (H_{\pm,t},J_{\pm,t}) $ for $ \pm s\ll 0,
    \\ (H_{s,t},J_{s,t}) $ satisfies \eqref{dissipative} for all $(s,t).
\end{cases}
\end{equation}
such a homotopy gives rise to a chain map $\mathrm{CF}^\bullet(H_-,J_-)\rightarrow \mathrm{CF}^\bullet(H_+,J_+)$, the continuation map, given by counting rigid elements of $\mathcal{M}(\gamma_+,\gamma_-;H,J)$, the moduli space of solutions to a parametrized Floer equation where $H,J$ now depend on $s$.

\subsection{Action filtration}
The action of an orbit $\gamma\in\mathfrak{X}_H$ is given by
\begin{equation}
    \mathcal{A}_H(\gamma)=\int_{0}^1H_t(\gamma(t))\,dt-\int_\gamma\lambda.
\end{equation}
The Floer differential, and continuation maps induced by monotone homotopies, increase action, so each Floer complex is naturally equipped with an exhausting $\mathbb{R}$-filtration under which the continuation maps become filtered chain maps. Given a possibly infinite interval $[a,b)$, we denote by $\mathrm{CF}_{[a,b)}^\bullet(H,J)$ the subquotient complex generated by orbits with action in $[a,b)$. 

For a strictly increasing function $a:\mathbb{Z}\rightarrow \mathbb{R}$, one can describe a decreasing discrete filtration $\mathcal{F}$ on the Floer cochains as
\begin{equation}
\label{eq:discrete_filtration}
    \mathcal{F}^p\mathrm{CF}^\bullet(H,J):=\{x\in \mathrm{CF}^\bullet(H,J)|\, \mathcal{A}_H(x)\geq a(p)\}
\end{equation}
where by $\mathcal{A}_{H}(x)$ we mean the minimal action that appears when writing $x$ as a sum of Hamiltonian orbits. Given such a filtration, we can define the associated graded, 
\begin{equation}
\label{eq:assoc_graded}
   \mathrm{gr}^p \mathrm{CF}^\bullet(H,J)=\mathcal{F}^p\mathrm{CF}^\bullet(H,J)/\mathcal{F}^{p+1}\mathrm{CF}^\bullet(H,J)=\mathrm{CF}_{[a(p),a(p+1))}^\bullet(H,J)
\end{equation}
for $p\in\mathbb{Z}$.

\subsection{The telescope model}
An acceleration datum consists of a cofinal sequence (with respect to the pointwise order on the Hamiltonians) of regular time-dependent Floer data $(H_n,J_n)$ $n\in\mathbb{Z}_{\geq 0}$, and regular monotone homotopies between $(H_n,J_n)$ and $(H_{n+1},J_{n+1})$ for each $n$. Given this data, we can define a 1-ray of filtered chain complexes,
\begin{equation}
    \mathcal{C}=\mathrm{CF}^\bullet(H_0)\rightarrow \mathrm{CF}^\bullet(H_1)\rightarrow \cdots \rightarrow \mathrm{CF}^\bullet(H_n)\rightarrow \cdots.
\end{equation}

The telescope $\mathrm{tel}^\bullet(\mathcal{C})$ is a filtered complex which is a model for the homotopy colimit. Following \cite{sun2026contactbigfibertheorems}, for a fixed $b$, the action filtered telescope $\mathrm{tel}^\bullet (\mathcal{C})_{(-\infty, b)}$ is a chain model for $\widecheck{\mathrm{SH}}_{(-\infty, b)}^\bullet(W)$ and the degreewise completed telescope $$\widehat{\mathrm{tel}}^\bullet(\mathcal{C}):=\varprojlim_{b\rightarrow \infty} \mathrm{tel}^\bullet (\mathcal{C})_{(-\infty,b)}$$ is therefore a chain model for the full V-shaped symplectic cohomology $\widecheck{\mathrm{SH}}^\bullet(W)$ defined in \cite{CieliebakFrauenfelderOancea2010}, which is isomorphic to Rabinowitz Floer homology $\mathrm{RFH}_\bullet$ if we take the ground ring to be a field, as shown in \cite[Theorem 1.5]{CieliebakFrauenfelderOancea2010} and \cite[Appendix A.2]{sun2026contactbigfibertheorems}. 

\subsection{Product operations on symplectic cohomology}
\label{pop_prod_def}
Below we describe various operations on symplectic cohomology as in \cite{ritter_topological_2013}, which count solutions to a Floer-type PDE whose domain is a punctured sphere. Importantly, the Hamiltonians used will be domain-dependent, interpolating between various fixed Hamiltonians at the cylindrical ends of the punctured surface.

Let $S=\mathbb{CP}^1\setminus \{z_0,z_1,...,z_d\}$ be a punctured sphere and $Z_+=[0,\infty)\times\mathbb{R}/\mathbb{Z},Z_-=(-\infty,0]\times\mathbb{R}/\mathbb{Z}$ be cylinders endowed with the complex structure coming from $s+\mathrm{i}t\mapsto (s,t)$. A positive (resp. negative) cylindrical end at $z_i$ is a holomorphic embedding $\epsilon_i:Z_\pm\rightarrow S$ such that $\epsilon_i(s,t)\to z_i$ as $s\to\pm\infty$. Suppose $S$ is equipped with cylindrical ends $\epsilon_i$ at $z_i$, positive for $i=1,...,d$ and negative for $i=0$. Choose integer weights $w_i\in\mathbb{N}$ such that $w_0=w_1+...+w_d$. 
\begin{rem}
    For our purposes, it is important that $w_i\in\mathbb{Z}$, though in general, any real weights can be used.
\end{rem}
Take $\beta\in\Omega^1(S)$ be a closed 1-form such that $\epsilon_i^*\beta=w_i\,\mathrm{d}t$ for $i=0,...,d$. Let $z\mapsto(H_z,J_z)$ be an $S$-dependent Floer datum such that for all $i=0,...,d$, and $\pm s\gg 0$,
\begin{equation}
\label{cyl_ends_consistency}
    w_iH_{\epsilon_i(s,t)}=H_{i,t},\quad J_{\epsilon_i(s,t)}=J_{i,t}
\end{equation}
for some regular $\mathbb{R}/\mathbb{Z}$-dependent $(H_i,J_i)$. In other words, this is an $S$-dependent family of Hamiltonians that interpolates between $d+1$ fixed Hamiltonians resp. almost complex structures. Given such data, we can define a generalised version of Floer's equation,
\begin{equation}
    \label{Floer_eq_pop}
\begin{cases}
    u:S\rightarrow \widehat{W}
    \\ (\mathrm{d}u-X_H\otimes\beta)j=J(\mathrm{d}u-X_H\otimes\beta)
    \\ \lim_{s\to\pm\infty}u(\epsilon_i(s,t))=\gamma_i(t).
\end{cases}
\end{equation}
In particular, when $d=2$, we can define (a choice of) the pair-of-pants product,
\begin{equation}
    \mu:\mathrm{CF}^\bullet(H_1,J_1)\otimes\mathrm{CF}^\bullet(H_2,J_2)\rightarrow  \mathrm{CF}^\bullet(H_0,J_0)
\end{equation}
where the coefficient of $\gamma_0$ in $\mu(\gamma_1,\gamma_2)$ counts elements of the moduli space $\mathcal{M}_{(S,\beta)}(\gamma_0,\gamma_1,\gamma_2;H,J)$ of solutions to \eqref{Floer_eq_pop}. By a standard argument, any such choice of product obeys the Leibniz rule and commutes with continuation maps up to chain homotopy. 

\subsubsection{Compatibility with action filtration}
If we additionally require that $H$ and $\beta$ satisfy the condition of split-monotonicity\footnote{Since we use the opposite sign convention for $\omega(-,X_H)=dH$ the sign of the inequality is reversed.}   (see \cite[p.~17]{abouzaid2024framede2structuresfloer}) i.e. for all fixed $x\in \widehat{W}$, 
\begin{equation}
\label{split_monotone}
    dH(x)\wedge\beta\leq 0,
\end{equation}
then the maps $\mu$ respect the action filtration in the sense that $$\mathcal{A}_{H_0}(\mu(x,y))\geq \mathcal{A}_{H_1}(x)+\mathcal{A}_{H_2}(y)$$ because (letting $\tilde{u}=(id,u):S\rightarrow S\times \widehat{W}$ be the graph of $u$)
\begin{equation}
    E_{top}(u)-E_{geo}(u)=-\int_S \tilde{u}^*(dH) \wedge \beta\geq 0.
\end{equation}

To describe split-monotone data, we need to ensure that for all $x\in\widehat{W}$,
\begin{equation*}
\label{eq:split_monotone_suff_cond}
    \min_{t\in\mathbb{R}/\mathbb{Z}}H_{0,t}(x)\geq 2\max_{t\in\mathbb{R}/\mathbb{Z}}H_{i,t}(x),\quad i=1,2
\end{equation*}
Then, setting $w_1=w_2=1$, fix some autonomous $K\in C^\infty(\widehat{W})$ such that $\max(H_1,H_2)\leq K\leq \frac{1}{2}H_0$, and choose $H$ such that for $i=0,1,2$, $H_{\epsilon_i(s,t)}$ is a monotone homotopy from $H_i$ to $K$, and $H=K$ outside the cylindrical ends. It follows from the above that the maps $\mu$ induce a well-defined associative product on $\widecheck{\mathrm{SH}}^\bullet(W)$. This product is also graded commutative.

\subsection{Units for $\widecheck{\mathrm{SH}}^\bullet$}
We describe a chain-level representative of the unit via a count of disks which defines the map $$\psi_{\mathbb{C}}:\mathbb{Q}\rightarrow \mathrm{CF}^\bullet(H)$$
analogously to Ritter's TQFT unit in \cite[Section 6]{ritter_topological_2013}.
\begin{defn}[PSS unit]
    Given an $\mathbb{R}/\mathbb{Z}$-dependent $H$, let $\psi_{\mathbb{C}}(1)=e^H\in\mathrm{CF}^0(H)$ be the element obtained by counting rigid disk solutions to the PSS map equation, which is a version of Floer's equation \eqref{Floer_eq_pop}, where $S=\mathbb{C}$, equipped with a single cylindrical end $e^{-2\pi(s+it)}$. The domain-dependent Hamiltonian and 1-form are given by $f(s)\,H$, respectively $f(s)\,\mathrm{d}t$, in the cylindrical end, extended by zero over the complement. Here, $f:\mathbb{R}_-\rightarrow\mathbb{R}$ is a function such that $f'(s)\leq 0$, $f(s)=0$ for $s\geq -1$ and $f(s)=1$ for $s\leq -2$. 
\end{defn}

The cohomology class $[e^H]$ is a unit for the pair-of-pants product, which can be shown by comparing it to the continuation maps (as in \cite[Lemma 6.2]{ritter_topological_2013}).

\subsection{Periodic Reeb flows and Morse-Bott perturbations}

\subsubsection{Reeb dynamics and inertia orbifolds}
By assumption, the Reeb flow is periodic \footnote{We will describe a converse construction, starting from an orbifold compactification $(\mathcal{X}, \mathcal{D})$ and producing a Liouville domain with periodic Reeb flow on the contact boundary in Section \ref{sec:positiveequivariant}}. The unparametrized Reeb orbits of period $\leq 1$ form families indexed by various strata of the orbifold $\partial W/S^1$. More precisely, there is an identification $$\{\text{Unparametrized Reeb orbits of } (\partial W,\alpha) \text{ of period }\leq 1\} \leftrightarrow \{\text{points on the inertia orbifold of }\partial W/S^1  \}$$

In other words, each Morse-Bott submanifold of Reeb orbits $\mathcal{N}\subset \mathbb{R}_+\times \partial W$ has an associated orbit space given by $\mathrm{proj}_2(\mathcal{N})/S^1\subset \partial W/S^1$. These orbit spaces stratify the orbifold $\partial W/S^1$, with the principal orbits forming the whole space, whereas the orbits with smaller period form various loci of orbifold points of $\partial W/S^1$. 

\begin{eg}
    Consider again the Milnor fiber $xy+z^p=1$ and the contact boundary which is Seifert fibered over a $\mathbb{P}^1_{1,p-1}$ teardrop orbifold. As such, the underlying manifold of the inertia orbifold is given by $$I\mathbb{P}^1_{1,p-1}=\mathbb{P}^1 \amalg \underbrace{pt\amalg \dots \amalg pt}_{p-2 \text{ copies}}$$
    The non-twisted sector $\mathbb{P}^1$ corresponds to the space of principal orbits, whereas the $p-2$ components correspond to the exceptional orbits $\delta_{1/p-1}, \dots, \delta_{(p-2)/(p-1)}$ of period $<1$.
\end{eg}

\subsubsection{V and J-shaped Hamiltonians and their orbits}

\begin{defn}[Autonomous J and V-shaped Hamiltonians]\label{def:profileautonomous}
A Hamiltonian $H:\widehat{W}\rightarrow \mathbb{R}$ is:
\begin{itemize}
    \item J-shaped if it is $C^2$ close to $h\circ \rho$, where $h$ is a piecewise linear function such that $-\epsilon\leq h(\rho)<0$ is constant for $\rho\leq 1$, and $h'(\rho)=\tau$ for $\rho>1$, some $\tau>0$ not in the period spectrum of $\partial W$, up to a convex corner smoothing at $1$. In other words, $H$ is a negative Morse function in $W$, linear with slope $\tau$ at $\infty$, and strictly convex between $\partial W$ and the linear region.
    \item V-shaped if it is $C^2$-close to $h\circ\rho$, where $h$ is a piecewise linear function such that $h'(\rho)\geq 0$ and $h(\rho)\geq b+\varepsilon$ for $\rho< 1-\delta$;  $h'(\rho)=-\eta$ for $1-\delta< \rho< 1$, $-\epsilon\leq h(1)< 0$; and $h'(\rho)=\mu$ for $\rho> 1$, for some $\mu,\eta>0$ not in the period spectrum  of $\partial W$, up to a corner smoothing which is convex at $1$ and concave at the other vertices. Moreover we require that $h(1)\in (-\epsilon,0)$. In other words, up to  $C^2$-small perturbation, in $W$, $H$ is a decreasing function of $\rho$; near $\partial W$, a negative, convex function of $\rho$; and at $\infty$, linear with slope $\mu$.
\end{itemize}
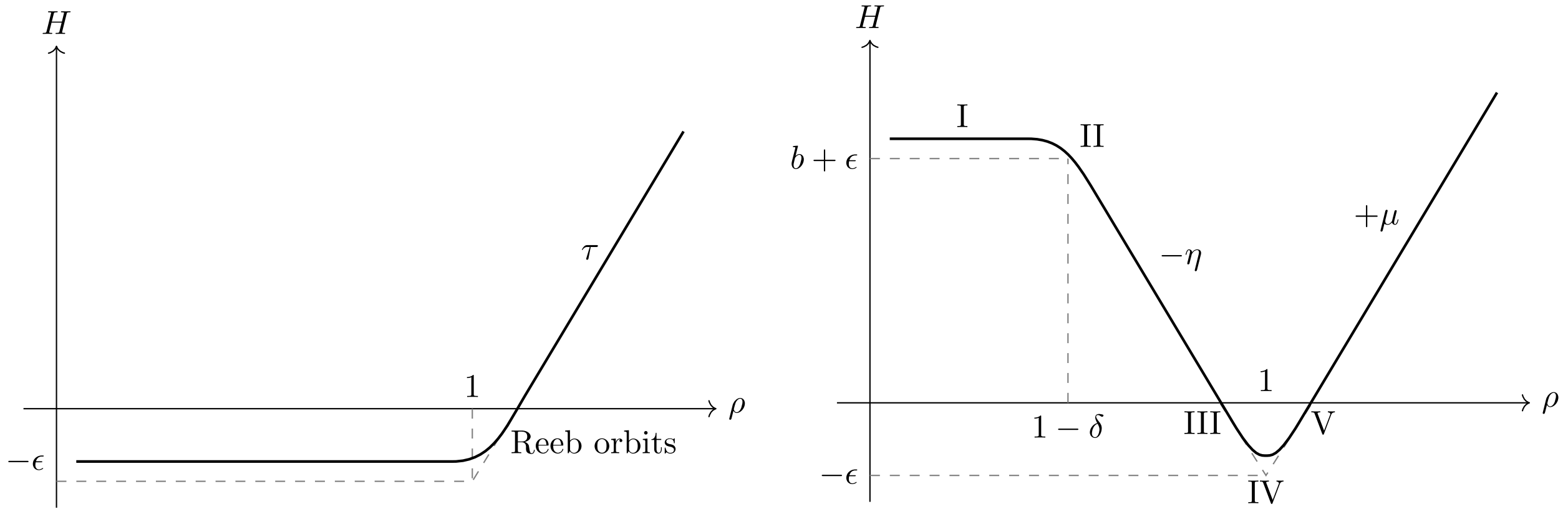
\begin{figure}[H]

\minipage{0.48 \textwidth}
\centering

\begin{tikzpicture}[scale=0.7]
  \draw[->] (-0.5,0) -- (10,0) node[right] {$\rho$};
  \draw[->] (0,-1.5) -- (0,5.5) node[above] {$H$};

  \draw[dashed,gray] (6.3,0) -- (6.3,-1.1);
  \draw[dashed,gray] (0,-1.1)--(6.3,-1.1)--(7.1,0.2);
  \draw[thick]
    (0.3,-0.8) -- (6,-0.8)
    .. controls (6.6,-0.8) and (6.8,-0.3) .. (7.1,0.2)
    -- (9.5,4.2);

  \node[above] at (6.3,0) {$1$};
  \node[left] at (0,-0.8) {$-\epsilon$};

  \node[left] at (8.4,2.37) {$\tau$};

  \node[right] at (6.7,-0.5) {Reeb orbits};
\end{tikzpicture}
\endminipage
\hfill
\minipage{0.48 \textwidth}
\centering


\begin{tikzpicture}[scale=0.7]
  \draw[->] (-0.5,0) -- (10,0) node[right] {$\rho$};
  \draw[->] (0,-1.5) -- (0,5.5) node[above] {$H$};

  \draw[dashed,gray] (3,0) -- (3,3.7);
  \draw[dashed,gray] (0,-1.1) -- (6,-1.1);
  \draw[dashed,gray] (0,3.7) --(3,3.7);

  \draw[dashed,gray] (5.3,0)--(6,-1.1)--(6.7,0);

  \draw[thick]
    (0.3,4)
    -- (2.4,4)
    .. controls (2.9,4) and (3.1,3.7) .. (3.4,3.2)
    -- (5.5,-0.3)
    .. controls (5.8,-0.8) and (5.9,-0.8) .. (6,-0.8)
    .. controls (6.1,-0.8) and (6.2,-0.8) .. (6.5,-0.3)
    -- (9.5,4.7);

  \node[below] at (3,0) {$1-\delta$};
  \node[above] at (6,0) {$1$};
  \node[left] at (0,-1.1) {$-\epsilon$};
  \node[left] at (0,3.7) {$b+\epsilon$};

  \node[right] at (4.2,2.2) {$-\eta$};
  \node[left]  at (8.2,2.8) {$+\mu$};

  \node[above] at (1.4,4) {$\mathrm{I}$};
  \node[above right] at (3,3.7) {$\mathrm{II}$};
  \node[left]  at (5.5,-0.3) {$\mathrm{III}$};
  \node[below] at (6,-1) {$\mathrm{IV}$};
  \node[right] at (6.5,-0.3) {$\mathrm{V}$};
\end{tikzpicture}
\endminipage
\caption{Left: a J-shaped Hamiltonian for symplectic cohomology. Right: a V-shaped Hamiltonian}
\end{figure}
 
\end{defn}

For $J$-shaped Hamiltonians, the Hamiltonian orbits are identified with:
\begin{itemize}
    \item Constant orbits in the interior of $W$
    \item Positively parametrized Reeb orbits for the contact structure $(\partial W, \alpha)$
\end{itemize}

In contrast, for V-shaped Hamiltonians, the orbits are identified with:
\begin{itemize}
    \item Constant orbits on $\partial W$ (type $\mathrm{IV}$)
    \item Positively resp. negatively parametrized Reeb orbits for the contact structure $(\partial W, \alpha)$ (types $\mathrm{V},\mathrm{III},$)
    \item Certain 'upper' orbits which are ultimately excluded based on their high action (types $\mathrm{I},\mathrm{II})$
\end{itemize}

The action of an orbit of type $\mathrm{I}$ or $\mathrm{II}$ is bounded below by some $b>0$. When defining the action-truncated symplectic cohomology, $\widecheck{\mathrm{SH}}^\bullet_{(-\infty,b)}(W)$, we can therefore choose a model defined using a cofinal family of V-shaped Hamiltonians with the same bound $b$, and neglect the upper orbits. Moreover, the value of $h(1)=-\epsilon$ will get closer and closer to $0$ in the limiting family.

The action of an orbit on a collar $\partial W \times \{\rho\}$ is given by $-\rho h'(\rho)+h(\rho)$. Hence, as in \cite[Section 3c]{seidel2008biased}, for positively parametrized orbits, a higher period corresponds to a bigger value of $\rho$ and a smaller action.

For a V-shaped Hamiltonian as above, the actions are as follows:
\begin{itemize}
    \item The constant orbits of type $\mathrm{IV}$ have action in $(-\epsilon,0)$, which can be chosen to be arbitrarily close to $0$.
    \item The positively parametrized orbits of type $\mathrm{V}$ and period $T$ have action in $(-T-\epsilon, -T)$, as in \cite[Section 5.4]{cieliebak2023loop}.
    \item The negatively parametrized orbits of type $\mathrm{III}$ and period $-T$ have action in $(T-\epsilon, T)$.
\end{itemize}

Hence, each non-upper orbit has action in an $\epsilon$-ball around a period $T$. In the next section, we will perform Morse-Bott perturbations of these autonomous Hamiltonians, which, provided the perturbation is small enough, do not affect the fact that the actions stay within a small distance relative to the periods.

\begin{defn}
    Given $b\in\mathbb{Z}_{\geq 0}$, a V-shaped Hamiltonian is $b$-admissible if all orbits of type $\mathrm{I},\mathrm{II}$ have action greater than $b+\epsilon$.
\end{defn}

Being $b$-admissible ensures that $\mathrm{CF}^\bullet_{(-\infty, b)}$ is generated by non-upper orbits. As the slope of a Hamiltonian grows, it becomes $b$-admissible for larger values of $b$. We will later fix $b$ and choose cofinal sequences of $b$-admissible Hamiltonians to compute $\widecheck{\mathrm{SH}}^\bullet_{(-\infty,b)}(W)$.

\begin{rem}In Section \ref{sec:estimates}, we describe a lower bound $\hbar$ on all Floer cylinders or pairs of pants which intersect the region $[\frac{1}{4},\frac{1}{2}]\times \partial W$. This lower bound depends only on the behaviour of $H$ there, and hence is entirely independent of $\epsilon$ which in turn depends on the behaviour of $H$ near $1$. We are thus free to choose $\epsilon \ll\hbar$, and we will impose this condition later on.\end{rem}

Finally, we note that if $u$ is either a Floer cylinder or pair of pants which is contained entirely in the symplectisation, then there is a relation between the periods $$T_{in}-T_{out}=\int u^* d\alpha\geq 0 \quad \text{resp. }\quad T_{in_1}+T_{in_2}-T_{out}=\int u^* d\alpha\geq 0 \qquad$$ since $d\alpha$ annihilates the Reeb and radial directions and $\int u^*d\alpha$ can be computed by projecting to the contact component (see for example \cite[Proposition 2.15, 2.17]{pascaleff2019symplectic}).

\subsubsection{Perturbations}
\label{sec:MB_perturbation}
We now describe how to perturb autonomous V and J-shaped Hamiltonians following \cite{bourgeois_morse-bott_2003}, \cite{van_koert_contact_2008} and \cite{ritter_filtrations_2025}. We will perturb $H$ separately near each Morse-Bott family $\mathcal{N}$ of Reeb orbits of period $T$, ensuring that the perturbed lower orbits have action in $[T-2\epsilon,T)$, and perturbed upper orbits have action $> b+\epsilon$.

In \cite{bourgeois_morse-bott_2003}, Bourgeois defines a function $f/S^1:\partial W/S^1 \rightarrow \mathbb{R}$ by inductively defining $S^1$-invariant functions $f_\mathcal{N}$ on the Morse-Bott manifolds $\mathcal{N}$, beginning with the manifold of Reeb orbits of smallest period. The induced functions $f_\mathcal{N}/S^1$ are Morse functions in the orbifold sense (see \cite[Sections 2.2, 2.3]{bourgeois_morse-bott_2003}, as well as \cite{10_2140_agt_2009_9_1105} for a review of Morse functions on orbifolds and their properties). Note that one needs to define these functions only on the manifolds of Reeb orbits of period $\leq 1$: for those of higher period, one can just repeat the same construction, since the spaces of Reeb orbits of period $\sigma$ are naturally identified with those of period $\sigma-1$ by rotating by the Reeb flow.

Whereas Bourgeois perturbs the contact form by replacing $\alpha$ by $(1+\lambda f)\alpha$, we will instead perturb our Hamiltonian. The effect is the same: our Hamiltonian orbits will break the Morse-Bott families of Reeb orbits into a collection of $S^1$ families of orbits, indexed by the critical points of $f/S^1$.

Define 
\begin{equation*}
    H'=H+ \sum \varepsilon_\mathcal{N}\rho_\mathcal{N}f_\mathcal{N},
\end{equation*}
where $\rho_\mathcal{N}$ is a bump function with support near $\mathcal{N}$ and $\varepsilon_\mathcal{N}>0$ is a small constant. The non-constant 1-orbits of $X_{H'}$ now come in $S^1$-families, with one circle family for each critical point of each $f_\mathcal{N}$.
\begin{rem}
    We are free to choose $\varepsilon_\mathcal{N}$ as small as we want.
\end{rem}

Finally, after performing a small time-dependent perturbation as in \cite[Remark A.1]{ritter_filtrations_2025} (or \cite[Proposition 2.2]{Cieliebak1996}), we end up with a Hamiltonian $\tilde{H}$ whose orbits are non-degenerate, with the circle families breaking into pairs of a maximum and minimum (whose index differs by one). The actions of the resulting lower orbits of $\tilde{H}$ can be made arbitrarily close to the periods of $(\partial W, \alpha)$ by taking the $\varepsilon_\mathcal{N}$ and the time-dependent perturbation to be arbitrarily small.

\begin{rem}
    After perturbation, the surviving 1-periodic orbits of $X_{\tilde{H}}$ are still preserved by the Reeb flow $R$.
\end{rem}

The Floer complex for a V-shaped $b$-admissible $\tilde{H}$ splits into a direct sum
\begin{equation*}
(\mathrm{CF}^\bullet_{[-\infty,b)}(\tilde{H}),d)=\left(\mathrm{CF}_{\lor}^\bullet(\tilde{H})\oplus \mathrm{CF}_{\wedge}^\bullet(\tilde{H}), \begin{pmatrix} d^{\lor, \lor} & d^{\lor, \wedge} \\ d^{\wedge, \lor} & d^{\wedge, \wedge}\end{pmatrix}\right)
\end{equation*}
where the summands are generated by the minima and maxima in the $S^1$-families, respectively, denoted type $\wedge$ and type $\lor$. For $J$-shaped $\tilde{H}$, there will be a summand corresponding to constant orbits in $W$.

The differential breaks into three pieces $d=d^{\lor, \wedge}+(d^{\lor, \lor}+d^{\wedge, \wedge})+d^{\wedge, \lor}$, where the superscript corresponds to whether the differential is connecting a minimum/maximum on the $S^1$ families (similar to the index filtration in \cite[Section 4]{bourgeois_exact_2009}).

\begin{prop}
\label{prop:perturbed_action_windows}
    If $H$ is a $b$-admissible V-shaped Hamiltonian, we can choose a non-degenerate perturbation $H$ such that the perturbed upper orbits have action $\geq b+\epsilon$, and perturbed lower orbits of period $T$ have action in $(-T-\epsilon,-T)$.
\end{prop}
\begin{proof}
    Both the autonomous perturbations of the Morse-Bott families, and the time-dependent perturbations, can be supported in a small neighbourhood of each Morse-Bott family of 1-periodic orbits, and can be made arbitrarily $C^2$-small in this neighbourhood, preserving the required (closed) conditions on the action.
\end{proof}

\subsubsection{Filtration by Period}
In this section, we show that the $\mathbb{R}$-filtration by action can be refined to a $\mathbb{Z}$-filtration by period, and that the operations respect this filtration. 

The Reeb flow, by assumption, is periodic. Therefore, the period spectrum is discrete and is contained entirely inside $\frac{1}{N}\mathbb{Z}$ for some large $N$. For example, we could take $N$ to be the product of all isotropy orders on the associated orbifold $\partial W/S^1$. As described in the previous section, we can choose perturbations $H$ of V-shaped Hamiltonians so that the actions of orbits of type $\mathrm{III},\mathrm{IV},\mathrm{V}$ arising from Reeb orbits of period $T$, all lie in windows of the form $(-T-\epsilon, -T)$ after the Morse-Bott perturbation.

From now on, we assume that $\hbar <\frac{1}{2N}$ is a fixed constant and $\epsilon < \hbar/8$. All our Hamiltonians will be required to satisfy $h(1)\in (-\epsilon,0)$ . We will later also impose the condition that $\hbar<\hbar_{cross}$, which is a positive constant determining the energy penalty upon crossing into the interior of $W$.

\begin{defn}[Action filtration compatible with period]
    Let
    \begin{equation}
        a(2p)=\frac{p}{N}-\frac{\hbar}{2},\quad a(2p+1)=\frac{p}{N}+\frac{\hbar}{2}
    \end{equation}
    for $p\in\mathbb{Z}$, and let $\mathcal{F}$ denote the decreasing action filtration on $\mathrm{CF}^\bullet(H,J)$ defined using the function $a$, as in Equation \eqref{eq:discrete_filtration}.
\end{defn}

We discuss now some relevant properties of this filtration and its compatibility with the product operation.

\begin{prop} Suppose $H$ is $b$-admissible. Then we have the following:
\label{prop:filtration_properties_list}
    \begin{enumerate}
        \item All non-upper orbits are contained in some action window $[a(2p), a(2p+1))$. In particular, $$\mathrm{CF}^\bullet_{[a(2p+1), b)}(H,J)= \mathrm{CF}^\bullet_{[a(2p+2), b)}(H,J)$$
        \item A shift by $-t$ does not affect the filtration up to orbits of action $>b$: $$\mathrm{CF}^\bullet_{[a(p), b)}(H,J)= \mathrm{CF}^\bullet_{[a(p)-t, b)}(H,J)$$ as long as $0<t<\hbar/4$.
        \item The product operation induces a well-defined operation $$\mathrm{CF}^\bullet_{[a(2p), a(2p+1))}(H_1,J_1)\otimes \mathrm{CF}^\bullet_{[a(2q), a(2q+1))}(H_2,J_2)\rightarrow \mathrm{CF}^\bullet_{[a(2p+2q), a(2p+2q+1))}(H_0,J_0)$$ as long as both $H_1$ and $H_2$ are $b$-admissible for some $b$ satisfying $$b\geq \max\{a(2p+2q+1)-a(2q),a(2p+2q+1)-a(2p)\} $$
    \end{enumerate}
\end{prop}

\begin{proof}
    The first claim follows by the fact that all actions of non-upper orbits are contained in some window $[p/N-\epsilon, p/N)$ for some $p$. The second claim follows by the choice of $\epsilon < \hbar/4$. 

    For the last claim, we argue as follows. First of all, two terms $x,y$ with actions $\geq a(2p), a(2q)$ respectively multiply to have action $\geq a(2p)+a(2q)=\frac{p+q}{N}-\hbar$. However, unless the product is zero, in which case the claim is trivial, a nonzero orbit of action $\geq \frac{p+q}{N}-\hbar$ necessarily has action $\geq \frac{p+q}{N}-\frac{\hbar}{2}=a(2p+2q)$ by our assumption that $\hbar < \frac{1}{2N}$.

    Now, we need to ensure that this is well-defined, up to terms $x_0,y_0$ with actions $\geq a(2p+1)$ and $\geq a(2q+1)$ respectively. The upper components of $x_0, y_0$ will contribute to action above $a(2p+2q+1)$ by our assumptions on $b$.
    
    By the first part of the proposition, any non-upper orbits contributing to $x_0,y_0$ will have action $\geq a(2p+2)$ resp. $\geq a(2q+2)$. More precisely, $\mathcal{A}(x_0)\geq \frac{p+1}{N}-2\epsilon, \mathcal{A}(y_0)\geq \frac{q+1}{N}-2\epsilon$. In particular, $$\mathcal{A}(x_0\cdot y+x\cdot y_0+x_0\cdot y_0)\geq \frac{p+q+1}{N}-4\epsilon\geq \frac{p+q}{N}+\frac{\hbar}{2}$$ so that indeed the product is well-defined up to corrections of action $\geq a(2p+2q+1)$.
\end{proof}
\begin{prop}
\label{prop:mu_respects_F}
    Let 
    \begin{equation}
        \mu:\mathrm{CF}^\bullet(H_1)\otimes \mathrm{CF}^\bullet(H_2)\rightarrow \mathrm{CF}^\bullet(H_0)
    \end{equation}
    be a chain-level product defined using split-monotone Floer data where each $H_i$ is $b_i$-admissible. Then $\mu$ respects the filtration $\mathcal{F}$ in a restricted sense: whenever $b_0\leq\min\{p/{2N}+b_2,q/2N+b_1\}$,
    \begin{equation}
        \mu(\mathcal{F}^p\mathrm{CF}_{(-\infty,b_1)}^\bullet(H_1)\otimes \mathcal{F}^q\mathrm{CF}_{(-\infty,b_2)}^\bullet(H_2))\subset\mathcal{F}^{p+q}\mathrm{CF}^\bullet_{(-\infty,b_0)}(H_0,J_0).
    \end{equation}
    Moreover, the chain homotopies witnessing compatibility of the product with continuation maps are filtered in this sense, provided they are defined using split-monotone Floer data.
\end{prop}
\begin{proof}
    The proof is similar to the previous proposition, in particular the fact that two orbits of actions $>a(p)$ and $>a(q)$ multiply to something with action $>a(p+q)$ is immediate by the same logic. Well-definedness follows from the assumption on $b_0$.
\end{proof}

\begin{prop}
\label{prop:unit_mult_cont_map}
Suppose that $H_\pm$ are $b$-admissible, $H_-\leq H_+$ and $2\max_{t\in\mathbb{R}/\mathbb{Z}}H_{-,t}< \min_{t\in\mathbb{R}/\mathbb{Z}}H_{+,t}$. Then there is a chain homotopy
\begin{equation}
\mu(e^{H_-},\;\cdot\;)\simeq\kappa:\mathrm{CF}^\bullet_{(-\infty,b)}(H_-)\rightarrow\mathrm{CF}_{(-\infty,b)}^\bullet(H_+)
\end{equation}
which is filtered with respect to $\mathcal{F}$, where $\mu$ is a chain-level multiplication arising from split-monotone Floer data on the pair-of-pants surface, and $\kappa$ is a continuation map arising from a monotone homotopy.
\end{prop}
\begin{proof}
The map $\mu(e^H,\;\cdot\;)$ counts configurations of a disc glued to a pair-of-pants along a 1-periodic orbit of $X_H$. Denoting by $P$ the domain for the pairs of pants and by $\mathbb{C}$ the disk for the PSS unit, by \cite[Theorem 16.12]{ritter_topological_2013}, composition of operations agrees on the chain level with gluings of surfaces: $$\psi_\mathbb{C}\circ \psi_P=\psi_{\mathbb{C}\#_{\lambda} P}$$where $\lambda\gg 0$ is a sufficiently large neck length. The surface $C\#_\lambda P$ is topologically a cylinder, however the Floer data is not a priori the same as that of the continuation cylinders.

A gluing argument, such as \cite[Theorem 16.10]{ritter_topological_2013}, gives a (not necessarily filtered) chain homotopy between $\psi_{\mathbb{C}\#_\lambda P}$ and $\kappa$, obtained by counting solutions to Floer's equation for a parametrised family of Floer data on the cylinder, which interpolates between the datum defining $\kappa$ and the datum obtained by gluing the Floer data defining $\mu$ on the pair-of-pants to the PSS datum on the disc. It remains to show that we can arrange for the solutions $u$ to the resulting family of Floer's equations all to have topological energy $E_{top}(u)\geq -\hbar/4$, which would, by Proposition \ref{prop:filtration_properties_list}, imply that the chain homotopy is indeed filtered. We do this next.

Our homotopy of Floer data on the cylinder will be constant near the cylindrical ends, so 
\begin{equation}
\label{eq:geo_top_energy}
E_{top}(u)-E_{geo}(u)=-\int_S \tilde{u}^* d_S(H\beta)\geq -\int_S\max_{x\in\widehat{W}}d(H(x)\beta)
\end{equation}
where $d_S$ denotes the projection of $d$ onto $\Omega(S)$. We will construct a three-stage homotopy from the data of $\mathbb{C}\#_\lambda P$ to $\mathbb{C}$ such that the integral on the right is bounded by $\hbar/4$.
We first show that the lower bound in \eqref{eq:geo_top_energy} is $\geq -\hbar/4$ at the endpoints of the homotopy, and at two intermediate data.
\begin{enumerate}
    \item First, recall that the data on $P$ defining $\mu$ satisfies positivity of $\eqref{eq:geo_top_energy}$ by split-monotonicity. The datum on $\mathbb{C}$ defining $e^{H_-}$ satisfies $d(H(x)\beta)=d(f(s)^2H_{-,t}(x)dt)=(f^2)'(s)H_{-,t}(x)ds\wedge dt\leq -(\hbar/4)(f^2)'(s)ds\wedge dt$. Integrating over $\mathbb{C}$, we have \eqref{eq:geo_top_energy} $\geq -\hbar/4$. \label{item:htpy_glued}
    \item Secondly, let $K$ be an autonomous Hamiltonian such that $-\hbar/4\leq H_-\leq K\leq \frac{1}{2}H_+$. Consider the glued 1-form on $\beta$ on $\mathbb{C}\#_\lambda P$ as in \ref{item:htpy_glued}, and a Hamiltonian $H$ on $\mathbb{C}\#_\lambda P$ which agrees with $K$ outside the cylindrical ends, and in particular over $\mathbb{C}$, and is a monotone homotopy interpolating from $H_-$ to $K$, and $K$ to $\frac{1}{2}H_+$ in the input and output cylindrical ends, respectively. The restriction to $P$ is split-monotone so \eqref{eq:geo_top_energy} is positive. On restriction to $\mathbb{C}$, we have $d(H(x)\beta)=K(x)d\beta\leq -(\hbar/4)d\beta$ (noting that $d\beta\leq 0)$. Integrating over $\mathbb{C}$, we have \eqref{eq:geo_top_energy} $\geq -\hbar/4$. \label{item:htpy_int}
    \item Thirdly, consider a Hamiltonian $H$ as in \ref{item:htpy_int}, and a 1-form $g(s)dt$, where $g$ is some decreasing function interpolating between $1$ and $2$, and $(s,t)$ are cylindrical coordinates on $P\#_\lambda \mathbb{C}$. In the cylindrical ends, this datum has the form of a monotone homotopy and \eqref{eq:geo_top_energy}$\geq 0$. Outside the ends, we have \eqref{eq:geo_top_energy}$\geq -\hbar/4$ by the same argument as \ref{item:htpy_cyl}. \label{item:htpy_cyl}
    \item Finally, consider a monotone homotopy defining $\kappa$: the Hamiltonian is a monotone homotopy between $H_\pm$, and the 1-form is $dt$. \label{item:htpy_kappa}
\end{enumerate}
We proceed to construct the homotopy in three steps:
\begin{enumerate}
    \item First, we interpolate from the data \ref{item:htpy_glued} to \ref{item:htpy_int} by fixing the glued $1$-form $\beta$ on $C\#_\lambda P$, and linearly interpolating between the Hamiltonians. Because the lower bound in \eqref{eq:geo_top_energy} is concave in $H$ for fixed $\beta$, the condition \eqref{eq:geo_top_energy}$\geq -\hbar/4$ is preserved along this homotopy.
    \item Secondly, we interpolate from \ref{item:htpy_int} to \ref{item:htpy_cyl} by fixing the Hamiltonian $H$ and linearly interpolating between the 1-forms $\beta$ and $g(s)\mathrm{d}t$. Because the lower bound in \eqref{eq:geo_top_energy} is concave in $\beta$ for fixed $H$, the condition \eqref{eq:geo_top_energy}$\geq -\hbar/4$ is preserved along this homotopy.
    \item Finally, we view \ref{item:htpy_cyl} as (not-necessrily monotone) homotopy with Hamiltonian $Hg(s)$ and 1-form $dt$. We fix the 1-form, and linearly interpolate the Hamiltonian from $Hg(s)$ to the monotone homotopy defining $\kappa$. As in the first stage, the condition \eqref{eq:geo_top_energy}$\geq -\hbar/4$ is preserved.
\end{enumerate}

Finally, the parametrized count of this homotopy shows that $$\mu(e^{H_-},-)-\kappa=dh+hd$$where $\mathcal{A}_{H_+}(h(x))\geq \mathcal{A}_{H_-}(x)-\frac{\hbar}{4}$, where $h$ counts rigid Floer solutions to the 1-parameter family of Floer's equations arising from the homotopy of Floer data we have constructed.
\end{proof}
\subsection{Grading and Conley-Zehnder index}

One always has a mod $2$ Conley Zehnder index, which uses any trivialization $\Psi$ of $TW$ over an orbit $\gamma$: this exists, since symplectic vector bundles over a circle are trivial. Two trivializations differ by a multiple of $2$, so the mod $2$ index is well-defined and can be used to grade symplectic cohomology: $$\mathrm{deg}_{\mathrm{mod}\,2}(\gamma)\equiv n-\mu_{CZ}(\gamma, \Psi)\, (\mathrm{mod}\,2)$$

Instead, to define a (fractional) grading, we will require $c_1^{TW}$ to be $N$-torsion, which means we can fix some trivialization $\eta$  of $(\mathrm{det}\,TW)^{\otimes N}$ for some integer $N>0$. Given an orbit $\gamma$, consider the complex vector bundle $TW^{\oplus N}$ over $\gamma$, whose determinant is trivialized by $\eta$. Choosing a trivialization of $TW^{\oplus N}$ over $\gamma$ whose determinant agrees with $\eta$, one can compute the Conley-Zehnder index of the path of symplectic matrices $d\phi_t^{\oplus N}$ which we denote by $\mu_{CZ}^{TW^{\oplus N}}$. The cohomological grading of an orbit is defined to be 
$$\deg(\gamma)= n-\frac{1}{N}\mu_{CZ}^{TW^{\oplus N}}(\gamma)$$

Note that, as explained in \cite[Section 3.1]{diogo2018morsebottsplitsymplectichomology}, for null-homologous orbits $\gamma$ one can always use a capping surface $c$ of $\gamma$ to trivialize $TW|_c$ and compute this index as $\mu_{CZ}(\gamma, c)$. Moreover, since $c_1^{TW}$ is torsion, this does not depend on the chosen capping surface.

\begin{rem}
    The convention for grading symplectic homology differs from that of cohomology by $n-$.
\end{rem}

In our computations, we will begin by studying the indices of Reeb orbits which come in Morse-Bott families. After perturbing our autonomous Hamiltonian $H$ to the non-degenerate $\tilde{H}$ where the Hamiltonian orbits are isolated, the usual Conley-Zehnder index is modified. If $\gamma$ belongs to a Morse-Bott family $\mathcal{N}\subset \partial W$ and a Hamiltonian orbit $\gamma'$  is associated to a critical point $q$ of $f/S^1$ on a $\mathcal{N}/S^1$, together with the data of whether it is a minimum or maximum in the $S^1$ family of orbits living over $q$, then (by \cite[Lemma 2.4]{bourgeois_morse-bott_2003}, see also \cite[Section 3.3]{Oancea2004ASO} and \cite[Appendix C.5]{ritter_filtrations_2025}):
\begin{equation*}
\mathrm{deg}(\gamma')=\mathrm{deg}(\gamma)-\left(\mu_{Morse}(q)+\mu_{S^1}-\frac{1}{2}\dim_{\mathbb{R}}\mathcal{N}/S^1\right)    
\end{equation*}
The index $\mu_{S^1}\in \{0,1\}$ encodes whether the orbit is of type $\wedge$ resp. $\lor$.

The simple principal orbits come in a Morse-Bott family $\partial W$. We denote by $-2I$ the cohomological degree of the Hamiltonian orbit associated to a principal orbit and a critical point with $\mu_{Morse}(q)=\mathrm{dim}_\mathbb{R} \partial W/S^1, \mu_{S^1}=1$, in other words $\mu_{Morse}(q)+\mu_{S^1}=2n-1$.

\subsection{Principal orbits}
Our next task is to choose a suitable chain-level representative of the principal orbit, which will be the lowest action component of our periodic element $\mathfrak{M}$. We will consider $J$-shaped Hamiltonians (as in the definition of symplectic cohomology) to construct this class, as we can safely ignore the upper orbits in that case. Note that there are continuation maps between the J-shaped and V-shaped Hamiltonians which define a canonical map $\mathrm{SH}\rightarrow \widecheck{\mathrm{SH}}$ (as in \cite{CieliebakFrauenfelderOancea2010}), and moreover one can include the J-shaped Hamiltonians in the telescope model for V-shaped symplectic cohomology.

\begin{prop}
\label{prop:M_k_consistent}
    Consider a Hamiltonian $H$ which is a sufficiently small perturbation of a $b$-admissible autonomous V-shaped Hamiltonian as in Section \ref{sec:MB_perturbation}. Then $$H^{\bullet-2kI}(\mathrm{gr}^{-2kN} \mathrm{CF}_{(-\infty,b)}(H)) \cong H^\bullet(\partial W)$$ and this isomorphism commutes with continuation maps whenever the slope of $H$ at infinity is larger than $k$.
\end{prop}
\begin{proof}
Recall that the non-degenerate Hamiltonians $H$ are defined as small perturbations of an autonomous V-shaped Hamiltonian, supported in disjoint neighbourhoods of each family. First of all, by our definition of the discrete action filtration, the only orbits in the action window for $[a(-2kN),a(-2kN+1))$ are those which arise from perturbing the Morse-Bott family $\mathcal{N}$ of Reeb orbits of period $k$, which in fact have action in $[k-\epsilon,k)$. Provided the parameter $\varepsilon_{\mathcal{N}}$ (which controls the perturbation near $\mathcal{N}$) and the time-dependent perturbation are sufficiently small, the Floer cohomology of the perturbed Hamiltonian in this action window computes the local Floer cohomology of $\mathcal{N}$ by \cite[Lemma 2.8]{mclean2012localfloerhomologyinfinitely}. 

To ensure coherence under continuation maps, we can, possibly after subdividing the monotone homotopy, guarantee that continuation maps arise from isolated deformations of the perturbed family $\mathcal{N}$ in the sense of \cite[Lemma 2.5]{mclean2012localfloerhomologyinfinitely}, which induce isomorphisms on local Floer cohomology.
\end{proof}

\begin{rem}
    From now on, we will assume our Morse-Bott perturbations to be small enough so that Proposition \ref{prop:M_k_consistent} holds.
\end{rem}

\begin{defn}[Principal orbit on the chain level]\label{def:principalorbit}

Let $H$ be as in Definition \ref{def:profileautonomous}. We denote by 
\begin{equation}
    \check{M}_1^H\in \mathrm{gr}^{-2N}\mathrm{CF}^{-2I}_{(-\infty,b)}(H)
\end{equation}
the unique representative of $1\in H^0(\partial W)$ under the quasi-isomorphism in Proposition \ref{prop:M_k_consistent}. If $H$ is $J$-shaped, then one could take $b=\infty$.
\end{defn}

\begin{rem}
    If we assume that the induced function $f/S^1:\partial W/S^1\rightarrow \mathbb{R}$ has a unique local maximum (alternatively, $-f/S^1$ has a unique local minimum), then we may take $\check{M}_1^H$ to be the unique generator of type $\lor$ corresponding to that maximum.
\end{rem}

We now similarly give a characterisation of a choice of chain-level representative of the unit, in terms of our Morse-Bott model, up to a positive-energy correction term.

\begin{prop}[Multiplication by principal orbit is filtered]\label{prop:M1isfiltered} 
Let $H$ be a $J$-shaped Hamiltonian and $\check{M}_1^H\in \mathrm{CF}^{-2I}(H)$ denote the principal orbit. If $H_-$ is $b$-admissible and $H_+$ is $b-1$-admissible, then multiplication by $\check{M}_1^H$ descends to the truncated complexes, $$m=\mu(\check{M}_1^{H_-},\;\cdot\;):\mathrm{CF}^\bullet_{(-\infty, b)}(H_-, J_-)\rightarrow \mathrm{CF}^{\bullet-2I}_{(-\infty, b-1)}(H_+,J_+)$$
and $m$ respects the filtration $\mathcal{F}$ up to a shift; $$m(\mathcal{F}^p\mathrm{CF}^\bullet_{(-\infty, b)}(H_-, J_-))\subset\mathcal{F}^{p-2N}\mathrm{CF}^\bullet_{(-\infty, b-1)}(H_+, J_+).$$\end{prop}\begin{proof}First of all, a non-upper orbit with action $<b$ satisfies $\mathcal{A}_{H_+}(\check{M}_1 \cdot x)\geq \mathcal{A}_{H_-}(x)-1-\epsilon$ but since the lower orbits are contained in open balls of radius $<\epsilon$ around the period spectrum, the error $\epsilon$ can be discarded. This shows that the map respects the filtration. An upper orbit has action $\geq b+\epsilon$, so the output will have action $\geq b-1$.\end{proof}

\begin{prop}
\label{chain_level_unit}
Let $H$ be a Morse-Bott perturbation of a $b$-admissible, V-shaped Hamiltonian or a J-shaped Hamiltonian. Then the image of $e^H$ in $\mathrm{gr}^0\mathrm{CF}^0_{(-\infty,b)}(H)$ is a representative of $1\in H^0(\partial W)$ (respectively $1\in H^0(W)$) under Proposition \ref{prop:M_k_consistent}.
\end{prop}
\begin{proof}

We prove this for V-shaped Hamiltonians, the claim being almost identical for J-shaped ones. The degree zero part of the associated graded $\mathrm{gr}^0\mathrm{CF}^0_{(-\infty,b)}(H)\cong CM^0(\partial W)$ has a single closed element of degree zero, which is the one that represents the unit. Moreover, $e^H=\psi_\mathbb{C}(1)$ is closed and of degree zero and the PSS map admits constant solutions. This then forces the component of $\psi_\mathbb{C}(1)$ in $\mathrm{gr}^0\mathrm{CF}^0$ to be the representative of the unit.

\end{proof}

\begin{rem}
    Strictly speaking, in the above proposition if one works over $\mathbb{Q}$, one would obtain that $e^H$ and $1$ agree up to a scalar. However, one can just as well use integer coefficients for this proposition and it would not affect the proof, hence how we get the stronger conclusion without the need to rescale. The only reason we stick to field coefficients is to use the isomorphism between symplectic homology and symplectic cohomology later on, when we invoke the long exact sequence from \cite{CieliebakFrauenfelderOancea2010}.
\end{rem}

\subsection{Liouville structure on orbifold divisor complements and orbifold Boothby-Wang}

Given a smooth $W$ with $\partial W$ a Besse contact manifold, performing a symplectic cut provides a natural compactification $\mathcal{X}$ with $\mathcal{D}=\partial W/S^1$. Conversely, we will now describe, starting from a nice compactification $(\mathcal{X}, \mathcal{D})$ a Liouville structure on $W$. We also compute the gradings in the symplectic chain complex using the first Chern class of $\mathcal{X}$ as well as the age grading of the orbifold points. This is an extension of the setup of \cite{diogo2018morsebottsplitsymplectichomology}, \cite{DiogoLisi2019}: there, $\partial W$ is a prequantization bundle over a symplectic manifold but here we allow also for the base to be a symplectic orbifold. 

We briefly review the Liouville structure of an ample orbifold divisor complement. This follows the usual construction for smooth affine varieties with a single modification: we require $\mathcal{D}\subset \mathcal{X}$ to be a nice suborbifold so that it admits a symplectic tubular neighbourhood. These are called \textit{normalizable} suborbifolds in \cite[Proposition 17]{Munoz2020} - the condition of being normalizable rules out pathological examples. In particular, quasismooth divisors in weighted projective varieties are all normalizable, so all of our applications fit within this framework. 

First, we recall the orbifold Boothby-Wang theorem for circle orbibundles over symplectic orbifolds:
\begin{prop}{\cite[Theorem 7.1.6]{BoyerGalicki2008}, \cite[Theorem 1.3]{Kegel_Lange_2021}}
Let $\mathcal{D}$ be a symplectic orbifold with cyclic quotient orbifold singularities equipped with a symplectic form $\omega_\mathcal{D}$ such that $K\omega_\mathcal{D}=c_1(\mathcal{L})$ for an orbifold line bundle $\mathcal{L}$, for some $K>0$. Then the total space of the principal $S^1$-orbibundle $Y:=\mathbf{S}\mathcal{L}$ is an orbifold that can be equipped with a contact structure such that $d\eta=-K\pi^*\omega_\mathcal{D}$. The Reeb flow induces a circle action on $Y$ realizing it as a Seifert bundle over $\mathcal{D}$. Furthermore, if the local uniformizing groups of $\mathcal{D}$ act faithfully on the fibers of $Y$, then $Y$ is itself a manifold.
\end{prop}
The version which originally appeared in Boyer-Galicki's \cite[Theorem 7.1.6]{BoyerGalicki2008} is phrased in the language of almost Kahler manifolds and Sasakian geometry; the statement is discussed in the symplectic case in \cite[Theorem 1.3]{Kegel_Lange_2021}, where a cohomological condition is used instead of the faithful action to deduce that the contact orbifold is actually a manifold.

The next proposition is a straightforward generalization of \cite[Lemma 2.2]{DiogoLisi2019}.

  \begin{prop}
    Let $(\mathcal{X}, \omega_\mathcal{X})$ be a closed symplectic orbifold and $\mathcal{D}$ a smooth symplectic divisor whose class in homology represents $K[\omega_{\mathcal{X}}]$. Assume that the orbifold locus of $\mathcal{X}$ coincides with the orbifold locus of $\mathcal{D}$, ensuring that $\mathcal{X}\setminus \mathcal{D}$ is a smooth manifold. Then, the completion of $\mathcal{X}\setminus \mathcal{D}$ admits the structure of a Liouville manifold with Liouville 1-form $\lambda$ such that its contact boundary is equivalent to the Besse contact manifold which is the sphere orbibundle $\mathbf{S}\mathcal{N}_\mathcal{D}$ with its canonical Boothby-Wang contact structure. 
  \end{prop}

  \begin{proof}
    As mentioned, the proof is more or less identical to \cite[Lemma 2.2]{DiogoLisi2019}, with a few modifications needed to make it work for orbifolds: namely, we need to use a symplectic orbifold tubular neighbourhood theorem, as well as the orbifold Gysin sequence. Other than that, it is exactly identical. We sketch the main idea of the construction for completeness.

    \textit{Step 1}: We first equip the total space of the line orbibundle  minus the zero section with a canonical symplectic orbifold structure as follows:
    take the orbifold connection 1-form $\alpha\in \Omega^1(\mathcal{L}\setminus \mathcal{D})$. This defines a symplectic form \begin{equation*}
    \omega_\mathcal{L}:=-d\left(\frac{1}{K}e^{-\rho^2}\alpha\right)=e^{-\rho^2}\pi^*\omega_\mathcal{D}+\frac{2\rho e^{-\rho^2}}{K}d\rho\wedge \alpha.
    \end{equation*}

    which is exact on $\mathcal{L}^\#=\mathcal{L}\setminus 0$. This is the so-called canonical symplectic structure on $\mathcal{L}^\#$ induced by $(\mathcal{D},\omega_\mathcal{D})$.

    The circle bundle of radius $\rho_0$ of $\mathcal{L}$ is denoted by $Y$ and admits the contact form $\alpha_Y:=-\frac{1}{K}e^{-\rho_0^2}\alpha$. This implies that \begin{equation*}
        \mathcal{L}^\#\simeq ((-\infty, \rho^2)\times Y, d(e^r\alpha_Y)).
    \end{equation*}
    
    $Y$ is now a prequantization bundle over $\mathcal{D}$.

    \textit{Step 2}: Now, we use the symplectic orbifold tubular neighbourhood theorem \cite[Proposition 17]{Munoz2020}. This provides us with a symplectic embedding \begin{equation*}
        (\mathbf{D}\mathcal{N}_\mathcal{D}, \omega_{\mathcal{N}})\xrightarrow{\varphi} (\mathcal{X},\omega).
    \end{equation*}
    \textit{Step 3}: Finally, we need to compare the $1$-forms and contact structures.

    On $\mathcal{X}\setminus \mathcal{D}$, $\omega|_{\mathcal{X}\setminus \mathcal{D}}$ is exact with a primitive $\lambda_0$. Similarly, on $\mathcal{N}_\mathcal{D}\setminus \mathcal{D} = \mathcal{L}^\#$, we have that $\omega_\mathcal{L}=d(\frac{-1}{K}e^{-\rho^2}\alpha)=d\lambda_\mathcal{L}$. We would like to match $\lambda_\mathcal{L}$ and $\lambda_0$ under $\varphi$, but the only relation we have is that 
    \begin{equation*}
        d(\varphi^*\lambda_0 - \lambda_\mathcal{L})=0 \text{ on } (\mathbf{D}\mathcal{N}_\mathcal{D})^\#\simeq \mathbf{S}\mathcal{N}_\mathcal{D}.
    \end{equation*}

    The orbifold Gysin sequence reads
    \begin{equation*}
    0\rightarrow H^1_{S^1}(Y)\xrightarrow{\pi^*} H^1(Y)\rightarrow H^0_{S^1}(Y)\xrightarrow{\omega_\mathcal{D}}H^2_{S^1}(Y),
    \end{equation*}
    and shows that $\pi^*$ is an isomorphism. Therefore, we can write $$\varphi^*\lambda_0-\lambda_\mathcal{L}=dg+\pi^*\mu$$
     for some function $g$. Modify $\lambda_0$ to $\lambda:=\lambda_0-dg$ where by abuse of notation we denote $g$ the function induced by $g$ with support on $\varphi(\mathbf{D}\mathcal{N}_D)\subset \mathcal{X}$. We have that $\varphi^*\lambda=\lambda_\mathcal{L}+\pi^* \mu$. Finally, one constructs a flow $\phi^t$ defined via the vector field $\omega_{\mathcal{D}}$-dual to $\mu$, where $\omega_{\mathcal{D}}=\iota^*\omega_{\mathcal{X}}$. By a Moser argument exactly as in  \cite[Lemma 2.2]{DiogoLisi2019}, one obtains $(\phi^1)^*\varphi^*\lambda = \lambda_{\mathcal{L}}$, which concludes this sketch.
    
  \end{proof}
\part{Positive equivariant symplectic cohomology}\label{part:positiveeqsh}

  \section{A convenient formula for the Robbin-Salamon index and the grading}
  \label{sec:grading_formulae}
  In this section, we will describe a formula for the transverse Conley-Zehnder (also known as Robbin-Salamon) index of the fibers of a prequantization bundle over an orbifold $\mathcal{D}$ using the filling $W$. The benefit of this formula is that it depends only on some topological information, namely the Chern class of $\mathcal{X}$, as well as the (fractional) multiplicity of an orbit and the weights of the cyclic quotient singularity, if the orbit is exceptional.
  \begin{rem}
      Note that under monotonicity assumptions, that is, if the symplectic form is proportional to $c_1^\mathcal{X}$, then $c_1^W$ is torsion. This ensures that $\langle c_1^W, B\rangle=0$ for any compact surface $B$ in $W$, meaning that the index of a null-homologous orbit does not depend on the choice of capping surface in $W$. We will assume this from now on (otherwise, one should keep track of the capping).
  \end{rem}

\gradingformula
  \begin{proof}
      We present two proofs.\\ 
      \textit{Proof 1}: The first proof is based on comparing two different trivializations, similar to the formula in \cite[Section 3]{diogo2018morsebottsplitsymplectichomology} and \cite{bae2024rabinowitzfloerhomologyprequantization}, combining formulas from \cite{mak_orbifold_2026}.
      A capping surface $c$ inside $W$ of $\gamma$ induces a trivialization of $\gamma^*TW$ and produces the index $\mu(\gamma)=\mu(\gamma,c)$, which does not depend on $c$ provided that $ c_1^{TW}$ is torsion. On the other hand, such an orbit can be capped in the normal bundle $\mathcal{N}_\mathcal{D}$ using a cover of the normal fiber. This produces an orbifold capping $c_{orb}$ of $\gamma$. While orbifold vector bundles over an orbifold disk are not necessarily trivial, one can desingularize (\cite[Section 4]{chen_new_2004}) $c_{orb}^*T\mathcal{X}$ to get a smooth vector bundle $|T\mathcal{X}|$ over the smooth disk $|c_{orb}|$ underlying $c_{orb}$. By \cite[Section 5, Remarks 5.10-5.14]{mak_orbifold_2026}, 
  \begin{equation*}
      \mu(\gamma,c)=-\mu(\gamma, |c_{orb}|)+\underbrace{2c_1^\mathcal{X}(c\# c_{orb})-2\iota}_{2\langle c_1 |T_\mathcal{X}|, |B|\rangle},
  \end{equation*}
  where $\iota$ is the age grading at the orbifold point which $\gamma$ asymptotes to. This is essentially \cite[Proposition 4.2.1]{chen_new_2004}.
  
  It remains to compute $\mu(\gamma, |c_{orb}|)$, using the desingularized fiber frame. The local model is as follows: $\mathbb{D}\times \mathbb{C}^{n-1}/ \mu_m$, with the action being 
    \begin{equation*}
        \zeta \cdot (z_1,\dots, z_n)=(\zeta z_1, \zeta^{b_2}z_2,\dots, \zeta^{b_{n}}z_n).
    \end{equation*} The weights $b_i$ satisfy $0\leq b_i <m$ (some might be zero in case the orbifold locus is not isolated). The Reeb flow rotates the first factor
    \begin{equation*}
        [z_1,...,z_n]\mapsto [e^{2\pi it}z_1,...,z_n]
    \end{equation*}
    and the exceptional orbit is given by $[e^{2\pi i t},0,...,0]$ with $ t\in[0,1/m]$  in the chart $\mathbb{D}\times \mathbb{C}^{n-1}$ 

    Consider the tangent bundle of $\mathbb{C}^n/\mu_m$ and restrict it to the $\ell$-fold cover of the orbifold disk fiber. The desingularization process of Chen and Ruan first restricts to the punctured disk ${\mathbb{D}}^\times$. Afterwards, a coordinate change is applied: 
    \begin{equation*}
     [z,v_1,\dots, v_n]\in (\mathbb{D}\times \mathbb{C}^{n})/\mu_m\mapsto (z^m,z^{-\overline{\ell}}v_1,z^{-\overline{b_2\ell}}v_2,\dots, z^{-\overline{b_n\ell}}v_n)\in \mathbb{D}\times \mathbb{C}^n   
    \end{equation*} 
    which extends over the origin, trivializing the desingularized tangent bundle over the smooth underlying disk. Here $z=z_1^{1/\ell}$ is a coordinate on the $\ell$-fold cover of the normal fibre. 

    In this coordinate frame, the linearized Reeb flow looks like\footnote{Note that this is exactly the same formula as in \cite[Proposition 3.4]{LiZhou2026}, except we allow for multiple covers}
    \begin{equation*}\label{eq:linearizedReeb}
    (u,u_1,\dots, u_n)\mapsto (e^{2\pi imt/\ell}u, e^{2\pi i(\ell-\overline{\ell})t/\ell}u_1, e^{-2\pi i(\overline{b_2\ell})t/\ell}u_2, \dots, e^{-2\pi i(\overline{b_n\ell})t/\ell}u_n).
    \end{equation*}
     The Conley-Zehnder index of this path of symplectic matrices (for $t\in [0,\ell/m]$) is 
    \begin{equation*}
    \begin{split}
        2\left(\frac{\ell-\overline{\ell}}{m}\right)+\sum_{\overline{lb_i}\neq 0, i\geq 2}^n\left(1+2\left\lfloor-\frac{\overline{b_i\ell}}{m}\right \rfloor\right)=2\left\lfloor\frac{\ell}{m}\right\rfloor-(n-1-\dim_\mathbb{C} \mathcal{N}/S^1)
    \end{split}
    \end{equation*}
    since the number of $\overline{lb_i}$ which are zero corresponds to the complex dimension of the relevant orbifold locus (it is zero when it is isolated). This completes the first proof, since $\lfloor \frac{l}{m}\rfloor=\lfloor B \cdot \mathcal{D}\rfloor$.

    \noindent\textit{Proof 2}: A second approach is given by considering the expected dimension of the moduli space of unparametrized holomorphic planes capping it and keeping it fixed. This is given as the SFT index (see Bourgeois' \cite{bourgeois_morse-bott_2003}, or \cite{bae2024rabinowitzfloerhomologyprequantization} Remark 2.4) 
    \begin{equation*}
        |\gamma|_{SFT}=\frac{1}{2}\dim W -3 + \mu_{CZ}(\gamma)-\frac{1}{2}\dim \mathcal{N}/S^1
    \end{equation*}
    However, since the Reeb orbits in our setup are multiple covers of circle fibers in the normal bundle of a divisor $\mathcal{N}_{\mathcal{D}}$, any holomorphic cap of $\gamma$ in $W$ can be extended to an orbifold sphere $B$ with some prescribed tangency with the divisor at infinity. In particular, generic orbits correspond to smooth spheres, but a cap of an exceptional orbit will compactify to an orbifold sphere.

    The expected dimension of the moduli space of holomorphic spheres in $\mathcal{X}$ hitting the divisor $\mathcal{D}$ at a prescribed orbifold point with some fractional multiplicity is given by \begin{equation*}
    2(c_1^\mathcal{X}|B|-|B|\cdot \mathcal{D}+\dim_{\mathbb{C}}X+1-3)=2(c_1^\mathcal{X}B -\iota - \lfloor B\cdot \mathcal{D}\rfloor+\dim_\mathbb{C}\mathcal{X}-2)
    \end{equation*}
    where $|B|$ is the normalization of $B$ and $\iota$ is an age grading number (coming from Chen-Ruan cohomology). This is the same as  $|\gamma|_{SFT}+\dim_\mathbb{R}\mathcal{N}/S^1$ (since in $|\gamma|_{SFT}$, the planes keep the orbit fixed, but in the above moduli problem the orbit can move in its Morse-Bott space),
  
    From this we can read off that \begin{equation*}
        \mu_{CZ}(\gamma)=2(c_1^\mathcal{X}B -\iota - \lfloor B\cdot \mathcal{D}\rfloor)+(\dim_\mathbb{C}W-\dim_\mathbb{C} \mathcal{N}/S^1 -1).
    \end{equation*}
  \end{proof}

  \begin{cor}\label{cor: principalorbitdegree} Suppose $c_1^\mathcal{X}= \lambda \mathcal{D}$ on $H_2(\mathcal{X};\mathbb{Q})$. Then, the $i$-th cover of a Hamiltonian orbit corresponding to a principal orbit which is a maximum in its Morse-Bott family has degree $-2i(\lambda-1)$ in symplectic resp. V-shaped symplectic cohomology.
  \end{cor}
  \begin{proof}
      This is clear for the covers of the principal orbit which are null-homologous by the above formulas. For the ones which are not null-homologous, one argues exactly as in \cite[Section 3.1]{DiogoLisi2019}: one shows that the McLean-Seidel fractional grading and the above grading coincide for covers of the principal orbit which are null-homologous. Moreover, one shows that for each $i$, the difference between the two is a linear function in $i$. This implies that they must agree for all $i$.
  \end{proof}
  \begin{rem}
      One could get away with the weaker condition that $c_1^\mathcal{X}=\lambda \mathcal{D}$ on the spherical part of $H_2(\mathcal{X};\mathbb{Q})$ if $\pi_1(W)\simeq H_1(W;\mathbb{Z})$, since then null-homologous orbits are contractible and the capping surface could be taken to be a disk.
  \end{rem}
  \begin{rem}
      The proof of the above Corollary does not carry over immediately to the exceptional orbit case. To show that the difference between the McLean-Seidel fractional grading and the fractional cap grading is linear in $i$, Diogo and Lisi \cite{DiogoLisi2019} use the constant trivialization of $T_{\pi(\gamma)} \mathcal{D}$ and lifting it to $Y$. However, for the exceptional orbits, there is no 'constant trivialization' coming from $T_{\pi(\gamma)} \mathcal{D}$ since those orbits project to an orbifold point, where the tangent orbibundle carries non-trivial stacky structure.
  \end{rem}

  We reformulate and simplify the formula in the case of an isolated cyclic quotient singularity.
  \begin{cor}\label{cor:explicitcyclicexample}
      Suppose there is an isolated $\frac{1}{w}(1,b_2,\dots, b_n)$ orbifold point on $\mathcal{X}$ (with weight $1$ corresponding to the normal direction of $\mathcal{D}$) and let $\delta_{i/w}$ denote the exceptional Reeb orbit which is the $i$-th cover of the exceptional orbit over that point (with $w \nmid i$) . If $\delta_{i/w}$ is null-homologous then \begin{equation*}
          \mu_{CZ}(\delta_{i/w})=2\left(c_1^\mathcal{X}B-\frac{\overline{i}+\overline{ib_2}+\dots+\overline{ib_n}}{w}-\left\lfloor \frac{i}{w}\right\rfloor\right)+n-1.
      \end{equation*}
      When $\langle c_1^\mathcal{X},-\rangle=\lambda \langle \mathcal{D}, -\rangle$ on $H_2(\mathcal{X};\mathbb{Q})$, this is equal to \begin{equation*}
      \begin{split}
          \mu_{CZ}(\delta_{i/w})&=2\left(\lambda \frac{i}{w}-\frac{\overline{i}+\overline{ib_2}+\dots+\overline{ib_n}}{w}-\left\lfloor \frac{i}{w}\right\rfloor\right)+n-1
          \\&=\frac{2}{w}\left((\lambda-1)i-\left(\overline{ib_2}+\dots+\overline{ib_n}\right)\right)+n-1.
          \end{split}
      \end{equation*}
  \end{cor}
  \begin{rem}\label{rem:extensionofgradingformula}
      The last formula applies not only to the nullhomologous exceptional orbits, but to any exceptional orbit which is equivalent to a multiple of the principal orbit in $H_1(W)$. That is, if $g^n\cdot \delta$ is null-homologous, one can use the above formula to compute the index (here $g$ denotes the loop of diffeomorphisms generated by the Reeb flow). On the other hand, the index of $\delta $ differs from that of $g^n\cdot\delta$ by $-2nI$. So the index of $\delta$ can be computed by the same formula.
  \end{rem}

  \begin{cor}
      In the mod $2$ grading of symplectic resp. V-shaped symplectic cohomology, a Hamiltonian orbit $\check{\delta}_{l/m}$ of type $\lor$ corresponding to an isolated cyclic quotient singularity on $\mathcal{D}$ contributes in degree $0\mod 2$.
  \end{cor}
  \begin{proof}
      The mod 2 index does not depend on the choice of trivialization, so we can just choose the desingularized normal fiber framing for the Reeb flow and then perform a Morse-Bott perturbation. The difference between the Reeb mod $2$ index and the mod 2 index of a Hamiltonian orbit $\check{\delta}_{l/m}$ of type $\lor$ is $1$, since the Morse-Bott family of orbits in this case is just a circle, which is then broken into two orbits using a small time-dependent perturbation.
      
      The formula for the mod $2$ index, after perturbation, becomes $-2\lfloor \frac{l}{m}\rfloor +n$ and the mod 2 cohomological grading is given by subtracting this from $n$, i.e. it is $2\lfloor \frac{l}{m}\rfloor=0\mod 2$.
  \end{proof}

\section{Lacunarity of the positive equivariant symplectic cohomology}\label{sec:positiveequivariant}
In this section, we show that under some assumptions on $\mathcal{D}$, the Morse-Bott spectral sequence computing positive equivariant symplectic cohomology degenerates on the first page.

\degeneration*
\begin{proof}
The first page of the Morse-Bott spectral sequence has contributions associated to each Morse-Bott family of orbits. In the case of the Seifert bundle $\mathbf{S}\mathcal{N}_D$, each orbit of period $\leq 1$ comes in a family parametrized by a component of the inertia orbifold $I\mathcal{D}$. The simple principal orbits are parametrized by $\mathcal{D}$ itself (and so are their multiple covers), whereas exceptional orbits correspond to the various orbifold loci of $\mathcal{D}$.

For each $N$ and each connected component $S\subset I\mathcal{D}$ parametrizing a family of orbits, there is a contribution to the first page of the spectral sequence given by the cohomology shifted by a shifting number: 
\begin{equation*}
    H^\bullet(S)\left[-\left(\mu_{CZ}(S)-(n-1)+\frac{1}{2}\dim_\mathbb{R} S\right)+1\right]
\end{equation*}

We first show that the spectral sequence degenerates on the $E_1$ page by studying the mod $2$ CZ index. The Conley-Zehnder index modulo $2$ does not depend on the chosen trivialization so we can use the normal fiber capping disk to compute it. Recall that 
\begin{equation*}
    \mu_{CZ}(\delta_{l/m}, c_{orb})=2\left\lfloor \frac{l}{m}\right\rfloor -(n-1)+\frac{1}{2}\dim_\mathbb{R} \mathcal{N}/S^1
\end{equation*}
In turn, the shifting number is equivalent to $1 \,(\mathrm{mod}\,2)$. By assumption, all strata of the inertia orbifold have only even cohomology, hence on the $E_1$ page all contributions are concentrated in degree $1 \,(\mathrm{mod}\,2)$. Therefore, there is no differential and the spectral sequence degenerates there.

When every orbit is nullhomologous, one has a $\mathbb{Z}$-grading and the group $A$ corresponds to the contributions of the orbits of period $\leq 1$. Each orbit of period $>1$ is equivalent to $g^m \gamma$ for some $\gamma$ with period $\leq 1$ and some $m\geq 1$. Since rotating by $g$ shifts the degree down by $-2I$, the claim follows.

\end{proof}

\begin{rem}
    In \cite{kwon_brieskorn_2016}, each component $S$ contributes $H_{p+q-(\mu_{CZ}-\frac{1}{2}\dim_\mathbb{R} S)}(S)$ to $E^1_{p,q}$ for positive equivariant symplectic homology. We work with cohomology instead, so we need to apply a suitable shift by $n-$. In particular, $H_i(S)$ contributes whenever $p+q-(\mu_{CZ}-\frac{1}{2}\dim_\mathbb{R} S)=i$. This, in turn, happens exactly when $\big(n-(p+q)\big)+\big(\mu_{CZ}-(n-1)+\frac{1}{2}\dim_\mathbb{R} S)-1=\dim_\mathbb{R} S-i $ so that the contribution to $E^{p,q}_1$ is given by $H^{(n-p-q)+(\mu_{CZ}-(n-1)+\frac{1}{2}\dim_\mathbb{R} S))-1}(S)$ which explains our shifting number.
\end{rem}

\begin{cor}
    If $(\mathcal{X},\mathcal{D})$ is a monotone compactification of a simply-connected $W$ with $\langle c_1^\mathcal{X},-\rangle=\lambda \langle \mathcal{D}, -\rangle$ on $H_2(\mathcal{X};\mathbb{Q})$ and $\lambda \in \mathbb{N}$ with $\lambda >1$, then the mean Euler characteristic is equal to
    \begin{equation*}
        \chi_m(W)=\frac{-1}{2(\lambda-1)} \chi(I\mathcal{D})
    \end{equation*}
\end{cor}
\begin{proof}
    The mean Euler characteristic does not change under passage to the first page of the spectral sequence. The proof then is a straightforward application of the index formulas from the previous section.
\end{proof}
\section{Application: orbifold compactifications of $\mathbb{C}^n$}
\label{sec:orb_div_comp}
The aim of this section is to prove the following theorem:
\orbcompclassification*

\subsection{Pinning down the Fano index and singularity type using symplectic cohomology}
The first task in our proof is to compute the Fano index using $n$ and the isotropy orders $w_i$ and moreover single out the type of singularity that can appear on $\mathcal{X}$ when $d=1$.

\begin{rem}
    By \cite[Lemma 5.18]{mclean2012growth}, different projective compactifications of an affine variety like $\mathbb{C}^n$ produce deformation equivalent Liouville structures, hence the symplectic cohomology (and its variants) remains the same. The setup in McLean's paper uses snc compactifications where all of the components are smooth. In our setup, we have an orbifold compactification: after resolving the cyclic quotient singularities by a blowup $Y\rightarrow X$, one obtains an snc compactification. We briefly explain why this does not affect the Liouville structure: for step 1 in \cite[Lemma 5.18]{mclean2012growth}, one just notes that each point on an orbifold has an analytic chart given by $\mathbb{C}^n$ modulo a finite cyclic group action, which makes the argument go through. For Step 2, however, it is not always possible to pull back $\mathcal{O}(\mathcal{D})$ to $Y$ since it is only a reflexive sheaf on the coarse space $X$. It is nevertheless $\mathbb{Q}$-Cartier, so after passing to some sufficient power making it a line bundle, the same argument as in \cite[Lemma 5.18]{mclean2012growth} applies. Note that taking an $N$-th power of $\mathcal{O}(\mathcal{D})$ and in turn taking a power of the section $s$ of $\mathcal{O}(\mathcal{D})$ which one uses to define the Liouville structure on $\mathcal{X}\setminus \mathcal{D}$ only multiplies the psh function $h=-\log ||s||^2$ by $N$.
    
    We will utilize the model of symplectic cohomology induced by a hypothetical compactification of $\mathbb{C}^n$ by an orbifold divisor to gain information, by comparing it to the standard computation of the symplectic (resp. positive equivariant symplectic) cohomology of $\mathbb{C}^n$.
\end{rem}

\begin{prop}\label{prop:Fanoindex}
    Suppose $\mathcal{X},\mathcal{D}$ is an orbifold compactification of $\mathbb{C}^n$ with isolated cyclic quotient singularities with isotropy orders $w_1,\dots, w_d$. Then, the Fano index of $\mathcal{X}$ is $n+1+\sum_{j=1}^d (w_j-1)$. When $d=1$, the singularity must be of type $\frac{1}{k}(1,...,1)$ for some positive integer $k$.
\end{prop}
\begin{proof}
On one hand we know that $\mathrm{SH}^\bullet_{+,S^1}(\mathbb{C}^n)$ is rank 1 in all negative odd degrees. Shifting this down by one for simplicity, the group $\mathrm{SH}^\bullet_{+,S^1}(\mathbb{C}^n)[-1]$ is rank $1$ in all negative even degrees.

On the other hand: if $\mathcal{X}\setminus \mathcal{D}=\mathbb{C}^n$, then after shifting down by one for simplicity again, by the formulae for the index in Section \ref{sec:grading_formulae}, the first page of the Morse-Bott spectral sequence computing $\mathrm{SH}^\bullet_{+,S^1}(\mathbb{C}^n)$ has contributions from the simple principal orbits in degrees
\begin{equation*}
-2r, -2r+2, \dots, -2r+2n-2
\end{equation*}
since $H^\bullet(D;\mathbb{Q})\simeq H^\bullet(\mathbb{P}^{n-1};\mathbb{Q})$. This is because $-K_\mathcal{X}=(r+1)\mathcal{L}$ and $\mathcal{L}=\mathcal{O}(\mathcal{D})$. This spectral sequence degenerates at the first page since $H^\bullet(I\mathcal{D};\mathbb{Q})$ is concentrated in even degrees. Thus, we get that $r\geq n$ since the contribution in degree $-2r+2n-2$ has to be negative.

The $i$-th cover of an exceptional orbit corresponding to an isolated orbifold point of type $\frac{1}{w}(1,b_2,\dots, b_n)$ contributes in degree (again after shifting down by one)
\begin{equation*}
-2\left((r+1)\frac{i}{w}-\frac{\overline{i}+\overline{b_2i}+\dots+\overline{b_ni}}{w}-\left\lfloor \frac{i}{w}\right\rfloor\right)=-2\left(r\frac{i}{w}-\frac{\sum_{j=2}^n \overline{b_ji}}{w}\right)
\end{equation*}

Each of these has to be negative, since the spectral sequence degenerates at this page and $\mathrm{SH}^\bullet_{+,S^1}(\mathbb{C}^n)[-1]$ is concentrated in negative even degrees. In particular, if $i>w$, then this will become less than $-2r$. This is because $\delta_{i/w}$ and $\delta_{(w+i)/w}$ contribute in degrees that differ by $-2r$. Hence, the contributions in degrees $-2,-4,\dots, -2r+2n$ must come from orbits of period less than one. Conversely, if some orbit of period less than one (i.e. $i<w$) contributed to a degree $-2A<-2r$, then in that degree, there will be more than one contribution: another will come from an orbit of period bigger than one, which comes from a sufficient rotation of the orbit that contributes in degree $- (\overline{2A} \,\mathrm{mod}\,2r)$. This contradicts the fact that $\mathrm{SH}^\bullet_{+,S^1}(\mathbb{C}^n)[-1]$ is rank one in all negative even degrees.
Therefore, we conclude that the contributions of the orbits of period less than one, in other words the contributions of the orbits $\delta_{i/w_j}, 1\leq i <w_j$, account for exactly the degrees $-2,-4,\dots, -2r+2n$. This immediately implies that 
\begin{equation*}
\sum_{j=1}^d (w_j-1)=r-n,
\end{equation*}
i.e. that $r+1=n+1+\sum_{j=1}^d (w_j-1)$.

In the case when $d=1$, we set $w_1=k$. The formula in that case is $r=n+k-1$. On the other hand, we know that $\frac{-2}{k}(r-b_2-\dots-b_n)$ is a negative integer, so that $r\geq k+b_2+\dots +b_n$, Therefore, since $r=n+k-1$, this implies that $n-1\geq b_2+\dots+b_n$. But the $b_i$ are positive integers, forcing them to all be equal to $1$ in this case.

\end{proof}

\begin{rem}
    Note that the contribution of $\delta_{i/w}$ (after shifting down by one as before) is \begin{equation*}
        -2\left(r\frac{i}{w}-\frac{\sum_{j=2}^n \overline{b_ji}}{w}\right)=-2i\left(r\frac{1}{w}-\frac{\sum_{j=2}^n \overline{b_j}}{w}\right)-2\sum_{j=2}^n \left\lfloor \frac{ib_j}{w}\right\rfloor.
    \end{equation*} 
    In other words, it is $i$ times the degree that $\delta_{1/w}$ contributes to, minus an error term bounded by $2i(n-1)$.
\end{rem}

\begin{rem}
    Our shifting down by $-1$ is to align with the cohomological gradings in symplectic cohomology. For example, the degree of $\check{M}_1$ in $\mathrm{CF}(H)$ (for a Hamiltonian with sufficiently high slope) will be $\deg \check{M}_1=-2(n+\sum (w_j-1))$, whereas the other ones corresponding to the family of principal orbits will be shifted upwards by some Morse index bounded by $2n-2$. Similarly, the exceptional orbits satisfy $\deg \check{\delta}_{i/w}=-2\left(r\frac{i}{w}-\frac{\sum_{j=2}^n \overline{b_ji}}{w}\right)$ which is necessarily $<0$. In the case $d=1$ we see that $\deg \check{\delta}_{1/k}=-2$ and $\deg \hat{\delta}_{1/k}=-1$, and these are the only orbits contributing in these degrees.
\end{rem}

\begin{prop}
    In the same setup as before, denote $w_{max}:=\max(w_i)_{1}^d$. Then, there is a shortest Reeb orbit $\delta_{1/w_{\max}}$ whose Conley-Zehnder index is $n+1$ and it is the unique Reeb orbit with this property.
\end{prop}
\begin{proof}
    We will use the vanishing $\mathrm{SH}^\bullet(\mathbb{C}^n)=0$ and recompute it via the model of symplectic cohomology which is suited to a hypothetical orbifold divisor compactification.

    Choose a directed system of autonomous Hamiltonians linear at infinity (using the Liouville structure on $\mathbb{C}^n$ induced from the compactification) successively approximating a quadratic Hamiltonian. Moreover, choose them so that on the interior of $\mathbb{C}^n$, they are $C^2$ close to a Morse function with a single local minimum. Perform a small Morse-Bott perturbation to each of these Hamiltonians, as in Section \ref{sec:Floer_theory} in such a way so that $H_i\equiv H_{i+1}$ on the part where $H_i$ is not linear. Following \cite[Sections 8.18]{ritter_topological_2013}, we can ensure that the continuation maps are inclusions of subcomplexes. In particular, we conclude that one can use a Hamiltonian of sufficiently high slope to compute $\mathrm{SH}^{-1}$ and $\mathrm{SH}^0$. This is because the cohomological degree of the principal orbit is $\deg \check{M}_1=-2(n+\sum_1^d (w_j-1))$ by Corollary \ref{cor: principalorbitdegree} and Proposition \ref{prop:Fanoindex}, so that all the new orbits after a sufficiently high slope will have degree $\ll 0$.

    In particular, we can take the slope to be slightly larger than $1$, which accounts for exactly the orbits in degrees $-2r, -2r+1, \dots, -1, 0$. The exceptional orbits occupy the degrees $-2r+2n, -2r+2n+1, \dots, -2,-1$, with exactly one orbit contributing in each of these degrees and with the orbits of type $\lor$ resp. $\wedge$ being in even resp. odd, degrees.

    Since all of these orbits are good, $$d\check{\gamma}=\langle d \check{\gamma}, \hat{\gamma}\rangle \hat{\gamma}=0$$
    Hence, all orbits of type $\wedge$ and degree less than $-1$ necessarily must have $d\hat{\gamma}\neq 0$ (otherwise, the symplectic cohomology would not vanish). In particular, for all of these orbits there is a Floer cylinder connecting them to another non-constant orbit. However, the orbit of smallest period $\hat{\delta}_{1/w_{\max}}$ (hence, highest action) cannot be connected to another non-constant orbit, and therefore is forced to have degree $-1$. In particular, there is at most one exceptional orbit of period $\frac{1}{w_{max}}$. The underlying Reeb orbit $\delta_{1/w_{\max}}$ has Conley-Zehnder index $n+1$, finalising the claim.
 
\end{proof}

\subsection{Using symplectic cohomology to construct a uniruling}
In this section, we will show that $\mathcal{X}$ is uniruled by certain rational curves when the $w_j$ are distinct. This will be done by an application of the vanishing of symplectic cohomology, as well as a particular model suited to the realization of $\mathbb{C}^n$ as a divisor complement. 

In order to do this, we will exhibit the existence of Floer cylinders connecting an orbit corresponding to an isolated cyclic quotient singularity to a constant orbit in the interior. We show that the moduli space of such Floer planes is nonempty for all admissible Hamiltonians and, after an application of a neck-stretching argument, we conclude that an honest $J$-holomorphic plane asymptotic to a Reeb orbit exists. Such holomorphic planes (by Gromov removal of singularities, as in \cite[Lemma 2.6]{DiogoLisi2019}) compactify to and biject with holomorphic spheres with a tangency with the divisor $\mathcal{D}$. The tangency is entirely determined by the Reeb orbit: in our setup, it will be a simple Reeb orbit corresponding to a cyclic quotient singularity of isotropy order $w_{\max}$, so the contact order is the fraction $\frac{1}{w_{\max}}$.

\subsubsection{Neck-stretching setup and almost complex structures}

We follow the account of neck-stretching in the Hamiltonian Floer setup as in \cite{pascaleff2019symplectic}, which itself is based on the work of Bourgeois and Oancea \cite{bourgeois_exact_2009}. 

Begin with a Hamiltonian $H$ which is $0$ on the interior and linear, with some finite slope, at infinity. We construct a one-parameter deformation of $(W_R,H_R,J_R)$ depending on a stretching parameter $R$ which simultaneously:

\begin{itemize}
    \item stretches the neck near a contact hypersurface which is a small pushoff of $\partial W$ under the contracting Liouville flow,
    \item degenerates the time-dependent $H_R$ to the autonomous $H_\infty=H$ which is zero on the interior,
    \item degenerates the time-dependent $J_R$ to a time-independent complex structure $J_W$ which is cylindrical at the end and the stretching region.
\end{itemize}

As in \cite[Lemma 2.7]{DiogoLisi2019}, we take the $J_s$ so that $(W_R,J_R)$ converge in the SFT sense to a split almost-complex manifold whose bottom level is $(W,J_W)$. 

The almost complex structure $J_W$ is obtained by first taking an admissible almost complex structure $J_\mathcal{X}$ on $\mathcal{X}$ as in \cite[Definition 2.6]{diogo2018morsebottsplitsymplectichomology} (adapted suitably to the orbifold setting). Then, one obtains an almost complex structure $J_W$ on $W$ which agrees with the restriction of $J_\mathcal{X}$ on a compact subset and is cylindrical outside of a slightly larger compact subset. Conversely, one might also begin with a $J_W$ that is Reeb-invariant and cylindrical at infinity, which extends to an almost complex structure $J_\mathcal{X}$.

Under this setup, finite energy $J_W$ planes asymptotic to a Reeb orbit compactify, by Gromov removal of singularities, to holomorphic spheres with a prescribed tangency with the divisor which depends on the Reeb orbit (see \cite[Lemma 2.6]{DiogoLisi2019}).

\begin{rem}
    We remark that it is not always the case that the integrable complex structure $J^{int}_\mathcal{X}$ can be obtained from this process: nevertheless, the complement $(\mathcal{X}\setminus \mathcal{D}, J^{int}_{\mathcal{X}}|_{\mathcal{X}\setminus \mathcal{D}})$, where $J^{int}_{\mathcal{X}}$ is the integrable complex structure, is \textit{asymptotically cylindrical} (see \cite[Section 3.2]{Bourgeois_Eliashberg_Hofer_Wysocki_Zehnder_2003}), for which a version of the compactness theorem also holds, as proved in \cite[Section 4]{Bao2015JHolomorphic}.
\end{rem}

This limiting family is chosen so that one can invoke the SFT compactness theorem \cite{Bourgeois_Eliashberg_Hofer_Wysocki_Zehnder_2003}. Moreover, the Hamiltonian orbits for each $H_R$ are naturally identified, for any $R$. In particular, if a moduli space $\mathcal{M}(\gamma_-,\gamma_+;H_R, J_R)$ of Floer cylinders contributing to the differential is nonempty for a sequence of $R$ growing arbitrarily large, then there must exist a subsequence of Floer cylinders converging to a generalized holomorphic building with:\begin{itemize}
       \item  a bottom level consisting of a $J_W$ holomorphic curve
        \item intermediate levels being holomorphic curves in the symplectisation
        \item a top level which is a Floer solution for the radial Hamiltonian and cylindrical complex structure in the symplectisation
    \end{itemize}
\subsubsection{Uniruling by planes via symplectic cohomology vanishing}
\begin{prop}\label{prop:uniruling}
    Assume $(\mathcal{X},\mathcal{D})$ is a compactification of $\mathbb{C}^n$ with isolated cyclic quotient singularities. Then: through every point $x$ in $\mathcal{X}\setminus \mathcal{D}$ there is a finite energy $J_W$-holomorphic plane $L_x$ passing through it and asymptotic at its end to the exceptional Reeb orbit corresponding to the $\frac{1}{w_{\max}}$ point on $\mathcal{D}$. This plane compactifies to a $J_\mathcal{X}$-holomorphic sphere $C_x$ which satisfies $C_x \cdot \mathcal{D}=\frac{1}{w_{\max}}$. Moreover, the class $[C_x]\in H_2(\mathcal{X};\mathbb{Q})$ is independent of $x$. In addition, the counts of rigid $J_w$ planes $L_x$ corresponds to a non-vanishing Gromov-Witten invariant $\langle [\mathrm{pt}_x], \mathbf{1}_{1/w_{\max}} \rangle_{[C]}^\mathcal{X}$.
\end{prop}

\begin{proof}

    We argue in a similar way as in the previous proposition, using the vanishing of symplectic cohomology. Take a Hamiltonian of slope slightly higher than $\frac{1}{w_{\max}}$: these already compute $\mathrm{HF}^{-1}=\mathrm{SH}^{-1}=0$ and $\mathrm{HF}^{0}=\mathrm{SH}^{0}=0$. Moreover, the only orbits that contribute are the constant orbit in degree $0$, as well as $\hat{\delta}_{1/w_{\max}}$ in degree $-1$ and $\check{\delta}_{1/w_{\max}}$ in degree $-2$.

    We now restrict to Hamiltonians of this fixed slope and perform the neck-stretching procedure from the previous section by taking a sequence of admissible Hamiltonians $H_s$. 
    We choose this sequence of $H_R$ so that in the interior they are equal to $\epsilon_R g$, where $\epsilon_R\rightarrow 0$ and $g$ is a Morse function, so that the critical point $e$ stays the same as one varies $R$. Since $\mathrm{HF}^0(H_R)=\mathrm{SH}^0(\mathbb{C}^n)=0$, restricting for the time being to $\mathbb{F}_2$ coefficients to not worry about orientation conventions, we deduce that $d\hat{\delta}_{1/w_{\max}}=e$ and the mod $2$ count of Floer trajectories is nonzero.

    Let $$\mathcal{P}(\delta_{1/w_{\max}}, J_W):=\{u:\mathbb{C}\rightarrow \hat{W}|\,\, \overline{\partial}_Ju=0, \, u\text{ asymptotes to }\delta_{1/w_{\max}} \}$$ denote the space of $J_W$ planes asymptotic to $\delta_{1/w_{\max}}$: it has virtual dimension $2n+2$. Imposing a point constraint at $x$ reduces this to $2$. Moreover, since the orbit $\delta_{1/w_{\max}}$ is simple, a generic $J_W$ suffices to regularize this moduli space. Let moreover $\mathcal{P}_x:=\mathrm{ev}_0^{-1}(x)$ and $\mathcal{C}$ denote the space of Floer trajectories from $\hat{\delta}_{1/w_{\max}}$ to $\delta_{1/w_{\max}}$, which is diffeomorphic to $\mathbb{R}$.

    The SFT compactness theorem ensures that for a sequence of Hamiltonians $H_R$ as before, there is a subsequence of Floer cylinders in $\mathcal{M}(e, \hat{\delta}_{1/w_{\max}}; H_R, J_R)$ which converge to a generalized holomorphic building. Since the Reeb orbit $\delta_{1/w_{\max}}$ is simple and of minimal period, the limiting punctured curve in the symplectization cannot have more than one output. Hence, the top level in the resulting building will consist of a Floer solution connecting the Hamiltonian orbit $\hat{\delta}_{1/w_{\max}}$ to the underlying Reeb orbit $\delta_{1/w_{\max}}$. There will be moreover some trivial cylinders in the symplectization which are all asymptotic to $\delta_{1/w_{\max}}$, together with a bottom level which is a $J_W$-holomorphic plane. This $J_W$-holomorphic plane corresponds to a finite energy once punctured $J_\mathcal{X}$-holomorphic sphere $L_x$, which by Gromov's removal of singularities theorem compactifies to a rational curve $C_x$. This suffices to show the existence of $J_W$ planes and now we proceed to describe the relation with the Gromov-Witten invariant.
    
    After stabilizing, the limiting cylinder will break into two pieces, hence belongs to the space $$\mathcal{P}_x \times_{\mathrm{ev}_{\infty}} \mathcal{C}$$

    Conversely, a gluing procedure as described in \cite[Proposition 5]{BourgeoisMohnke2004}, alternatively \cite[Section 5]{pardon2019contact} (with similar results appearing in \cite{HutchingsTaubesGluing1}, \cite{HutchingsTaubes2009Gluing}; see also \cite[Section 7]{abbondandolo2024closedcharacteristicsminimalaction} for an argument very similar to ours in the context of uniformly convex domains in $\mathbb{R}^{2n}$), turns a two-level building like the above into a genuine Floer trajectory for a sufficiently large stretching parameter $R$ (note that both pieces that are glued together are transverse: one is a simple $J_W$ holomorphic plane and the other is a trivial cylinder). We can therefore compare the mod $2$ counts of rigid Floer trajectories after dividing out by reparametrizations: $$|\# \mathcal{M}(e,\hat{\delta}_{1/w_{\max}};H_R,J_R)/\mathbb{R}|=|\# \mathcal{P}_x\times_{\mathrm{ev}_\infty} \mathcal{C}/(\mathbb{R}_1\times \mathbb{R}_2)|$$

    The space on the right is identified with the space of unparametrized planes $\mathcal{P}_x/\mathbb{C}^\times$ since the fiber product happens over the circle $S^1\simeq S_{\delta_{1/w_{\max}}}$ underlying the Reeb orbit. We conclude that the algebraic count of rigid unparametrized $J_W$ planes connecting $x$ to $\delta_{1/w_{\max}}$ is nonzero. Moreover, this count agrees (up to a sign) with the algebraic count of rigid $J_\mathcal{X}$ spheres which have a point insertion, as well as a twisted sector insertion $\mathbf{1}_{1/w_{\max}}$: this is due to the fact that the exceptional orbit lives over the $\frac{1}{w_{\max}}$ point. In other words, the correlator $\langle [\mathrm{pt}_x], \mathbf{1}_{1/w_{\max}}\rangle_{[C]}^\mathcal{X}\neq 0$

    Note that the moduli space of curves in the class $[C]$ with two marked points and one insertion in $\mathbf{1}_{1/w_{\max}}$ has virtual dimension $$2(c_1^\mathcal{X}[C]+n-3+2-\iota)=2n$$so that after imposing an insertion with the point class, the virtual dimension becomes zero, as expected.

     By deformation invariance of Gromov-Witten invariants, we conclude that there is a $J_\mathcal{X}^{int}$-holomorphic sphere contributing to the invariant. A $J_\mathcal{X}^{int}$-holomorphic stable curve in the class $[C]$ contributing to the invariant can have at most one component: the component of the curve that passes through the generic point is not contained in $\mathcal{D}$ and hence meets $\mathcal{D}$ only at the $\frac{1}{w_{max}}$ point. Any other nontrivial component must be disjoint from $\mathcal{D}$ but that contradicts the ampleness of $\mathcal{D}$.
\end{proof}

\begin{cor}
    There is a uniruling by rational curves $C\simeq \mathbb{P}(1,w_{\max})$ passing through the point with maximal isotropy order and satisfying $\mathrm{orbdeg}\,C:=-K_\mathcal{X} \cdot C+n-\iota_{age}=n+1$
\end{cor}
\begin{proof}
    The existence of a uniruling follows from the non-vanishing Gromov-Witten invariant. Moreover, the orbifold degree is given by $$\mathrm{orbdeg}\,C=d_{g^{-1}}\,C=c_1^\mathcal{X}C-\iota_{age}+n=1+n$$
    following the conventions of \cite{LiZhou2026} for the definition of $d_{g^{-1}}$.
\end{proof}

\begin{rem}
    An alternative proof, without using GW invariants and/or gluing, might be possible by showing that there is a sequence of admissible $J_\mathcal{X}$ approximating $J_\mathcal{X}^{int}$ and using Gromov compactness. The stable curve that appears at the limit is forced, as shown at the last part of the proof of Proposition \ref{prop:uniruling}, to have exactly one component.
\end{rem}
\subsection{Showing $|\mathcal{O}(k\mathcal{D})|$ is base point free }

We now wrap up the proof the last part of Theorem \ref{thm:orb_com_classification} where $d=1$, by finding a smooth representative of the linear system $|\mathcal{O}(k\mathcal{D})|$. As mentioned in the introduction, the existence of such a smooth divisor is enough to prove that $\mathcal{X}=\mathbb{P}(1^n,k)$: if such a thing exists, it would (by Kobayashi-Ochiai's theorem), be isomorphic to $\mathbb{P}^{n-1}$. The final blow is then provided by \cite[Lemma 13]{Liu_2017}.

\begin{rem}
     Since \cite[Lemma 13]{Liu_2017} only applies for $n\geq 3$, we need to also show the result holds for $n=2$. But in this case, orbifold del Pezzo surfaces with $\frac{1}{k}(1,1)$ points have been classified by Prince and Cavey \cite{cavey2020del} and the only one which is a rational homology $\mathbb{P}^2$ is indeed $\mathbb{P}(1,1,k)$.
\end{rem}

\begin{prop}
    In the notation of the previous proposition, the linear system $|\mathcal{O}(k\mathcal{D})|$ is base-point free and admits a smooth representative $Z\simeq \mathbb{P}^{n-1}$.
\end{prop}
\begin{proof}
    We assume from now on that $k>1$: the case $k=1$ is the smooth case that follows immediately, since $\mathcal{D}$ itself would then be the desired representative.
    
    Our first step is to compute the Hilbert series. The Grothendieck-Riemann-Roch formula (\cite{Toen1999}) shows that \begin{equation*}
        \chi(\mathcal{X}, \mathcal{O}(m\mathcal{D}))=P(m)+\epsilon(m)
    \end{equation*}
    where $P(m)$ is polynomial in $m$ and $\epsilon(m)$ is an error term depending only on the analytic type of $m\mathcal{D}$ near the orbifold point. As such, $\epsilon(m)=\epsilon(m+k)$. So, $g(m):=\chi(\mathcal{X}, \mathcal{O}(m\mathcal{D}))-\chi(\mathcal{X}, \mathcal{O}((m-k)\mathcal{D}))$ is polynomial in $m$ and of degree $n-1$. By Kodaira vanishing, we know that $g(-1)=g(-2)=\dots = g(-(n-1))=0$. We also know that $g(0)=1$, which shows that $g(m)=\binom{n+m-1}{n-1}$. In particular, for $m=0,1,\dots, k-1$ we have that $\chi(\mathcal{X}, \mathcal{O}(m\mathcal{D}))=h^0(\mathcal{O}(m\mathcal{D}))=\binom{n+m-1}{n-1}$. Then, we get $h^0(\mathcal{O}( k\mathcal{D}))=h^0(\mathcal{O})+\binom{n+k-1}{n-1}=1+\binom{n+k-1}{n-1}$. Continuing like this, we conclude that the Hilbert function is
    \begin{equation*}
        HS_t(\mathcal{X})=HS_t(\mathbb{P}(1^n,k))=\frac{1}{(1-t)^n(1-t^k)}
    \end{equation*}
    and hence in particular that the leading term of the polynomial part is $\mathcal{D}^n/n!$ from which we read off that $\mathcal{D}^n=\frac{1}{k}$.

    Now, let us resolve $X=|\mathcal{X}|$ by $\pi:Y\rightarrow X$, replacing the singular point by an exceptional $E\simeq \mathbb{P}^{n-1}$ with normal bundle $\mathcal{O}(-k)$. Denote by $D$ the underlying coarse space of $\mathcal{D}$. We have that \begin{equation*}\pi^*D=F+\frac{1}{k}E\end{equation*} where $F$ is the strict transform of $D$, since by assumption $\mathcal{D}$ is log smooth at the singular point (i.e. it is a smooth orbifold divisor in a smooth orbifold).

    Moreover, by the discrepancy formula for a $\frac{1}{k}(1,\dots,1)$ singularity, \begin{equation*}-K_Y=\pi^*(-K_X)+\frac{k-n}{k}E=(n+k)F+2E.\end{equation*}

    By adjunction $\mathcal{O}_{\mathbb{P}^{n-1}}(n)=-K_E=(-K_Y-E)|_E=((n+k)F+E)|_E$ which implies that $F|_E=\mathcal{O}_E(1)$.

    The family of rational curves $C$ obtained from Proposition \ref{prop:uniruling} lift via strict transform to rational curves $\tilde{C}$ in $Y$. Moreover, $\pi_* \tilde{C}=C$ and hence \begin{equation*}\tilde{C}\cdot \underbrace{(kF+E)}_{\pi^*(k\mathcal{D})}=C\cdot k\mathcal{D}=1.\end{equation*}

    Since the rational curves $C$ are not contained in $F$ they must intersect it non-negatively, and hence $F\cdot \tilde{C}=0$ and $E\cdot \tilde{C}=1$. 

    Our next goal is to show that $|F|$ is base-point free and as a first step, we show that $F$ is nef. Firstly, there cannot exist an effective divisor in the linear system $|mF-E|$: if one $V\in |mF-E|$ existed, we can find a rational curve $\tilde{C}_x$ for some $x\notin V$ and hence the intersection number will have to be non-negative, contradicting $\tilde{C}\cdot (mF-E)=-1$. Therefore, \begin{equation*}H^0(Y,mF-E)=0, \forall m\geq 0.\end{equation*}
    
    Moreover, since $-F+\pi^*D=\frac{1}{k}E$, by \cite[Lemma 7.30]{KollarSingularities2013} it follows that $\pi_*\mathcal{O}_Y(F)=\mathcal{O}(D).$ On global sections, this implies that \begin{equation*}h^0(\mathcal{O}(F))=h^0(\mathcal{O}(\mathcal{D}))=n.\end{equation*}

    Consider the long exact sequence in cohomology associated to the ideal sheaf sequence \begin{equation*}\mathcal{O}(F-E)\rightarrow \mathcal{O}(F)\rightarrow \mathcal{O}_E(1).\end{equation*}
    From before, we have that $H^0(Y,F-E)=0$ and so there is an injection \begin{equation*}
        H^0(Y,F)\rightarrow H^0(E,F|_E),
    \end{equation*}
    which must be an isomorphism, since $F|_E=\mathcal{O}_{\mathbb{P}^{n-1}}(1)$ and so both groups are of the same rank. Therefore, we conclude that $|F|$ has no basepoints on $E$. 

    Suppose now for the sake of contradiction that $F$ is not nef. Hence, there is some curve $\Gamma$ such that $\Gamma \cdot F<0$. This forces $\Gamma$ to be contained in the base locus of $|F|$ which we just proved is disjoint from $E$. However, $F=\pi^*L-\frac{1}{k}(E)$ so for such curves, $F\cdot \Gamma = \pi^*L \cdot \Gamma\geq 0$ since $L$ is ample, contradiction!

    Now, by Kawamata-Viehweg vanishing and the fact that adding a nef divisor to a big and nef divisor results in a big and nef divisor, we see that \begin{equation*}
        H^1(Y, mF-E)=H^1(Y, K_Y+\underbrace{(kF+E)+(n+m)F}_{big,nef})=0
    \end{equation*}
    In particular, from the long exact sequence associated to $\mathcal{O}(mF-E)\rightarrow \mathcal{O}(mF)\rightarrow \mathcal{O}_E(m)$ we conclude that
    \begin{equation*}
        H^0(Y,mF)\simeq H^0(E, \mathcal{O}(m))
    \end{equation*}
    and in fact since the restriction map is a map of rings, the section rings
    \begin{equation*}
        \bigoplus_{m\geq 0} H^0(Y,\mathcal{O}(mF))\simeq \bigoplus_{m\geq 0} H^0(\mathbb{P}^{n-1}, \mathcal{O}(m))
    \end{equation*}
    are isomorphic and generated in degree $1$. 

    We use the base-point-free theorem of Kawamata \cite[Theorem 3.1.1]{KawamataMatsudaMatsuki1987} to show that $F$ is asymptotically base-point-free: since $-K_Y+aF=2(kF+E)+(n-k+a)F$ is big and nef for $a\geq k-n$, a sufficiently high multiple of $F$ is basepoint free. But, since the section ring of $F$ is generated in degree $1$, it follows that $F$ itself is basepoint free.

    Finally, we argue that $\pi^*(\mathcal{O}(k\mathcal{D}))=\mathcal{O}(kF+E)$ is basepoint free. Again by Kawamata-Viehweg vanishing, $H^1(Y,kF)=0$ so that 
    \begin{equation*}
        0\rightarrow H^0(Y,kF)\rightarrow H^0(Y,kF+E)\rightarrow H^0(E,\mathcal{O}_E)\rightarrow 0
    \end{equation*}
    is exact. This shows that a generic section of $kF+E$ does not vanish anywhere on $E$, so the base locus is disjoint from $E$. Moreover, if $x\in Y\setminus E$, there must be a section $f\in H^0(Y,kF)$ which is nonzero at $x$ (since $|F|$ is base point-free) and hence $s_E f$ is nonzero at $x$, where $s_E\in H^0(Y,E)$ is the defining section for $E$. Thus, $kF+E$ is basepoint-free and Bertini's theorem guarantees that a generic section will be smooth. Moreover, a generic section $Z$ is completely disjoint from $E$, showing that $\pi(Z)$ is a smooth representative of $|\mathcal{O}(kD)|$ in $X$. 
\end{proof}

\part{A periodic structure on V-shaped symplectic cohomology}\label{part:Vshapedperiodicity}

In this part of the paper, we embark on the proof of Theorem \ref{thm:mainthm}. The structure is as follows:
\begin{itemize}
    \item We first describe a monotonicity-type lemma that will ensure that any Floer cylinder or pair of pants which leaves the symplectization region $\partial W\times [1,\infty)$ and goes deep inside the interior must consume a certain amount of energy $\hbar$. This constant depends only on the behaviour of the Hamiltonian before a certain level set, so we can arrange that all our Hamiltonians satisfy this property.
    \item Then, we study the Reeb rotation $\mathcal{S}$ applied to the generators of symplectic resp. V-shaped symplectic cohomology. This rotation does not necessarily commute with the differential on $\mathrm{CF}^\bullet(H)$, but by the monotonicity lemma, it does on $\mathrm{gr}\,\mathrm{CF}^\bullet(H)$. Notably, a rotated Hamiltonian orbit for $H$ becomes a Hamiltonian orbit for a certain $H^{(1)}$ which on the symplectisation region agrees with the pullback of $H$ along the rotation diffeomorphism.
    \item We relate $\mathcal{S}$ to the multiplication by $\mathfrak{M}$ on the associated graded and show that, in a suitable sense, they agree. 
    \item We then construct a cofinal set of Hamiltonians which we use to show that multiplying by $\mathfrak{M}$ is a filtered quasi-isomorphism.
\end{itemize}

\section{Energy Estimates}\label{sec:estimates}
\subsection{Review of monotonicity results for cylinders}

We recall some monotonicity lemmas in the case of cylinders as in \cite[Section 3.5]{hein2010conleyconjectureirrationalsymplectic}. Here, $W$ is a shell which one could take to be of the form $H^{-1}[a,b]$ for some small $a,b$.

\begin{lem}{\cite[Lemma 3.5]{hein2010conleyconjectureirrationalsymplectic}}
Let $W$ be a bounded open set with smooth boundary and at least two boundary components and let $K$ be an autonomous Hamiltonian on $\bar{W}$. Assume that $K$ is constant on each boundary component and does not have one-periodic orbits in $\bar{W}$.
Then there exists a constant $C_1>0$, depending only on the almost complex structure $J$, the open set $W$ and $K$ such that:
\begin{enumerate}
\item For $T\leq 1$, any path $\gamma\colon[0,\,T]\to \bar{W}$, which connects two distinct boundary components of $W$, satisfies
\[
\int_0^T\left\|\dot{\gamma}(t)-X_K(\gamma(t))\right\|^2\ dt >C_1.
\]
\item Any loop $\gamma\colon S^1\to \bar{W}$ satisfies
\[
\int_{S^1}\left\|\dot{\gamma}(t)-X_K(\gamma(t))\right\|^2\ dt >C_1.
\]
\end{enumerate}
\end{lem}

There is also the following:
\begin{lem}{\cite[Lemma 3.6]{hein2010conleyconjectureirrationalsymplectic}, \cite[Lemma 2.3]{Usher2009} }
\label{lem:usher}
Let $W$ be a bounded open set with smooth boundary and at least two boundary components and let $K$ be an autonomous Hamiltonian on $\bar{W}$.
 Let $S$ be a connected subset of the cylinder $S^1\times\mathbb{R}$ and let $u\colon S\to \bar{W}$ satisfy the Floer equation with Hamiltonian $K$. 
Assume that $u(\partial S)\subseteq \partial W$. If $u(S)$ intersects two distinct boundary components of $W$, then there exists a constant $C_2$, depending only on the domain $W$, the Hamiltonian $K$ and the complex structure $J$ on $W$, such that
\[
Area(S)+E(u)\geq C_2.
\]
\end{lem}

The proof of this lemma uses Gromov's trick, replacing $u$ by its graph $\tilde{u}$ which is $J$-holomorphic. The quantity $Area(S)+E(u)$ is then identified with $\int_{S} \tilde{u}^*\tilde{\omega}$, for which a lower bound exists by a monotonicity argument as in \cite{Sikorav1994}.

A combination of these lemmas results in:
\begin{prop}{\cite[Proposition 3.2]{hein2010conleyconjectureirrationalsymplectic}}
Let $K$ be a non-degenerate Hamiltonian and let $U$ and $V$ be open sets that are both bounded by level sets of $K$. Assume furthermore that $K$ does not have one-periodic orbits in $\bar{V}\setminus U$ and is autonomous on this shell.
Let $u\colon S^1\times\mathbb{R}\to M$ be a Floer trajectory that intersects $\partial U$ and $\partial V$.
Then there is a constant $c>0$, only depending on the open sets $U$ and $V$, the restriction of the Hamiltonian $K$ and the almost complex structure $J$ to $\bar{V}\setminus U$, such that $E(u)>c$.
\end{prop}

\subsection{The pair of pants case}

To deduce a lower bound for the energy of pairs of pants passing through the symplectisation region, we will need slight modifications of some of the estimates of \cite{hein2010conleyconjectureirrationalsymplectic}. We will use instead \cite[Proposition 3.3]{groman2026boundarydepthdeformationssymplectic}, which intuitively states that finding a path on $\mathrm{im}(u)$ which covers a lot of length implies the energy of $u$ has to be big. We use the notation $M$ for our target symplectic manifold. Note that our convention for the Hamiltonians $\omega(-,X_H)$ is opposite to Groman's, which accounts for the differences in sign in how we quote the Proposition.

\begin{prop}{(Graph packing with a boundary buffer \cite[Proposition 3.3]{groman2026boundarydepthdeformationssymplectic})} 

\label{prop:Groman_graph_packing}
    Fix a bounded $C^0$-set $B$ of compatible almost complex structures
    and a constant $q>0$. Restrict attention to the cylinder and to the
    fixed finite collection of punctured domains used for the product used in this paper, as well as exterior subdomains thereof.
    Given $d_*,b_*>0$, there exist constants $c_B,C_B>0$ and
    \[
    \tau_0=\tau_0(d_*,b_*)>0,
    \]
    depending also on the fixed geometric data, $B$, and $q$, such that
    the following holds.

    Let $N\subset M$ be a codimension-zero region with boundary and let $S$
    be one of the possibly punctured domains just listed. Let
    \[
    u:S\longrightarrow M
    \]
    solve Floer's equation for a possibly domain-dependent datum $(H,J)$ with
    $J\in B$. Assume that
    \[
    u(\partial S)\subset M\setminus\operatorname{int}N,
    \]
    where $\partial S$ denotes the genuine finite boundary of $S$; punctures
    are not part of $\partial S$. Write $\widehat H:=H-f$ for a function
    $f:S\to\mathbb R$ depending only on the domain (allowing $f=0$).
    Let $\gamma\subset S$ be a path compactly contained in $S$, and set
    \[
    d:=d_M\bigl(u(\gamma(0)),u(\gamma(1))\bigr),
    \qquad
    b:=d_M\bigl(u(\gamma),M\setminus\operatorname{int}N\bigr),
    \] Suppose that
    \[
    d\geq d_*,
    \qquad
    b\geq b_*,
    \qquad
    0<\tau\leq\tau_0,
    \]
    and that
    \[
    \|X_H\|_{C^0(S\times N)}\le q\tau,
    \qquad
    \|d_S\widehat H\|_{C^0(S\times N)}\le q\tau^2,
    \qquad
    d_S\widehat H\wedge\beta \leq 0
    \quad\text{on }S\times N.
    \]
    Then
    \[
    E^{\mathrm{geo}}(u)
    \geq
    c_B\,d\,\tau
    -
    C_B\tau^2\bigl(\ell(\gamma)+1\bigr),
    \]
    where the geometric energy is taken over the (possibly punctured)
    domain $S$.
    \end{prop}
This proposition is similar to Lemma \ref{lem:usher}, with the term $C_B \tau^2(\ell(\gamma)+1)$ taking place of the area. The idea is to combine a lower bound coming from the usual monotonicity for $J$-holomorphic curves applied to the graph $\tilde{u}$, as in Lemma 2 before, together with an upper bound on the area of $S$ to conclude that the geometric energy of $u$ is bounded from below, by choosing $\tau$ very small. We sketch its proof, as in its current form it looks quite technical:

\begin{proof}(\textit{Sketch})
\textbf{Step 1}: Perform the Gromov trick, replacing $u$ by its graph $\tilde{u}$. This is $J$-holomorphic, so the monotonicity lemma of \cite{Sikorav1994} can be applied: there are constants $a,c,\tau_*>0$ such that for all $0<\tau<\tau_*$, every graph ball of radius $a\tau$ centered on $\tilde{u}(S)$ and disjoint
from $\tilde{u}(\partial S)$ has graph energy at least $c\tau^2$ (importantly, one needs that $X_H$ is small to do this).

\textbf{Step 2}: Proceed to choose $\tau_0$ so small so that whenever $0<\tau\leq \tau_0$ and $d\geq d_*, b\geq b_*$, a ball in the graph centered on $\tilde{u}(\gamma)$ is contained in $S\times \mathrm{int}(N)$.

\textbf{Step 3:} The endpoints of $\gamma$ are separated by distance $d$. Thus, one can pack at least $c_1\frac{d}{\tau}$ disjoint balls along it (choosing $c_1$ so that the above is an integer). By monotonicity, these have energy bounded from below by $(c_1\frac{d}{\tau})(c\tau^2)=c_2d\tau$.

\textbf{Step 4:} The area of the neighbourhood of $\gamma$ where the balls are is at most $C(\ell(\gamma)+1)$. 

\textbf{Final step:} Finally, one relates the geometric energy to the graph energy. Let $S'$ be the union of the preimages of the disjoint balls under $\tilde{u}$. There is an identity $$\int_{S'} \tilde{u}^*\omega=E_{geo}(u)+Area(S')\tau^2+\int_S -d_S H\wedge \beta$$

On one hand, the graph energy is bounded from below by $c_2d\tau$. On the other hand, by the assumption that $||d_S \hat{H}||\leq  q\tau^2$, the above is bounded from above by $E_{geo}(u)+C'\tau^2(\ell(\gamma)+1)$. From this, the claim follows.
     
\end{proof}

We now use the above together with an argument originating from \cite{hein2010conleyconjectureirrationalsymplectic} to obtain an energy estimate adapted to our setting, which applies to Floer solutions with domain the cylinder or the pair-of-pants.

\begin{prop}
\label{prop:main_energy_estimate}
    Let $(S,\beta)$ be a punctured surface, equipped with a closed 1-form, cylindrical ends and weights $w_i$ as in Section \ref{pop_prod_def}. Let $B$ be a $C^0$-bounded set of compatible almost-complex structures on $M$ and set $J\in B$. Fix the following data:
    \begin{itemize}
        \item A function $f_S:S\rightarrow [\frac{1}{w_0},\frac{1}{\min_{i=1}^d w_i}]$ such that $df_S\wedge\beta\geq 0$ and $f_S\circ\epsilon_i=\frac{1}{w_i}$ for each $i$.  
        \item A closed interval $[a,b]\subset \mathbb{R}_{\geq 0}$. We will apply Proposition \ref{prop:Groman_graph_packing} to $\rho^{-1}[a,b]$ and the symplectisation coordinate $\rho:\widehat{W}\rightarrow\mathbb{R}_{\geq 0}$.
    \end{itemize}
    
    There exist constants $\sigma_0>0$ and $\delta>0$ and a continuous function $\hbar_{cross}:(0,\sigma_0]\rightarrow\mathbb{R}_+$ with the following property: if \begin{itemize}
        \item $H|_{\rho^{-1}[a,b]}=-\sigma \rho f_S +k$, for some $k:S\rightarrow\mathbb{R}$ such that $||dk||<\delta$ and $dk\wedge\beta\leq 0$
        \item $u:S\rightarrow M$ is a Floer solution for $(H,J)$ such that $u(\partial S)\subset \widehat{W}\setminus \rho^{-1}(a,b)$, and $[a,b]\subset\rho(u(S))$
    \end{itemize} then $E_{geo}(u)\geq \hbar_{cross}(\sigma)$.
\end{prop}
\begin{proof}
\textbf{Step 1:} First we set up the problem so that we can apply Proposition \ref{prop:Groman_graph_packing}. We define four nested subsets,
    \begin{equation*}
        N_j=\rho^{-1}\left[\frac{(6-j)a+jb}{6},\frac{ja+(6-j)b}{6}\right]
    \end{equation*}
    for $j=0,1,2,3$. When applying the previous Proposition we will take $N=N_0=\rho^{-1}[a,b]$. We specify all the constants that are in Proposition \ref{prop:Groman_graph_packing}. Let 
    \begin{equation*}
        d_*=\min_{j\in\{1,2\}}d_M\left(\partial N_j,N_{j+1}\right), \qquad 
        b_*=d_M\left(\partial N_0,N_1\right).
    \end{equation*}
    Fix $q>0$ and let $c_B,C_B,\tau_0$ be the constants obtained from Proposition \ref{prop:Groman_graph_packing}.
    
    Let $S_{thick}$ denote the complement of the cylindrical ends in $S$, which is a connected, compact, codimension-zero submanifold with boundary. Let $T_{\min}$ denote the minimal period of a Reeb orbit in $(\partial W,\alpha)$. Set 
    \begin{equation*}
    \begin{split}
        \tau=\min\left(\tau_0,\frac{c_B d_*}{2C_B\max\{\mathrm{diam}(S_{thick})+1,2\}}\right), \quad \delta=\frac{1}{2}q\tau^2, \quad
        \\\sigma_0=\min\left(\frac{q\tau}{||X_\rho||_{C^0(S\times N_0)}},\frac{q\tau^2}{2b\max\{1,||df_S||_{C^0(S\times N_0)}\}},\frac{T_{\min}}{2}\right),
        \end{split}
    \end{equation*}
    and fix $\sigma\in (0,\sigma_0]$. We have that $X_{\rho}$ has no $\sigma$-periodic orbits in $N_0$.

    By assumption, $H|_{N_0}=-\sigma f_S\rho+k$. At each cylindrical end $X_{H}\otimes \beta=-\sigma X_{\rho}\otimes dt$ and moreover there are inequalities $||X_H||_{C^0(S\times N_0)}\leq q\tau$, $||d_SH||_{C^0(S\times N_0)}\leq q\tau^2$ and $d_SH\wedge\beta\leq 0$ on $S\times N_0$, as required by the assumptions of Proposition \ref{prop:Groman_graph_packing}.
    
    \textbf{Step 2:} Now suppose $u$ is a Floer solution for $(H,J)$ such that $[a,b]\subset \rho(u(S))$. We split into two cases. 
    
    \textbf{Case 1:} First, suppose that $u(S_{thick})$ intersects both $N_3$ and $\partial N_2$. Then there exists a path $\gamma$ in $S_{thick}$ such that $u(\gamma(1))\in\partial N_2$ and $u(\gamma(0))\in N_3$, and $\ell(\gamma)\leq \mathrm{diam}(S_{thick})$. We may restrict $\gamma$ so that the image of $u\circ\gamma$ remains entirely contained in $N_1$ (ensuring that $u\circ\gamma$ first meets $\partial N_2$ at time $1$). Then by Proposition \ref{prop:Groman_graph_packing}, $E_{geo}(u)\geq \frac{1}{2}c_Bd_*\tau$. 
    
    \textbf{Reduction to cylinder case:} If the first case does not hold, the image of some cylindrical end $u\circ \epsilon_i$ must intersect $N_2$, else the image of $u$ is either disjoint from $N_3$ or from $\partial N_2$. We can now reduce to the case of a cylinder and use cylindrical coordinates and follow the argument of Hein in \cite{hein2010conleyconjectureirrationalsymplectic}.

    \textbf{Case 2}: Let $\eta_s(t)=u(\epsilon_i(s,t))$ for $s\in\mathbb{R}$. Fix $(s_0,t_0)\in\mathbb{R}\times\mathbb{R}/\mathbb{Z}$ such that $\eta_{s_0}(t_0)\in N_2$. Either $\eta_{s_0}$ connects $N_2$ to $\partial N_1$, and we can choose $s_1=s_0$, or we can choose $s_1\neq s_0$ such that the image of $\eta_s$ is contained in $N_1$ for all $s\in [s_0,s_1]$, and the image of $\eta_{s_1}$ intersects $\partial N_1$. In either case, we have $\eta_{s_1}(t_1)\in \partial N_1$ and $\eta_{s_1}(t)\in N_1$ for $t\in [t_0,t_1]$, some $t_1\in\mathbb{R}/\mathbb{Z}$. Let $$\gamma(t)=\epsilon_i((1-t)s_0+ts_1,(1-t)t_0+tt_1).$$ We have $\ell(\gamma)\leq 1+|s_1-s_0|$, $u(\gamma(0))\in N_2$, $u(\gamma(1))\in \partial N_1$, and the image of $u\circ\gamma$ is contained in $N_1$.
    
    By an energy-distance inequality (see e.g. \cite[Proposition A.1]{groman2026localfloercohomologyindicator}) for each $s\in(s_0,s_1)$,
    \begin{equation}
    \label{eq:energy_distance}
        \int_{t=0}^1||\dot{\eta}_s(t)+\sigma X_\rho\circ \eta_s(t)||^2dt\geq C(\sigma)
    \end{equation}
    where
    \begin{equation}
    \label{eq:energy_distance_constant}
    C(\sigma)=\frac{\min_{x\in N_0,i}d_M(x,\phi_\rho^{-\sigma }(x))^2}{\max_{|t|\leq \sigma_0 }||d\phi_\rho^{-t}||^2_{C^0(N_0)}},
    \end{equation}
    which is a constant depending only on $\rho$, $\sigma$, $w_0$ and the geometry of $M$ (because by assumption $\eta_s\subset N_0$ for $s\in[s_0,s_1]$), which is positive because $X_\rho$ has no $\sigma$-periodic orbits in the compact region $N_0$ by assumption on $\rho$. 

    By Proposition \ref{prop:Groman_graph_packing} applied to $\gamma$ we have $E_{geo}(u)\geq c_Bd\tau-C_B\tau^2(|s_1-s_0|+2)$. Meanwhile integrating \eqref{eq:energy_distance} yields $E_{geo}(u)\geq |s_1-s_0|C(\sigma)$. Combining these two inequalities, we obtain an energy bound in the second case
    \begin{equation}
        E_{geo}(u)+\frac{C_B\tau^2}{C(\sigma)}E_{geo}(u)\geq c_Bd\tau-C_B\tau^2(|s_1-s_0|+2)+C_B\tau^2|s_1-s_0|=c_Bd\tau-2C_B\tau^2
    \end{equation}
    which yields
  
    \begin{equation}
        E_{geo}(u)\geq \frac{c_Bd_*\tau-2C_B\tau^2}{1+\frac{C_B\tau^2}{C(\sigma)}}\geq\frac{\frac{1}{2}c_Bd_*\tau}{1+\frac{C_B}{C(\sigma)\tau^2}}.
    \end{equation}

    We may take $\hbar_{cross}(\sigma)$ to be the minimum of the bounds obtained in the cases above. The bounds and therefore $\hbar$ all depend only on the geometry of $S$ and $M$, $f_S$, $\beta$, $B$, $\rho$, $\sigma$ and $[a,b]$. Only the third bound depends, continuously, on $\sigma$.
\end{proof}

\begin{rem}
    We will use this result to confine holomorphic curves to the symplectisation region. We will always be able to choose split-monotone Floer data of the required form: in particular, we will require our acceleration data to have a fixed form in a small neck region, up to translation by a constant.
\end{rem}

\section{The rotation operator $\mathcal{S}$}
\label{sec:low_energy_chain_map}
In this section, we construct and analyse an operator $\mathcal{S}$ on the associated graded of the symplectic cochains coming from the circle action on the symplectisation. We first recall a standard construction of the Seidel representation arising from a global Hamiltonian circle action. We would like to apply this construction to the action on the symplectisation. To do this, we use a monotonicity lemma to confine low-energy curves to the symplectisation, then show that the resulting map behaves well with respect to the differential, continuation maps and product on the associated graded (but does not necessarily extend to the full Floer cochain complex!). After this, we consider the chain map given by multiplication by a cocycle $\mathfrak{M}$ under the assumption that its lowest action component is the element $\check{M}_1$, as in Definition \ref{def:principalorbit}. We wish to show that this produces a filtered quasiisomorphism. In order to do so, we study its effect on the associated graded, where it is just given by multiplication by $\check{M}_1$, and we relate it to the rotation operator $\mathcal{S}$.

\subsection{Loops of Hamiltonian Diffeomorphisms}
\label{Floer_pullback_def}
We now recall the standard construction of the Seidel representation $S_g$ associated to a loop of Hamiltonian diffeomorphisms $g$ on a symplectic manifold $(M,\omega)$. For our applications, we will take $M$ to be a subset of the symplectisation $\mathbb{R}_+\times\partial W$, and $g$ to be the circle action generated by the Reeb flow.

Let $g:\mathbb{R}/\mathbb{Z}\rightarrow \mathrm{Ham}(M,\omega)$ be a loop of Hamiltonian diffeomorphisms with $g_0=id$. Suppose $g_t=\phi_{K^g}^t$, for some autonomous Hamiltonian $K^g$.

Given $\mathbb{R}/\mathbb{Z}$-dependent Floer datum $(H,J)$, following \cite[Section 3.2]{ritter_floer_2014}, we can define the pullback $(g^*H,g^*J)$ by letting
\begin{equation*}
\label{eq:g_pullback}
    (g^*H)_t=H_t\circ g_t-K^g,\quad (g^*J)_t=dg_t^{-1}\circ J_t\circ dg_t.
\end{equation*}

There is a bijection
\begin{equation*}
    \mathfrak{X}_{g^*H}\rightarrow \mathfrak{X}_H,\;\gamma\mapsto g\cdot \gamma,
\end{equation*}
where $(g\cdot\gamma)(t)=g_t(\gamma(t))$.

Furthermore there is an orientation-preserving bijection
\begin{equation*}
    \mathcal{M}(\gamma_-,\gamma_+;g^*H,g^*J)\rightarrow \mathcal{M}(g\cdot\gamma_-,g\cdot\gamma_+;H,J),\;u\mapsto g\cdot u,
\end{equation*}
where $(g\cdot u)(s,t)=g_t(u(s,t))$.

Similarly, given $(S,\beta)$, and $S$-dependent Floer data $(H,J)$ as in Section \ref{pop_prod_def}, we can define $(g^*H,g^*J)$ as follows (this construction is slightly non-standard). Let $\tau:S\rightarrow\mathbb{R}/\mathbb{Z}$ be a primitive for $\beta$. Such a $\tau$ exists because $[\beta]\in H^1(S;\mathbb{Z})$ by our choice of weights $w_i$. We are free to rotate the cylindrical ends $\epsilon_i$ so that $\tau\circ\epsilon_i=w_it$. Now let
\begin{equation*}
    (g^*H)_z:=H_z\circ g_{\tau(z)}-K^g,\quad (g^*J)_z:=dg_{\tau(z)}^{-1}\circ J_z\circ dg_{\tau(z)}.
\end{equation*}    

There is a bijection
\begin{equation*}
    \mathcal{M}_{(S,\beta)}(\gamma_0,...,\gamma_d;g^*H,g^*J)\rightarrow \mathcal{M}_{(S,\beta)}(g^{w_0}\cdot\gamma_0,...,g^{w_d}\cdot\gamma_d;H,J),\;u\mapsto g\cdot u,
\end{equation*}
where $(g\cdot u)(z)=g_{\tau(z)}(u(z))$. 
\begin{rem}
    For $w\in\mathbb{Z}$, we can write $g^w_t=g_{wt}$, and let $K^{g^w}=wK^g$, so that $(g^w)^*H=(g^*)^wH$. Note that in the cylindrical ends, this construction agrees with the usual $\mathbb{R}/\mathbb{Z}$-dependent version \eqref{eq:g_pullback}: we have
    \begin{equation*}
        w_i(g^*H)_{\epsilon_i(s,t)}=((g^{w_i})^*H_i)_t,\quad (g^*J)_{\epsilon_i(s,t)}=((g^{w_i})^*J_i)_t.
    \end{equation*}
\end{rem}

\begin{prop}
By construction, $H$ is split-monotone if and only if $g^*H$ is.
\end{prop}

\begin{proof}
For any fixed $x\in M$ and $z\in S$,
\begin{equation*}
    \label{g_pullback_split_monotone}
    d((g^*H)(x))_z-d(H(y))_z=\left(\iota_{\dot{g}_{\tau(z)}(x)}d(H_z-K^g)\right)(y)\beta_z,
\end{equation*}
where $y=g_{\tau(z)}(x)$. The derivatives on the left-hand side are in the $S$-direction, and on the right-hand side, the $M$-direction. It follows that 
\begin{equation*}
    d((g^*H)(x))\wedge\beta=d(H(y))\wedge\beta,
\end{equation*}
so the left-hand side is non-positive for all $x\in M$ if and only if the right-hand side is non-positive for all $y\in M$. 
\end{proof}

\subsection{Confining curves of low energy to the symplectisation}
The loop $g$ arising from the Reeb flow is only defined in the symplectisation and does not necessarily extend into the interior. For this reason, we will need to establish sufficient conditions on a Floer datum so that low-energy solutions are confined to the symplectisation in order to define the action of $g$ on these solutions.
\begin{defn}
\label{def:g}
Assume without loss of generality that the principal orbit in $(\partial W,\alpha)$ has period $1$. Consider the loop of Hamiltonian diffeomorphisms of
$\mathbb{R}_+\times\partial W$ which we denote by $g:\mathbb{R}/\mathbb{Z}\rightarrow \mathrm{Ham}(\mathbb{R}_+\times \partial W)$, so that $g_t=(id,\phi^t)$, where $\phi^t$ denotes the time-$t$ Reeb flow on $\partial W$. Define $K^g=\rho-1$ and denote the Maslov index of the principal orbit by $\mu(dg_t)=2I$.
\end{defn}

The following are sufficient conditions on our choices of $S$-dependent Floer data so that Floer solutions which escape a certain subset of the symplectisation region have energy bounded below by a uniform constant. We require the Hamiltonian and almost-complex structure to be close to a fixed form in two neck regions which separate the symplectisation region from the skeleton (the blue regions in Figure \ref{fig:Hamiltonians}). The level sets $\frac{1}{16},\frac{1}{8},\frac{1}{4},\frac{1}{2}$ bounding these regions are arbitrary, but it is important that they are fixed.

\begin{defn}[Neck-Admissible almost complex structures]
\label{allowed_J}
Fix a compatible almost-complex structure $J_{cyl}$ on $\widehat{W}$, which is cylindrical and of contact type on $\rho^{-1}[\frac{1}{4},\infty)$. Let $B$ denote a fixed $C^0$ ball centred on $J_{cyl}$. An $S$-dependent almost-complex structure $J$ is neck-admissible if for all $z\in S$,
\begin{itemize}
\item $J_z\in B$,   
\item $J_z$ is of contact type at $\infty$.
\end{itemize}    
\end{defn}

\begin{defn}[Neck-Admissible Hamiltonians]
\label{allowed_H}
Let $(S,\beta)$ be either the cylinder or the pair-of-pants surface, equipped with fixed weights $w_1=w_2=1$, and a fixed choice of closed 1-form as in Section \ref{pop_prod_def}. Fix the following data:
\begin{itemize}
    \item A normalising function $f_S$ as in Proposition \ref{prop:main_energy_estimate}: in the cylinder case, we take $f_S=1$; in the pair-of-pants case, we have $f_S:S\rightarrow [\frac{1}{2},1]$.
    \item Fix constants $\delta>0$, $\sigma_0>0$ and a function $\hbar:(0,\sigma_0]\rightarrow\mathbb{R}_+$, which satisfy the statement of Proposition \ref{prop:main_energy_estimate}, simultaneously for $S$ a cylinder or pair-of-pants, the ball of almost-complex structures $B$, and both the intervals $[\frac{1}{16},\frac{1}{8}]$ and $[\frac{1}{4},\frac{1}{2}]$
    \item Fix also once and for all $\sigma\in (0,\sigma_0]$ and let $\hbar<\min\{\hbar_{cross}(\sigma), \frac{1}{2N}\}$.
\end{itemize} 

An $S$-dependent Hamiltonian $H$ is neck-admissible if there is some $b\in\mathbb{Z}$, functions $k_i:S\rightarrow \mathbb{R}$ such that $dk_i\wedge\beta\leq 0$ and $||dk_i||<\delta$ for $i=1,2$, and so that for all $z\in S$, 
\begin{itemize}
\item $H_z|_{\rho^{-1}[\frac{1}{16},\frac{1}{8}]}=-\sigma f_S(z)\rho+k_1(z)$,
\item $H_z|_{\rho^{-1}[\frac{1}{4},\frac{1}{2}]}=(-\sigma f_S(z) +w)\rho+k_2(z)$,
\item $H_z|_{\partial W}<0$,
\item $H_z-\mu(z)\rho$ is constant outside a compact set, for some $\mu$.
\end{itemize}
\end{defn}

The function $\hbar$ depending on $\sigma$ serves as a lower bound on the energy of curves which cross into the interior of $W$.
\begin{rem}
    The freedom to choose $b\in\mathbb{Z}$ is necessary, because we will later require pairs of Floer data to be related by pulling back by $g$ in the region $\rho^{-1}[\frac{1}{4},\infty)$, which changes the slope by $1$. Because we only consider a finite collection of operations, we can choose $\hbar$ and $\sigma_0$ once and for all.
\end{rem}

We now show that low-energy Floer solutions for neck-admissible $(H,J)$ cannot cross from the symplectisation region to the filling.

\begin{prop}[Curve confinement]
\label{admissible_curve_confinement}
Let $(S,\beta)$ and $\hbar>0$ be as in Definition \ref{allowed_H}. Suppose $(H,J)$ is a neck-admissible $S$-dependent Floer datum, and $u:S\rightarrow\widehat{W}$ is a solution of Floer's equation \eqref{Floer_eq_pop} such that the image $u(S)$ intersects $\rho^{-1}[\frac{1}{2},\infty)$. Then, either $u(S)\subset \rho^{-1}[\frac{1}{4},\infty)$, or $E_{geo}(u)>\hbar$.
\end{prop}

\begin{proof}
Suppose $u$ is a Floer solution, $E_{geo}(u)\leq \hbar$, and $u(S)$ intersects $\rho^{-1}[\frac{1}{2},\infty)$. By Proposition \ref{prop:main_energy_estimate}, applied to the region $\rho^{-1}[\frac{1}{16},\frac{1}{8}]$, $u(S)\subset \rho^{-1}[\frac{1}{16},\infty)$ since if $u$ were to reach inside $\rho^{-1}[0,\frac{1}{16})$ it would necessarily consume more energy than $\hbar$. We can therefore define $g^{-b}\cdot u$, and Proposition \ref{prop:main_energy_estimate} applies to $g^{-b}\cdot u$ as well in the region $\rho^{-1}[\frac{1}{4},\frac{1}{2}]$, so $g^{-b}\cdot u(S)\subset \rho^{-1}[\frac{1}{4},\infty)$. Since $g^{-b}$ preserves level sets of $\rho$, the Floer solution $u(S)\subset \rho^{-1}[\frac{1}{4},\infty)$ as required.
\end{proof}

\begin{rem}
    We cannot just apply Proposition \ref{prop:Groman_graph_packing} to the region $\rho^{-1}[\frac{1}{4},\frac{1}{2}]$ since we need to retain control of both $H$ and a $w$-related version $H^{(w)}$ as in \ref{def:g_related}. In other words, in Figure \ref{fig:Hamiltonians}, the Hamiltonians $H$ and $H^{(w)}$ can be guaranteed to be small simultaneously only in the leftmost blue region.
\end{rem}

Now we define a sufficient condition on pairs of Floer data so that Floer solutions confined to the symplectisation are related by $g^w$. This will allow us to later define an action of $g$ which intertwines the low-energy parts of the differential, continuation maps and pair-of-pants product on the Floer complexes.
\begin{defn}
\label{def:g_related}
A pair of $S$-dependent Hamiltonians and almost-complex structures $(H,J)$, $(H^{(w)},J^{(w)})$ is $w$-related if
\begin{equation*}
\label{eq:w_related}
    H^{(w)}|_{\rho^{-1}[\frac{1}{4},\infty)}=(g^w)^*H|_{\rho^{-1}[\frac{1}{4},\infty)},\quad J^{(w)}|_{\rho^{-1}[\frac{1}{4},\infty)}=(g^w)^*J|_{\rho^{-1}[\frac{1}{4},\infty)}.
\end{equation*}
\end{defn}

\begin{rem}
Note that the sets of neck-admissible $w$-related pairs of neck-admissible Hamiltonians (for fixed $w$) and almost-complex structures are convex and contractible respectively.
\end{rem}

\begin{prop}
We are free to choose neck-admissible $w$-related pairs of regular, split-monotone Floer data.
\end{prop}
\begin{proof}
    We are free to make separate, $S$-dependent perturbations of $(H,J)$ supported inside $\rho^{-1}[0,\frac{1}{16})$, $\rho^{-1}(\frac{1}{8},\frac{1}{4})$ and $\rho^{-1}[\frac{1}{2},\infty)$ such that any Floer solution with image intersecting one of these regions is regular. We can similarly adjust $(H^{(w)},J^{(w)})$ independently in the first two regions. In the third, this condition is satisfied automatically, because the Cauchy-Riemann operators are related by $dg_t^w$. Any Floer solution intersects one of these regions, because they contain all 1-periodic orbits.
\end{proof}

\subsection{A cofinal class of Hamiltonians $H_{b,\mu}$}
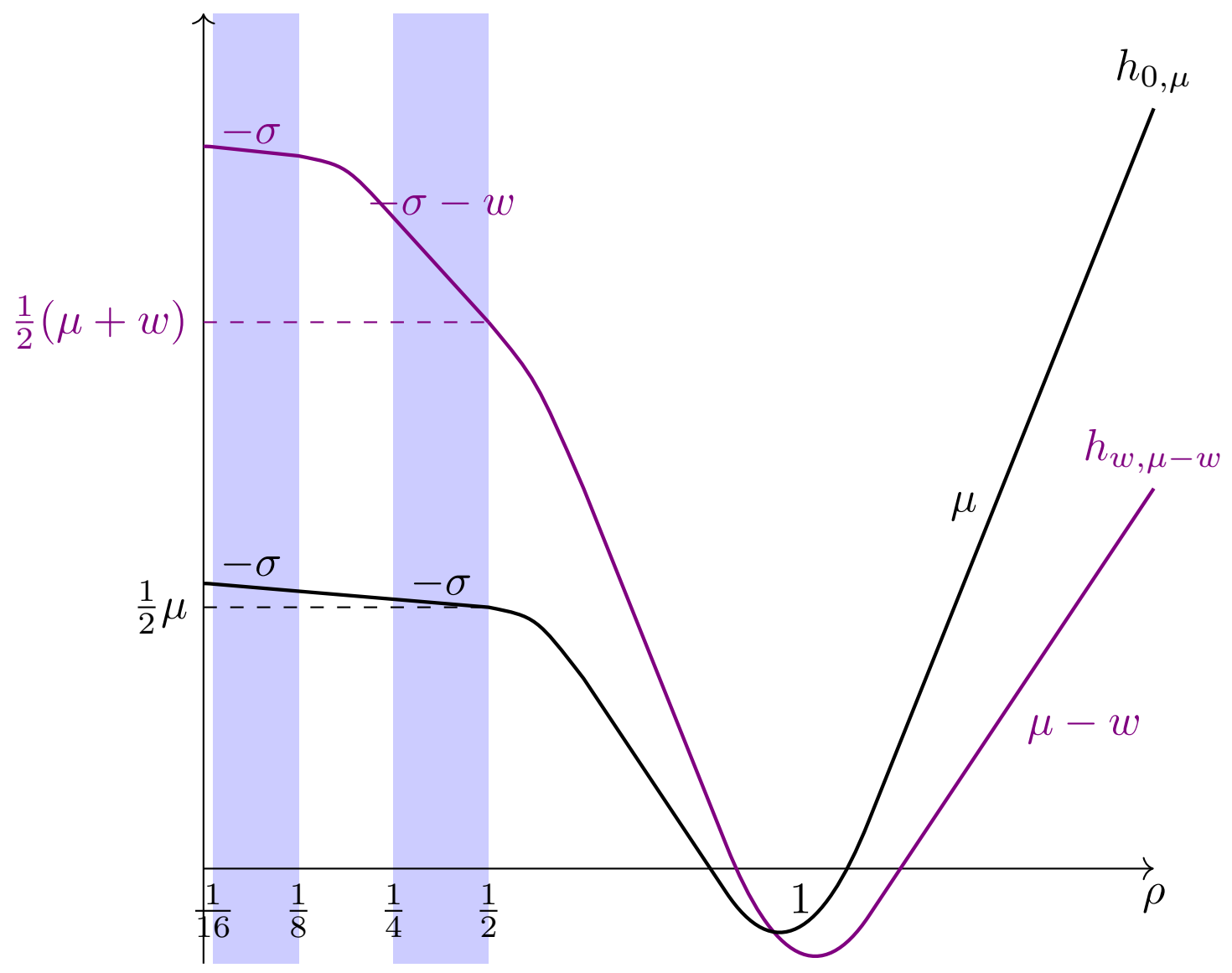
\begin{figure}[H]
\caption{The functions $h_{w,\mu-w}$ and $h_{0,\mu}$ after smoothing, for some $w\in\mathbb{Z}\cap [0,\mu]$, with labelled slopes.}
\label{fig:Hamiltonians}
\[\begin{tikzpicture}[scale=0.75]
    \filldraw[draw=none,opacity=0.2,fill=blue] (0.1,0) rectangle (1,10);
    \filldraw[draw=none,opacity=0.2,fill=blue] (2,0) rectangle (3,10);

    \draw[->] (0,0)--(0,10);
    \draw[->] (0,1)--(10,1);
    \node at (10,1)[anchor=north] {$\rho$};
    \draw[violet, thick] (3,6.75)--(2,7.85).. controls (1.5,8.4).. (1,8.5)--(0.1,8.595)..controls (0.05,8.6)..(0,8.6);
    \draw[violet, thick] (3,6.75).. controls (3.5,6.15)  .. (4,5)--(5.5,1.25).. controls (6,0) and (6.5,-0.25).. (7,0.5)--(10,5);
    \draw[thick] (0,4).. controls (0.05,4).. (0.1,3.995)--(3,3.75).. controls (3.5,3.65)  .. (4,3)--(5.5,0.75).. controls (6,0) and (6.5,0.25).. (7,1.5)--(10,9);
    \node at (10,9)[anchor=south]{$h_{0,\mu}$};
    \node at (10,5)[anchor=south,violet]{$h_{w,\mu-w}$};
    \node at (2,1)[anchor=north]{$\frac{1}{4}$};
    \node at (3,1)[anchor=north]{$\frac{1}{2}$};
    \node at (1,1)[anchor=north]{$\frac{1}{8}$};
    \node at (0.1,1)[anchor=north]{$\frac{1}{16}$};
    \node at (6,1)[anchor=north west]{$1$};

    \node at (0,3.75)[anchor=east]{$\frac{1}{2}\mu$};
    \draw[dashed] (0,3.75)--(3,3.75);

    \node at (0,6.75)[anchor=east,violet]{$\frac{1}{2}(\mu+w)$};
    \draw[dashed,violet] (0,6.75)--(3,6.75);
    
    \node at (0.5,4.2){$-\sigma$};
    \node at (2.5,4){$-\sigma$};

    \node at (0.5,8.75)[violet]{$-\sigma$};
    \node at (2.5,8)[violet]{$-\sigma-w$};

    \node at (8,4.5)[anchor=south]{$\mu$};
    \node at (8.5,2.5)[violet,anchor=west]{$\mu-w$};

\end{tikzpicture}\]
\end{figure}

We now define a cofinal class of time-dependent Hamiltonians $H_{b,\mu}$ for V-shaped symplectic cohomology, such that the orbits with action $\leq b$ are confined to the symplectisation region (in particular, these give rise to admissible Hamiltonians on the cylinder defining the Floer differential). These Hamiltonians are Morse-Bott perturbations of certain profile functions of the symplectisation coordinate $\rho$. We will use these Hamiltonians to define $\mathcal{S}$. The only main differences between our definition of these Hamiltonians and the standard one in \cite{CieliebakFrauenfelderOancea2010} and \cite{uebele2019periodicreebflowsproducts} are the conditions involving $\sigma$, which are used to ensure that the resulting Hamiltonians on the cylinder meet our admissibility condition, Definition \ref{allowed_H}.
\begin{prop}
\label{H_eta_mu}
Fix  $\mu>0$ not in the period spectrum of $(\partial W,\alpha)$, and $b\in\mathbb{Z}\cap [0,\mu]$. Let $\hbar,\sigma>0$ be as in Definition \ref{allowed_H}.

Consider a piecewise linear function $h_{b,\mu}=h$ such that
\begin{equation}
    \begin{cases}
        h'(\rho)=-\sigma, & \rho\in (0,\frac{3}{16})
        \\ h'(\rho)=-\sigma-b, &\rho \in (\frac{3}{16},\frac{3}{4})
        \\ h(\frac{1}{2})\geq b+\frac{1}{2}\mu,
        \\ h(1)\in (-\epsilon,0),
        \\ h'(\rho)=\mu, & \rho> 1
    \end{cases}
\end{equation}
up to a concave corner smoothing at $\frac{3}{16}$ and $\frac{3}{4}$, and a convex smoothing at $1$, without changing $h$ in the regions $[\frac{1}{16},\frac{1}{8}]$ and $[\frac{1}{4},\frac{1}{2}]$.
(see Figure \ref{fig:Hamiltonians}).

Let $H_{b,\mu}$ be a non-degenerate $\mathbb{R}/\mathbb{Z}$-dependent Hamiltonian which is a sufficiently small Morse-Bott perturbation of $h_{b,\mu}\circ \rho$, supported in sufficiently small disjoint neighbourhoods of the $1$-periodic orbits, disjoint from $\rho^{-1}[\frac{1}{16},\frac{1}{8}]$, $\rho^{-1}[\frac{1}{4},\frac{1}{2}]$.

Then $H_{b,\mu}$ is $b$-admissible and neck-admissible (viewed as an $s$-independent Hamiltonian on the cylinder $(\mathbb{R}/\mathbb{Z})_t\times\mathbb{R}_s$). Furthermore, all lower orbits of $H_{b,\mu}$ are contained in $\rho^{-1}(\frac{1}{2},\infty)$.
\end{prop}
\begin{proof}
The conditions imposed on $h=h_{b,\mu}$ imply that the upper orbits of the autonomous Hamiltonian $h\circ \rho$ satisfy $h\circ\rho>b$, and hence have action $>b$, and all orbits contained in $\rho^{-1}[0,\frac{1}{2}]$ are upper orbits. This condition is preserved by a sufficiently small Morse-Bott perturbation. Neck-admissibility follows immediately from the conditions imposed in the regions $\rho^{-1}[\frac{1}{16},\frac{1}{8}]$ and $\rho^{-1}[\frac{1}{4},\frac{1}{2}]$.\end{proof}

\begin{rem}
It is clear that for fixed $b$, the set of piecewise linear functions $h_{b,\mu}$ is convex. We can therefore always interpolate between Hamiltonians $H_{b,\mu}$, by first interpolating linearly between functions $h_{b,\mu}$, then smoothing, then adding a Morse-Bott perturbation.
\end{rem}

\begin{prop}
For each $b\in\mathbb{Z}_{\geq 0}$, the set of Hamiltonians of the form $H_{b,\mu}$ (for varying $\mu$) described in Proposition \ref{H_eta_mu} is cofinal.
\end{prop}
\begin{proof}
It suffices to show that the piecewise linear Hamiltonians of the form $h_{b,\mu}\circ \rho$ are cofinal. Given $K$ linear at $\infty$ on $\widehat{W}$ such that $K|_{\partial W}<0$, choose $\delta\in (0,\frac{1}{4})$ and $h_{b,\mu}(1)\in (-\epsilon,0)$ such that $K|_{\rho^{-1}[1-\delta,1+\delta]}<h_{b,\mu}(1)$. Choose $c>0$ such that $K|_{\rho^{-1}[c,\infty)}$ is a linear function of $\rho$. As $\mu\to\infty$, $\min_{[0,c]\setminus [1-\delta,1+\delta]}h_{b,\mu}\to\infty$, so we can choose $\mu>b$, larger than the slope of $K$ at $\infty$, such that $h_{b,\mu}\circ\rho\geq K$ on $[0,c]$. 
\end{proof}

Next, we describe how to choose pairs of time-dependent Hamiltonians $(H,H^{(w)})$ where $H$ and $H^{(w)}$ are as in the previous Proposition \ref{H_eta_mu}, so that $g^w$ relates the 1-periodic orbits with action below a given bound.

\begin{prop}
\label{hw}
Let $w,b\in\mathbb{Z}$ and $0\leq w\leq b\leq \mu-w$. Let $H_{b-w,\mu}$ be as in Proposition \ref{H_eta_mu}. We can choose $H_{b,\mu-w}$ also as in Proposition \ref{H_eta_mu}, such that $H_{b-w,\mu}$ and $H_{b,\mu-w}$ are $w$-related.

Furthermore there is a bijection
\begin{equation*}
    \mathcal{A}_{H_{b,\mu-w}}^{-1}(-\infty,b)\rightarrow\mathcal{A}_{H_{b-w,\mu}}^{-1}(-\infty,b-w),\;\gamma\mapsto g^w\cdot\gamma.
\end{equation*}
\end{prop}
\begin{proof}
Given a piecewise linear Hamiltonian $h_{b-w,\mu}$, we can choose $h_{b,\mu-w}$ also piecewise linear by letting
\begin{equation*}
    h_{b,\mu-w}(\rho)=\begin{cases}
        h_{b-w,\mu}(\rho)-w(\rho-1), & \rho\geq \frac{3}{16}
        \\ h_{b-w,\mu}(\rho)+\frac{13}{16}w, & \rho\leq \frac{3}{16}.
    \end{cases}
\end{equation*}
We can ensure that this relation still holds in $[\frac{1}{4},\infty)$ after corner smoothing, so the autonomous Hamiltonians $h_{b,\mu-w}$ and $h_{b-w,\mu}$ are $w$-related. Given a non-degenerate perturbation $H=H_{b-w,\mu}$ of $h_{b-w,\mu}\circ\rho$ as in Proposition \ref{H_eta_mu}, we can choose a non-degenerate perturbation $H^{(w)}=H_{b,\mu-w}$ of $h_{b,\mu-w}\circ\rho$ which is $w$-related, because any such $H^{(w)}$ has no 1-periodic orbits in $\rho^{-1}[\frac{1}{4},\frac{1}{2}]$, and all orbits in $\rho^{-1}[\frac{1}{4},\infty)$ are automatically non-degenerate, so we are free to add a small perturbation supported in $\rho^{-1}[0,\frac{1}{4})$. 

By Proposition \ref{H_eta_mu}, all orbits in $\mathfrak{X}_H$, $\mathfrak{X}_{H^{(w)}}$ with action in the relevant intervals are contained in $\rho^{-1}[\frac{1}{2},\infty)$. In this region $H^{(w)}$ agrees with $(g^w)^*H$, so the map $\gamma\mapsto g^w\cdot\gamma$, which decreases action by $w$, restricts to a well-defined bijection.
\end{proof}

We proceed to show that low-energy Floer solutions for $(H,J)$ and $(H^{(w)},J^{(w)})$ in $\widehat{W}$, asymptotic to the lower orbits (upper orbits are excluded on the basis of their high action) can be identified via composing with $g$. The point is that low-energy Floer solutions are confined in the symplectisation region, whence we can argue as in Section \ref{Floer_pullback_def}.

\begin{prop}
\label{curve_confinement}
Let $(S,\beta)$ be as in Definition \ref{allowed_H}.

Suppose we have for $i=0,...,d$, neck-admissible $w_i$-related pairs of $\mathbb{R}/\mathbb{Z}$-dependent Floer data $(H_i,J_i)$ and $(H_i^{(w_i)},J_i^{(w_i)})$, where $H_i=H_{b_i-w_i,\mu_i}$ and $H^{(w_i)}=H_{b_i,\mu_i-w_i}$ as in Proposition \ref{hw}. Let $(H,J)$ and $(H^{(1)},J^{(1)})$ be a $1$-related pair of admissible split-monotone $S$-dependent Floer data which agree with $(H_i,J_i)$ and $(H^{(w_i)}_i,\allowbreak J^{(w_i)}_i)$ respectively in the cylindrical ends, as in Equation \eqref{cyl_ends_consistency}.

Suppose $\gamma_i\in\mathcal{A}_{H_i^{(w_i)}}^{-1}(-\infty,b_i)$ for $i=0,...,d$ and 
\begin{equation*}
\mathcal{A}_{H_0^{(w_0)}}(\gamma_0)-\sum_{i=1}^d \mathcal{A}_{H_i^{(w_i)}}(\gamma_i)\leq \hbar.
\end{equation*}

Then there is an orientation-preserving bijection
\begin{equation*}
    \mathcal{M}_{(S,\beta)}(\gamma_0,...,\gamma_d;H^{(1)},J^{(1)})\rightarrow \mathcal{M}_{(S,\beta)}(g^{w_0}\cdot\gamma_0,...,g^{w_d}\cdot\gamma_d;H,J),\;u\mapsto g\cdot u.
\end{equation*}

If the moduli spaces are regular and zero-dimensional, this bijection preserves the sign of the induced map on orientation lines at each point.
\end{prop}

\begin{proof}
    Let $\hbar$ be as in Proposition \ref{admissible_curve_confinement}. Let $u$ be an element of the right or left-hand side. By split-monotonicity and our assumption on the action of the orbits, $E_{geo}(u)\leq E_{top}(u)\leq \hbar$ and the image of $u$ intersects the region $\rho^{-1}[\frac{1}{2},\infty)$, because the asymptotes are contained there. Therefore the image of $u$ is in fact entirely contained in $\rho^{-1}[\frac{1}{4},\infty)$ by Proposition \ref{admissible_curve_confinement}. In this region, $(H^{(1)},J^{(1)})$ agrees with $(g^*H,g^*J)$, so the left and right-hand moduli spaces are indeed related by $u\mapsto g\cdot u$.

    We may identify the orientation lines for $\gamma_i$ with those of $g^{w_i}\cdot\gamma_i$ via conjugation by $dg_t^{w_i}$, and the maps between orientation lines induced by $u$ and $g\cdot u$, by $dg_t$, so the signs are preserved.
\end{proof}

\subsection{The low-energy rotation operator $\mathcal{S}$ and its properties}
We now use the previous sections to construct a map on the associated graded pieces of symplectic cohomology which is given by the action of $g$ on the generators and show that it behaves well with respect to the differential, continuation maps, and product in action windows of width $<\hbar$; in particular, the operations on the associated graded. We denote these maps $\mathcal{S}$. They will be defined between action-truncated complexes $\mathrm{CF}_{(-\infty,b)}^\bullet(H)$ for certain $b$-admissible $H$, and do not respect the higher energy contributions to the various operations. This construction follows that in \cite{uebele2019periodicreebflowsproducts}, but instead of using index-positivity to confine holomorphic curves to compact subsets of $\mathbb{R}_+\times\partial W$, we use action truncation. 
\begin{defn}
\label{S_g_def}
Let $(H,J)$ and $(H^{(1)},J^{(1)})$ be a $1$-related pair of neck-admissible $\mathbb{R}/\mathbb{Z}$-dependent Floer data, where, in the notation of Proposition \ref{H_eta_mu}, 
\begin{equation*}
H=H_{b-1,\mu},\quad H^{(1)}=H_{b,\mu-1},
\end{equation*}
which can be achieved by Proposition \ref{hw}. We can define a linear isomorphism (\textbf{not} necessarily a chain map!)
\begin{equation*}
\label{eq:S_g}
\mathcal{S}:\mathrm{CF}_{(-\infty,b)}^\bullet(H^{(1)},J^{(1)})\rightarrow \mathrm{CF}_{(-\infty,b-1)}^{\bullet-2I}(H,J)
\end{equation*}
by letting $\mathcal{S}(\gamma)=g\cdot\gamma$. Recall that by definition all generators of the left-hand complex have action in the interval $(-\infty,b)$,
and the underlying orbits are contained in $\rho^{-1}[\frac{1}{2},\infty)$ and this map is well-defined. In other words, these conditions ensure that there are no upper orbits in $\mathrm{CF}_{(-\infty,b)}^\bullet(H^{(1)})$.

From now on we will always assume that our Floer data satisfies the above conditions, so that the map \eqref{eq:S_g} is well-defined whenever we invoke it. 
\end{defn}

\begin{prop}
    The map $\mathcal{S}$ respects the period filtration $\mathcal{F}$ up to a shift,
    \begin{equation}
        \mathcal{S}(\mathcal{F}^p\mathrm{CF}_{(-\infty,b)}^\bullet(H^{(1)},J^{(1)}))=\mathcal{F}^{p-2N}\mathrm{CF}^{\bullet-2I}_{(-\infty,b-1)}(H,J).
    \end{equation}
\end{prop}
\begin{proof}
    Because $\mathcal{A}_H(g\cdot\gamma)=\mathcal{A}_{H^{(1)}}(\gamma)-1$, $\mathcal{S}$ decreases action by $-1$, and therefore decreases the filtration $\mathcal{F}$ by $-2N$.
\end{proof}
\begin{rem}
It is clear that, on the level of the associated graded, $\mathcal{S}(e^{H^{(1)}})=\check{M}_1^{H}$.
\end{rem}

\begin{rem}
    The filtration $\mathcal{F}$ has windows on the associated graded $\mathrm{gr}^p=\mathcal{F}^p/\mathcal{F}^{p+1}$ which are either empty, or have action width $<\hbar$, so that all the curves which contribute to the zeroth page of the action spectral sequence have energy $<\hbar$.
\end{rem}

We now apply Proposition \ref{curve_confinement} to show that $\mathcal{S}$ intertwines the operations on the zeroth page:
\begin{prop}[Compatibility with the differential]
    $\mathrm{gr}^p\mathcal{S}$ is a chain map.
    \label{S_g_chain_map}
\end{prop}
\begin{proof}
Let $\gamma_\pm\in\mathfrak{X}_{H^{(w)}}$ be in a given nonzero window $\mathcal{F}^p/\mathcal{F}^{p+1}$. Then $
\mathcal{A}_{H^{(1)}}(\gamma_+)-\mathcal{A}_{H^{(1)}}(\gamma_-)< \hbar$. By Proposition \ref{curve_confinement}, there is an orientation-preserving bijection
\begin{equation*}
    \mathcal{M}(\gamma_+,\gamma_-;H^{(1)},J^{(1)})\xrightarrow{\sim} \mathcal{M}(g\cdot\gamma_+,g\cdot\gamma_-;H,J),\; u \mapsto g\cdot u.
    \label{low_energy_bijection}
\end{equation*}
\end{proof}

\begin{prop}[Compatibility with continuation maps]
\label{S_g_cont_maps}
$\mathcal{S}$ commutes with continuation maps in the following way. 

Suppose $H_+,H^{(1)}_+\geq H_-,H^{(1)}_-$ respectively. Then the following diagram commutes up to filtered homotopy:
\begin{equation*}
\begin{tikzcd}
\mathrm{gr}^p \mathrm{CF}_{(-\infty,b)}^{\bullet}(H_-^{(1)},J_-^{(1)}) \arrow[r,"\mathrm{gr}^p\kappa"]  \arrow[d,"\mathrm{gr}^{p}\mathcal{S}" '] & \mathrm{gr}^p\mathrm{CF}_{(-\infty,b)}^{\bullet}(H_+^{(1)},J_+^{(1)}) \arrow[d,"\mathrm{gr}^{p}\mathcal{S}"]
\\ \mathrm{gr}^{p-2N}\mathrm{CF}_{(-\infty,b-1)}^{\bullet-2I}(H_-,J_-)\arrow[r,"\mathrm{gr}^{p-2N}\kappa" '] & \mathrm{gr}^{p-2N}\mathrm{CF}_{(-\infty,b-1)}^{\bullet-2I}(H_+,J_+)
\end{tikzcd}   
\end{equation*}
where the horizontal maps are those induced on the zeroth page of the spectral sequence by (filtered) continuation maps $\kappa$ and arising from some choices of monotone homotopies.
\end{prop}

\begin{proof}
By a standard argument, we are free to choose the monotone homotopies giving rise to the continuation maps $\kappa$, and we may subdivide these homotopies and replace the continuation maps $\kappa$ with the corresponding compositions of continuation maps (the resulting maps are related by a filtered chain homotopy). We may therefore assume that the differences $H_+-H_-$ and $H_+^{(1)}-H_-^{(1)}$ are arbitrarily small on restriction to the regions $\rho^{-1}[\frac{1}{16},\frac{1}{8}]$ and $\rho^{-1}[\frac{1}{4},\frac{1}{2}]$. Under this assumption, we can choose a $w$-related pair of admissible, regular, monotone homotopies $(H,J)$ and $(H^{(1)},J^{(1)})$ between $(H_\pm,J_\pm)$ and $(H_\pm^{(1)},J_\pm^{(1)})$ respectively, to define the continuation maps $\kappa$.

Then, by the same argument as Proposition \ref{S_g_chain_map}, whenever $
\mathcal{A}_{H^{(1)}_+}(\gamma_+)-\mathcal{A}_{H^{(1)}_-}(\gamma_-)< \hbar$, there is an orientation-preserving bijection
\begin{equation*}
    \mathcal{M}(\gamma_+,\gamma_-;H^{(1)},J^{(1)})\xrightarrow{\sim}\mathcal{M}(g\cdot\gamma_+,g\cdot\gamma_-;H,J),\; u\mapsto g\cdot u,
\end{equation*}
so under our extra assumptions on $H,H^{(1)}$, the diagram strictly commutes.
\end{proof}
\begin{rem}
    The preceding two Propositions show that $\mathrm{gr}^p\mathcal{S}$ is independent of choices, up to homotopy.
\end{rem}
\begin{prop}[Compatibility with the product]
\label{S_g_product}
The maps $\mathrm{gr}^p\mathcal{S}$ intertwine the pair-of-pants product in the following way. Suppose $b_0\leq \min\{a(p)+b_2,a(q)+b_1\}$.

Given $H_i,H_i^{(1)}$ for $i=1,2$, we can choose $H_0,H_0^{(2)}$ such that the following diagram commutes up to filtered homotopy:
\begin{equation}
\label{eq:S_g_product}
\begin{tikzcd}
    \mathrm{gr}^p\mathrm{CF}^{\bullet}_{(-\infty,b_1)}(H_1^{(1)},J_1^{(1)})\otimes \mathrm{gr}^q\mathrm{CF}_{(-\infty,b_2)}^{\bullet}(H_2^{(1)},J_2^{(1)}) \arrow[r,"\mathrm{gr} \mu"]  \arrow[d,"\mathrm{gr}^p\mathcal{S}\otimes\mathrm{gr}^q\mathcal{S}" '] & \mathrm{gr}^{p+q}\mathrm{CF}_{(-\infty,b_0)}^{\bullet}(H_0^{(2)},J_0^{(2)})\arrow[d,"\mathrm{gr}^{p+q}\mathcal{S}^2"]
    \\ \mathrm{gr}^{p-2N}\mathrm{CF}_{(-\infty,b_1-1)}^{\bullet-2I}(H_1,J_1)\otimes \mathrm{gr}^{q-2N}\mathrm{CF}_{(-\infty,b_2-1)}^{\bullet-2I}(H_2,J_2) \arrow[r,"\mathrm{gr}^{p+q-4N}\mu" '] & \mathrm{gr}^{p+q-4N}\mathrm{CF}_{(-\infty,b_0-2)}^{\bullet-4I}(H_0,J_0),
\end{tikzcd}
\end{equation}
where the horizontal maps are induced by the (filtered) chain-level products $\mu$ arising from some choices of split-monotone Floer data on the pair-of-pants surface.
\end{prop}
\begin{rem}
This statement could be generalised to the higher operations, and with arbitrary choices of weights on the inputs, but we will only use this case.
\end{rem}
\begin{proof}
The chain-level product indeed descends to the associated graded by Proposition \ref{prop:mu_respects_F}. By a standard argument, the product commutes with continuation maps up to filtered homotopy. By Proposition \ref{S_g_cont_maps}, we may therefore assume, after replacing the input Hamiltonians by larger ones, that $H_1=H_2$ in the regions $\rho^{-1}[\frac{1}{16},\frac{1}{8}]$ and $\rho^{-1}[\frac{1}{4},\frac{1}{2}]$, and choose output Hamiltonians such that the inequalities \eqref{eq:split_monotone_suff_cond} are satisfied, and $H_0+\sigma\rho=2(H_i+\sigma\rho)$ in the two regions for $i=1,2$. We can therefore choose split-monotone $1$-related Hamiltonians $H,H^{(1)}$ on the pair-of-pants $S$, which have the form $-f_S\sigma\rho+w\rho$, up to a constant shift, in the two regions (outside the two regions, we can use convex interpolation). We use these Hamiltonians, and some choice of regular 1-related almost-complex structures, to define $\mu$.

Let $\gamma_0\in\mathfrak{X}_{H_{0}^{(2)}}$, $\gamma_1\in\mathfrak{X}_{H_{1}^{(1)}}$, and $\gamma_2\in\mathfrak{X}_{H_{2}^{(1)}}$.
Arguing as before, any curve which contributes on the associated graded has energy $ <\hbar$ 
so that by Proposition \ref{curve_confinement}, there is an orientation-preserving bijection
\begin{equation*}
    \mathcal{M}_{(S,\beta)}( \gamma_0,\gamma_1,\gamma_2;H^{(1)},J^{(1)})\xrightarrow{\sim}\mathcal{M}_{(S,\beta)}(g^2\cdot\gamma_0,g\cdot\gamma_1,g\cdot\gamma_2;H,J),\;u\mapsto g\cdot u,
\end{equation*}
and the coefficients of $\gamma_0$ in $\mu(\gamma_1,\gamma_2)$ and of $g^2\cdot\gamma_0$ in $\mu(g\cdot\gamma_1,g\cdot\gamma_2)$ are obtained by counting these moduli spaces, so under our extra assumptions on $H,H^{(1)}$, the diagram strictly commutes.
\end{proof}

\begin{rem}
    Note that the multiplication operator is not strictly speaking compatible with the filtration globally, for all $p,q$. However, the multiplication by the principal orbit is.
\end{rem}

\begin{prop}[$\mathrm{gr}^p\mathcal{S}$ agrees with multiplication by the principal orbit, up to homotopy]
\label{M_1_low_energy_product}
For any $H_-$ of the form in Proposition \ref{H_eta_mu}, there exists $H_+\geq H_-$ such that
\begin{equation*}
    m=\mu(\check{M}_1^{H_-},\;\cdot\;):\mathrm{CF}_{(-\infty,b)}^{\bullet}(H_-,J_-)\rightarrow \mathrm{CF}_{(-\infty,b-1)}^{\bullet-2I}(H_+,J_+),
\end{equation*}
respects the filtration up to a shift: $m(\mathcal{F}^p\mathrm{CF}_{(-\infty,b)}^\bullet(H_-))\subset \mathcal{F}^{p-2N}\mathrm{CF}_{(-\infty,b-1)}(H_+)$, where $\mu$ is a choice of filtered chain-level product as in Proposition \ref{S_g_product}, and the following diagram commutes up to filtered homotopy:
\begin{equation*}
\label{m_low_energy_cont_map}
\begin{tikzcd}
    \mathrm{gr}^{p}\mathrm{CF}^{\bullet}_{(-\infty,b)}(H_-^{(1)},J_-^{(1)}) \arrow[r,"{\mathrm{gr}^p\kappa}"]  \arrow[d,"\mathrm{gr}^p\mathcal{S}" '] & \mathrm{gr}^{p}\mathrm{CF}^{\bullet}_{(-\infty,b)}(H_+^{(2)},J_+^{(2)}) \arrow[d,"\mathrm{gr}^p\mathcal{S}^2"]
    \\ \mathrm{gr}^{p-2N}\mathrm{CF}^{\bullet-2I}_{(-\infty,b-1)}(H_-,J_-) \arrow[r,"\mathrm{gr}^{p-2N}  m" '] & \mathrm{gr}^{p-4N}\mathrm{CF}^{\bullet-4I}_{(-\infty,b-2)}(H_+,J_+),
\end{tikzcd}
\end{equation*}
where $\kappa$ is a continuation map.
\end{prop}
\begin{proof}
By Proposition \ref{prop:M1isfiltered}, the operator $m$ is filtered.
By Proposition \ref{chain_level_unit} there is a chain-level representative of the unit $e^{H_-^{(1)}}$ such that $e^{H_-^{(1)}}-\mathcal{S}^{-1}(\check{M}_1^{H_-})\in\mathcal{F}^1\mathrm{CF}^0_{(-\infty,b)}(H^{(1)})$. First we apply Proposition \ref{S_g_product}, restricted to the case where one of the inputs to $\mu$ is $\check{M}_1^{H_-}$, $e^{H_-^{(1)}}$ respectively, then by Proposition \ref{prop:unit_mult_cont_map}, we can replace $\mu(e^{H^{(1)}_-},\;\cdot\;)$ by $\kappa$.
\end{proof}

\section{Proof of the main theorem}
\label{sec:main_thm}
In this section, we prove the main theorem about periodicity. The strategy is to define a chain level operator given by multiplication by $\mathfrak{M}$, which is a chain map provided that $\mathfrak{M}$ is closed. Then, we show that this chain map is a filtered quasi isomorphism by comparing the map it induces on the associated graded to the map $\mathcal{S}$ from the previous section.

To prove the main theorem, we will use the standard fact that filtered continuation maps induce a quasi-isomorphism on the associated graded of the telescope complex.
\begin{prop}[Continuation maps are filtered quasi-isomorphisms on the telescope]
    \label{prop:cont_maps_filtered_qiso}
    Suppose there are two cofinal sequences of Hamiltonians $H_{\pm,k}$ giving rise to filtered 1-rays $\mathcal{C}_{\pm}$ resp. (together with a choice of almost-complex structures and monotone homotopies) such that $H_{+,k}\geq H_{k,-}$ for all $k$. Moreover, suppose we have a choice of continuation map $\kappa:\mathrm{CF}^\bullet(H_{-,k})\rightarrow\mathrm{CF}^\bullet(H_{+,k})$ arising from a monotone homotopy, for each $k$. Let $\kappa:\mathcal{C}_-\rightarrow\mathcal{C}_+$ denote the map of 1-rays arising from these continuation maps, together with chain homotopies arising from homotopies of monotone continuation data.
    
    The filtration $\mathcal{F}$ on the action-truncated complexes $\mathrm{CF}^\bullet_{(-\infty,b)}(H_{\pm,k})$ extends to the truncated telescopes $\mathrm{tel}(\mathcal{C}_\pm)_{(-\infty,b)}$, and $\kappa$ respects this filtration: $\kappa(\mathcal{F}^p\mathrm{tel}(\mathcal{C}_-)_{(-\infty,b)})\subset \mathcal{F}^p\mathrm{tel}(\mathcal{C}_+)_{(-\infty,b)}$. The induced maps 
     \begin{equation}
         \mathrm{gr}^p\mathrm{tel}(\kappa):\mathrm{gr}^p\mathrm{tel}(\mathcal{C}_-)_{(-\infty,b)}\rightarrow \mathrm{gr}^p\mathrm{tel}(\mathcal{C}_+)_{(-\infty,b)}
     \end{equation}
     on the associated graded are quasi-isomorphisms.
\end{prop}
\begin{proof}
The telescope differential, and $\kappa$, are filtered because they arise from monotone homotopies, which increase action and therefore $\mathcal{F}$, for fixed $b$. 

We can inductively choose an increasing sequence $(s_k)_{k\in\mathbb{Z}_{\geq 0}}$ such that $H_{-,s_k}\geq H_{+,k}$. We can choose further continuation maps $\kappa':\mathrm{CF}^\bullet(H_{+,k})\rightarrow\mathrm{CF}^\bullet(H_{-,s_k})$, and upgrade this choice to a map of 1-rays $\kappa':\mathcal{C}_+\rightarrow\mathcal{C}_-^s$, where $\mathcal{C}_-^s$ is the subray of $\mathcal{C}_-$ induced by $s$. By restricting $\kappa$ to the subrays induced by $s$, we obtain a map of 1-rays $\kappa^s:\mathrm{tel}(\mathcal{C}_-^s)\rightarrow\mathrm{tel}\mathcal({C}_+^s)$.

By a standard argument  there are chain homotopies connecting $\mathrm{tel}(\kappa')\circ\mathrm{tel}(\kappa)$, and $\mathrm{tel}(\kappa^s)\circ \mathrm{tel}(\kappa')$, to the compression maps induced by $s$ on $\mathrm{tel}(\mathcal{C}_\pm)$ respectively. We can ensure that all the maps we have defined arise from monotone homotopies, and therefore increase action. The induced maps on the action-truncated telescopes therefore increase the filtration $\mathcal{F}$. The chain homotopies therefore descend to the associated graded. The compression maps descend to compression maps, which induce the identity on cohomology, and $\mathrm{gr}^p\mathrm{tel}(\kappa^s)$ induces the same map on cohomology as $\mathrm{gr}^p\mathrm{tel}(\kappa)$. We therefore have a left and right inverse for $H(\mathrm{gr}^p\mathrm{tel}(\kappa))$ as required.  
\end{proof}

\mainthm*

\begin{proof}
    It suffices to show that the multiplication map
    \begin{equation}
        [\mathfrak{M}]\;\cdot\;:\widecheck{\mathrm{SH}}^{\bullet}_{(-\infty,b)}(W)\rightarrow \widecheck{\mathrm{SH}}^{\bullet-2I}_{(-\infty,b-1)}(W)
    \end{equation}
    is an isomorphism for $b\gg 0$. This is well-defined because we can choose representatives of $[\mathfrak{M}]$ with action arbitrarily close to $-1$. Now fix $b\in\mathbb{Z}_{\geq 0}$. We will construct a chain model for the multiplication map as follows. Choose cofinal sequences $(H_{k,\pm},J_{k,\pm})$ for $k\in\mathbb{Z}_{\geq 0}$ where
    \begin{itemize}
        \item $H_{-,k}=H_{b,\mu_k}$ as in Proposition \ref{H_eta_mu},
        \item $\mu_{k+1}>\mu_k\geq k$,
        \item $H_{\pm,k}$ are as in Proposition \ref{M_1_low_energy_product},
        \item $H_{-,k+1}\geq H_{+,k}$.
    \end{itemize}
    
    Now choose $1$-related and $2$-related sequences $H_{k,-}^{(1)}$ and $H_{k,+}^{(2)}$ respectively, as in Proposition \ref{hw}. We will take the Floer complexes associated to $H_{\pm,k}$ to be the domain resp. target of the maps $m$, which agree to first order with multiplication by $\mathfrak{M}$, and the complexes associated to $H_{-,k}^{(1)}$ and $H_{+,k}^{(2)}$ to be domain resp. target of the continuation maps $\kappa$. These maps are then related using Proposition \ref{M_1_low_energy_product}. Together with appropriate choices of almost-complex structures and continuation maps, the sequences $H_{\pm,k}^{(i)}$ give rise to 1-rays $\mathcal{C}_\pm^{(i)}$. The resulting action-truncated telescopes $\mathrm{tel}(\mathcal{C}_\pm^{(i)})_{(-\infty,b)}$ are filtered by $\mathcal{F}$, and by definition, are models for $\widecheck{\mathrm{SH}}^\bullet_{(-\infty,b)}(W)$.
    
    Because $\mathfrak{M}-\check{M}_1^{H_-}\in\mathcal{F}^{-2N+1}\mathrm{CF}^{-2I}_{(-\infty,b-1)}(H_-)$, Proposition \ref{M_1_low_energy_product} holds with $\mathfrak{M}$ in place of $\check{M}_1^H$. Together with a choice of chain homotopies witnessing the compatibility of the product with continuation maps, arising from a choice of split-monotone homotopies of Floer data, the maps $m$ assemble into a map of 1-rays $m:\mathcal{C}_-\rightarrow\mathcal{C}_+$. By Proposition \ref{prop:mu_respects_F}, $\mathrm{tel}(m)$ respects the filtration $\mathcal{F}$ on the action-truncated telescopes, up to a shift: $\mathrm{tel}(m)(\mathcal{F}^p\mathrm{tel}(\mathcal{C}_-)_{(-\infty,b-1)})\subset \mathcal{F}^{p-2N}\mathrm{tel}(\mathcal{C}_+)_{(-\infty,b-2)}$. By definition, $\mathrm{tel}(m)$ descends on cohomology to multiplication by $[\mathfrak{M}]$. Similarly, as in Proposition \ref{prop:cont_maps_filtered_qiso}, we obtain a filtered chain map $\mathrm{tel}(\kappa):\mathrm{tel}(\mathcal{C}_-^{(1)})_{(-\infty,b)}\rightarrow \mathrm{tel}(\mathcal{C}_+^{(2)})_{(-\infty,b)}$.
    
    To verify that $\mathrm{tel}(m)$ is a quasi isomorphism, by \cite[Lemma 5.6]{Campos_2024} it suffices to show that the map $\mathrm{gr}^p\mathrm{tel}(m)$ between the associated graded pieces of the telescopes is a quasi isomorphism for all $p\in\mathbb{Z}$. This is a simple case of the Lemma because the filtration $\mathcal{F}$ on the action-truncated telescopes is bounded above, hence complete. It suffices to show that $H(\mathrm{gr}^p\mathrm{tel}(m))$ is an isomorphism. The induced map on cohomology is the colimit of the induced maps on cohomology at each level $k$. By Proposition \ref{M_1_low_energy_product}, on the level of cohomology we have a commuting diagram
    \begin{equation*}
    \label{1_ray_diagram}
    \begin{tikzcd}
    H^\bullet(\mathrm{gr}^{p}\mathrm{tel}(\mathcal{C}_-^{(1)})_{(-\infty,b)})\arrow[rrr,"H(\mathrm{gr}^p\mathrm{tel}(\kappa))"]  \arrow[d,"H(\mathrm{gr}^p\mathcal{S})" '] &&& H^\bullet(\mathrm{gr}^{p}\mathrm{tel}(\mathcal{C}_+^{(2)})_{(-\infty,b)})\arrow[d,"H(\mathrm{gr}^p\mathcal{S}^2)"]
    \\ H^{\bullet-2I}(\mathrm{gr}^{p-2N}\mathrm{tel}(\mathcal{C}_-)_{(-\infty,b-1)}) \arrow[rrr,"H(\mathrm{gr}^{p-2N} \mathrm{tel}(m))"] &&&H^{\bullet-4I}(\mathrm{gr}^{p-4N}\mathrm{tel}(\mathcal{C}_+)_{(-\infty,b-2)}),
    \end{tikzcd}
    \end{equation*}
    The vertical maps are isomorphisms by construction. The top map is an isomorphism by Proposition \ref{prop:cont_maps_filtered_qiso}, so the bottom map is an isomorphism as required.
\end{proof}

\section{Constructing a cocycle $\mathfrak{M}$}\label{sec:constructingmathfrakM}

In this section, we construct a cocycle $\mathfrak{M}$ as in Theorem \ref{thm:mainthm} for many types of compactifications. The existence of such a cocycle is expected for virtually any example, if the transversality assumptions of \cite[Proposition 3]{bourgeois_exact_2009} could be somehow circumvented. Instead of dealing with transversality issues, we use our action and index computations to rule out various Floer cylinders to construct the class $\mathfrak{M}$.

First, we establish sufficient conditions such that if one takes a J-shaped Hamiltonian, as in the definition of symplectic cohomology (not V-shaped!), the principal orbit is closed. This is slightly easier than for a V-shaped Hamiltonian, since we do not have to worry about the upper orbits. We can then take $\mathfrak{M}$ to be the resulting cocycle.
\begin{prop}\label{prop:closednessofprincipalorbit}
    Suppose $(\mathcal{X},\mathcal{D})$ is a compactification of $W$ with $\mathcal{D}$ smooth as an orbifold and whose isotropy locus is in codimension one. Assume moreover that $\langle c_1^\mathcal{X},-\rangle=\lambda \langle \mathcal{D}, -\rangle$ on $H_2(\mathcal{X};\mathbb{Q})$ with $\lambda>1$.
    Let $H$ be a J-shaped Hamiltonian of slope $\geq 1$ at $\infty$. Then, $\check{M}_1^H\in \mathrm{CF}^{-2I}(H)$ is closed.
\end{prop}
\begin{proof}
    The assumptions imply that $c_1^W$ is torsion so there is a well-defined grading. Moreover, the class of the principal orbit is torsion in $H_1(W;\mathbb{Z})$\footnote{To see this, one can either use the orbifold Gysin sequence or alternatively use the existence of closed surface in $\mathcal{X}$ avoiding the singularities and hitting the divisor at finitely many points and with non-zero algebraic intersection number.} so the grading formulas can be applied to any orbit in the same homotopy class as $\check{M}_1$. 

    Suppose there is a Floer cylinder between $\check{M}_1$ and $\gamma$. Since the action filtration is compatible with the period and $\check{M}_1$ corresponds to the minimum in the Morse-Bott family of period $1$ orbits, it implies that the period of $\gamma$ is less than or equal to that of $\check{M}_1$. 
    
    Comparing mod $2$ gradings, one sees that $\gamma$ has to correspond to an odd index critical point in its Morse-Bott family. This rules out the possibility that $\gamma$ is constant, since the degree of $\check{M}_1$ is strictly negative and the degree of a constant orbit of odd Morse index is strictly positive and integral. We continue with the case that the period of $\gamma$ is $\frac{i}{w}<1$ for non-constant $\gamma$. We can then write: $$\deg(\check{M}_1)+1=-2(\lambda-1)+1=\deg(\gamma)=-2\left(\frac{i}{w}(\lambda-1)-\sum \frac{\overline{ib_j}}{w}\right)+\underbrace{(\dim_\mathbb{R} \mathcal{N}/S^1+1-\mu_{Morse}(q)-\mu_{S^1})}_{\geq 1}$$
    and therefore $$2\left(\frac{i}{w}-1\right)(\lambda-1)\geq 2\left(\left(\frac{i}{w}-1\right)(\lambda-1)-\iota_{age}\right)\geq 0$$ which is impossible, since the period of $\gamma$ is $\frac{i}{w}<1$ and $\lambda>1$.

    Second, if the period of $\gamma$ is equal to that of $\check{M}_1$, it means it must be part of the Morse-Bott family of principal orbits. But by definition, $\check{M}_1$ corresponds to an orbit which is closed in the local Floer homology of principal orbits and so there can be no contributions to the differential in this case either.

\end{proof}

We also show a slightly different criterion for the existence of $\mathfrak{M}$:
\begin{prop}\label{prop:closednessusinglocalFloer}
Suppose that: \begin{itemize}
    \item $\langle c_1^\mathcal{X},-\rangle=\lambda \langle \mathcal{D}, -\rangle$ on $H_2(\mathcal{X};\mathbb{Q})$ with $\lambda\neq 1$
    \item $c_1^W=0$
    is torsion
\end{itemize} and finally that for a J-shaped Hamiltonian of sufficiently high slope $\tau \geq 1$ at $\infty$, the first page of the Morse-Bott spectral sequence for $H$ 
computes the degree $-2I$ part of symplectic cohomology in the homology class of $\check{M}_1$, i.e. $$\bigoplus_{p+q=-2I} E_{1,[\check{M}_1]}^{p,q}\simeq \mathrm{SH}^{-2I}_{[\check{M}_1]}(W)$$  Then, the principal orbit $\check{M}_1^H$ can be promoted to a cocycle $\mathfrak{M}=\check{M}^H_1+c\in\mathrm{CF}^{-2I}(H)$, where $c\in\mathrm{CF}^{-2I}_{[-1,\infty)}(H)$ so that all of its terms have higher action than $\check{M}_1$.
\end{prop}
\begin{rem}
    Note that we are not assuming that the first page degenerates completely and computes the whole of the symplectic cohomology, but rather that just the degree $-2I$ component degenerates. We will use this criterion when computing the symplectic cohomology of $W$ which compactify to K3 surfaces in Section \ref{sec:K3}.
\end{rem}

\begin{proof}
    By the assumption on $\lambda$, $I\neq 0$ and so for a sufficiently high $\tau$, all orbits of degree $-2I$ homologous to $\check{M}_1$ will be accounted for in $\mathrm{CF}^{-2I}(H)$. This follows by the formulas for the Conley-Zehnder index in Section \ref{sec:grading_formulae}, in particular \ref{cor:explicitcyclicexample} and \ref{cor: principalorbitdegree}.
    The proof follows by analysing the properties of Morse-Bott spectral sequence and some homological algebra. We filter the chain complex as $$\mathcal{F}^p \mathrm{CF}(H):=\mathrm{CF}_{[c(p),\infty)}(H)$$ where the increasing function $c:\mathbb{Z}\rightarrow \mathbb{R}$ is defined as in \cite[Section 8.1]{kwon_brieskorn_2016} (note that we are using the opposite action convention and are working with cohomology rather than homology). This is a decreasing filtration with $\mathcal{F}^{p+1}\mathrm{CF}(H)\subset \mathcal{F}^{p}\mathrm{CF}(H)$ and it exhausts the complex in finitely many steps. By assumption, the first page of the associated spectral sequence computes the symplectic cohomology in the fixed degree $-2I$. Therefore, denoting by $[\check{M}_1]^i$ the class of the principal orbit on the $i$-th page, we have: $$d_0\check{M}_1=d_1[\check{M}_1]^1=\dots=0$$ where $d_i:E^{p,q}_i\rightarrow E^{p+i, q-i+1}_i$ is the differential on the $i$-th page. 
    
    In particular, we see that $$d \check{M}_1 \in \mathcal{F}^{p+2}\mathrm{CF}^{p+q+1}+d\mathcal{F}^{p+1}\mathrm{CF}^{p+q}$$
    so that there must be an element $c_1 \in \mathcal{F}^{p+1}C^{p+q}$ whose terms have higher action than $\check{M}_1$ and such that $d(\check{M}_1+c_1)\in \mathcal{F}^{p+2}\mathrm{CF}^{p+q+1}$. Continuing like this inductively, noting that eventually $\mathcal{F}^{p+m}\mathrm{CF}=0, m\gg 0$, one obtains the desired cocycle $\mathfrak{M}:=\check{M}_1+c$.
\end{proof}

\part{Examples and computations of symplectic cohomology}\label{part:computations}

\section{Symplectic cohomology of Milnor fibers of surface singularities}\label{sec:surfacecomputations}
\subsection{The $A_{p-1}$ Milnor fiber}
\label{sec:Ap-1}
Consider the Milnor fiber of the surface singularity $xy+z^p=0$. This compactifies to 
\begin{equation*}
\{xy+z^p+w^p=0\} \subset\mathbb{P}(p-1,1,1,1)
\end{equation*}
with compactifying divisor a hyperplane section which is an orbifold $\mathcal{D}=\mathbb{P}^1_{1,p-1}$. This is a monotone log-Fano compactification with index $\lambda=2$ since $-K_\mathcal{X}=\mathcal{O}_\mathcal{X}(2)$ whereas $\mathcal{D}$ is a hyperplane section i.e. a section of $\mathcal{O}_{\mathcal{X}}(1)$. Moreover, $\mathcal{X}$ has a single orbifold point of type $\frac{1}{p-1}(1,1)$.

As such, for the canonical contact structure on $\mathbf{S}\mathcal{N}_\mathcal{D}$ there is one exceptional orbit and an $S^2$ family of principal orbits. For the Hamiltonian orbits after a Morse-Bott perturbation, we denote by $\check{\delta}_{i/p-1}$ (resp. $\hat{\delta}_{i/p-1}$) the Hamiltonian orbit of type $\lor$ (resp. $\wedge$) corresponding to the $i$-th cover of the exceptional Reeb orbit, and moreover $\check{m}_i, \check{M}_i$ (resp. $\hat{m_i},\hat{M}_i$) the two Hamiltonian orbits of type $\lor$ (resp. $\wedge$) corresponding to the minimum and maximum of the $S^2$ family of Reeb orbits. 

Our first step is to compute the positive equivariant symplectic cohomology. This will give us enough knowledge to compute $\mathrm{SH}^1$ and $\mathrm{SH}^2$. Then, we will apply Theorem \ref{thm:monotonecompactifications} and show that in fact in degrees $\leq 2$, symplectic cohomology is $-2$-periodic i.e. $\mathrm{SH}^i\simeq \mathrm{SH}^{i-2}$. In other words, periodicity, together with knowledge of $\mathrm{SH}$ in a finite interval, is enough to determine the whole of $\mathrm{SH}$.

\subsubsection{Positive equivariant symplectic cohomology and the computation of $\mathrm{SH}^1$ and $\mathrm{SH}^2$}

From the formulae of Section \ref{sec:grading_formulae}, we have that 
\begin{equation*}
    \mu_{CZ}(\delta_{i/p-1})=2\left((\lambda-1)\frac{i}{p-1}-\frac{\overline{i}}{p-1}\right)+(n-1)=2\left\lfloor \frac{i}{p-1}\right\rfloor+1
\end{equation*}
Moreover, the index of the family of principal orbits is $\mu_{CZ}(\gamma_i)=2i$. In the first page of the Morse-Bott spectral sequence for positive equivariant symplectic cohomology, a component $S\subset I\mathcal{D}$ parametrizing orbits contributes $H^\bullet(S)$ shifted in degree by $-(\mu_{CZ}-(n-1)+\frac{1}{2}\dim_\mathbb{R} S)+1$. Hence, we can write the contributions to the first page (where the spectral sequence degenerates) as:
\begin{table}[H]
\begin{center}
\begin{tabular}{llccccl}
    \toprule
     \dots & \(SH^{-2}_{+,S^1}\)  &\(SH^{-1}_{+,S^1}\)          &\(SH^0_{+,S^1}\)          &\(SH^1_{+,S^1}\)          & \(SH^2_{+,S^1}\) & \(SH^{3}_{+,S^1}\) \\
    \midrule
     \dots&                   &\(\delta_{p/p-1},\dots,\delta_{2p-3/
     p-1}\)              &             &\(\delta_{1/p-1},\dots,\delta_{p-2/p-1}\) \\
     
  \dots& &\(m_2\)  & & \(m_1\) & \\
  
     \dots&  &\(M_1\)  & &  & \\
    \bottomrule
  \end{tabular}
\end{center}
\end{table}
This is supported only on odd degrees $\leq 1$. Namely, $\mathrm{SH}^{1}_{+,S^1}(W)=\mathbb{Q}^{p-1}$ and $\mathrm{SH}^{1-2i}_{+,S^1}(W)=\mathbb{Q}^p, i\geq 1$. This coincides with the computation of \cite{kwon_brieskorn_2016}. By the Gysin sequence for positive equivariant symplectic cohomology, we conclude that $\mathrm{SH}^1_+\simeq \mathrm{SH}^1_{+,S^1}=\mathbb{Q}^{p-1}$, as well as $\mathrm{SH}^2_+(W)=0$. Moreover, by the LES for positive symplectic cohomology, we see that $H^2(W)$ surjects onto $\mathrm{SH}^2(W)$ and that $\mathrm{rk}\,\mathrm{SH}^2(W)=\mathrm{rk}\,\mathrm{SH}^1(W)$. Because $H^2(W)=\mathbb{Q}^{p-1}$, the rank of $\mathrm{SH}^2$ is at most $p-1$. However, the fact that there are $p-1$ Lagrangian spheres in $W$ and the compatibility under Viterbo restriction shows that in fact this rank has to be exactly $p-1$.
\subsubsection{Using the $2$-periodicity and the Cieliebak-Frauenfelder-Oancea long exact sequence to wrap up}

We now introduce the V-shaped symplectic cohomology as an intermediary in the computation. Denote by $[\check{M}_1] \in \widecheck{\mathrm{SH}}^{-2}(W)$ the class of the principal orbit which induces the periodicity $\widecheck{\mathrm{SH}}^{\bullet}(W)\simeq \widecheck{\mathrm{SH}}^{\bullet-2}(W)$. Recall the long exact sequence in V-shaped symplectic cohomology \cite{CieliebakFrauenfelderOancea2010}:
\begin{equation*}
\begin{tikzcd}
\underbrace{\mathrm{SH}_3}_0 \arrow[r] & \mathrm{SH}^1 \arrow[r]                & \widecheck{\mathrm{SH}}^1 \arrow[lld] \\
\mathrm{SH}_2 \arrow[r]                & \mathrm{SH}^2 \arrow[r]                & \widecheck{\mathrm{SH}}^2 \arrow[lld] \\
\mathrm{SH}_1 \arrow[r]                & \underbrace{\mathrm{SH}^3}_0 \arrow[r] & \widecheck{\mathrm{SH}}^3            
\end{tikzcd}
\end{equation*}
We have that $\mathrm{SH}^2(W)\simeq H^2(W)$ and dually  $H^2(W,\partial W)\simeq H_2(W)\simeq \mathrm{SH}_2(W)$. The connecting map $\mathrm{SH}_2(W)\rightarrow \mathrm{SH}^2(W)$ is identified, by \cite{CieliebakFrauenfelderOancea2010}, with the map $H^2(W,\partial W)\rightarrow H^2(W)$. The kernel of this map is $H^1(\partial W)=0$ and its cokernel is $H^2(\partial W)=0$ (note that this is over $\mathbb{Q}$; over the integers, this cokernel is $\mathbb{Z}/p\mathbb{Z}$), so the map is an isomorphism. This implies that 
\begin{equation*}
    \mathrm{SH}^1\simeq \widecheck{\mathrm{SH}}^1, \quad \mathrm{SH}_1\simeq \widecheck{\mathrm{SH}}^2.
\end{equation*}
These are both of rank $p-1$. By the $2$-periodicity of V-shaped symplectic cohomology, which is Theorem \ref{thm:mainthm}, it follows that $\widecheck{\mathrm{SH}}^\bullet(W)=\mathbb{Q}^{p-1}$ in all degrees. Finally, in degrees $\leq 0$, the sequence of \cite{CieliebakFrauenfelderOancea2010} shows that $\mathrm{SH}^\bullet\simeq \widecheck{\mathrm{SH}}^\bullet, \bullet \leq 0$ and we can conclude that 
\begin{equation*}
    \mathrm{SH}^\bullet(W)\cong\mathbb{Q}^{p-1}, \bullet \leq 2
\end{equation*}
recovering the result of \cite{EtguLekili2017}.

\subsection{Rational homology balls $B_{p,q}$}
A particularly interesting example consists of the rational homology smoothings of $T$-singularities of type $\frac{1}{p^2}(1,pq-1)$. These are frequently denoted $B_{p,q}$, and there is a wide literature on them (see e.g. \cite{LM14} where it is shown that the symplectic cohomology is non-vanishing, as well as \cite{EvansSmith2018Markov} which studies the symplectic embeddings of these into $\mathbb{P}^2$, and \cite{adaloglou2026nearbylagrangianconjecturepinwheels} which proves a nearby Lagrangian conjecture result).

We will describe them as follows: consider the hyperquotient T-singularity $$xy+z^p=0\in \mathbb{A}^3/\frac{1}{p}(a,p-a,1)$$where $a\equiv q^{-1} \mathrm{mod}\, p,\,0\leq a <p$.
Add another parameter $t$ and compactify into $$\{xy+z^p+\varepsilon t=0\}\subset \mathbb{P}(a,p-a,1,p)$$
At $\varepsilon=0$, one gets the T-singularity in an affine chart, but, partially smoothing by setting $\varepsilon\neq 0$, one obtains a compactification of $B_{p,q}$. One can eliminate the variable $t$ and thus identify $$B_{p,q}=\mathbb{P}(a,p-a,1)\setminus \mathcal{D}_p,\quad \mathcal{D}_p\subset \mathbb{P}(a,p-a,1)$$

The $B_{p,q}$ have the rational homology of a point, but they are not simply-connected: $\pi_1 B_{p,q}=\mathbb{Z}/p$. Moreover, the first Chern class is $p$-torsion, which produces a $\frac{1}{p}\mathbb{Z}$ grading on symplectic cohomology.

The compactification has $-K_\mathbb{P}=\mathcal{O}(p+1)$ whereas the divisor is of degree $p$, and thus the monotonicity constant is $\lambda=\frac{p+1}{p}$ and the principal orbit has degree $\frac{-2}{p}$. Writing $b=p-a$, the two singular points on $\mathbb{P}(a,b,1)$ are of type $\frac{1}{a}(1,b^{-1})$ resp. $\frac{1}{b}(1,a^{-1})$, ensuring that the weight $1$ corresponds to the normal direction to $\mathcal{D}_p$.

\subsubsection{Positive equivariant symplectic cohomology}
There are three types of orbits: the principal orbits, coming in an $S^2$ family, as well as two exceptional orbits of type $\frac{1}{a}$ resp. $\frac{1}{b}$. In $H_1(W;\mathbb{Z})=\mathbb{Z}/p$, they satisfy $a[\delta_{1/a}]=b[\delta_{1/b}]=[M_1]$, but since $a,b$ are coprime to $p$, this shows that the class of the principal orbit generates  $H_1(W;\mathbb{Z})$ and our grading formulas can be applied.

\begin{prop}
    The positive equivariant symplectic cohomology of $B_{p,q}$ is given by $$\mathrm{SH}^{1+\frac{2j}{p}}_{+,S^1}(B_{p,q})=\begin{cases}
        \mathbb{Q}^p, \quad j< 0\\
        \mathbb{Q}^{p-1}, \quad j=0\\
        \mathbb{Q}^{p-j}, \quad 1\leq j \leq p-1\\
        0, \quad j\geq p      
    \end{cases}$$
\end{prop}

\begin{proof}
    The Morse-Bott spectral sequence degenerates at the first page, by Theorem \ref{thm:degeneration}. The principal orbits of period $m$ contribute in degrees $$1-\frac{2m}{p},\quad 3-\frac{2m}{p}$$

    On the other hand, the exceptional orbits of period $\frac{m}{w}, w\in \{a,b\}$, contribute in degree $$-2\left(\frac{2}{p}\frac{m}{w}-\frac{\overline{mv_w}}{w}\right)+1, \quad v_w\equiv p^{-1} \,\mathrm{mod}\, w$$

    Denote by $g=[M_1]$ in $H_1(W;\mathbb{Z})$. Write $l=\overline{mv_w}$. Since $pl\equiv m\,\mathrm{mod}\,w$ we can write $\,m=pl-wj>0$. The contribution thus lies in degree $1+\frac{2j}{p}$ and is in the homology class $-jg$.

    Now, fixing $j\in \mathbb{Z}$, the contributions of the exceptional orbits in degree $1+\frac{2j}{p}$ amount to \begin{gather*}
        \#\{l\in \{1,\dots, a-1|\,pl>aj\}+\#\{l\in \{1,\dots, b-1|\,pl>bj\}=\\
        =(a-1)+(b-1)-\#\{l\in \{1,\dots, a-1|\,pl\leq aj\}-\#\{l\in \{1,\dots, b-1|\,pl\leq bj\}=\\=\begin{cases}
        p-2,\quad j\leq 0\\
            p-2-\underbrace{\big(\lfloor \frac{aj}{p}\rfloor+\lfloor \frac{bj}{p}\rfloor\big)}_{j-1}=p-j-1,\quad 1\leq j \leq p-1
        \end{cases}       
    \end{gather*}
  Adding up these contributions with the ones from the principal orbits leads to the conclusion.  
\end{proof}

\subsubsection{The full symplectic cohomology}
\begin{prop}
The symplectic cohomology of $B_{p,q}$ is given by $$\mathrm{SH}^{k-\frac{2j}{p}}_{j[g]}(B_{p,q})=\begin{cases}
    \mathbb{Q}^{p-1}, \quad k\leq 1\\
    \mathbb{Q}^{j-1}, \quad k=2, j \neq 0 \text{ and } 0 \text{ when } j=0\\
    \mathbb{Q}^j, \quad k=3\\
    0,\quad k\geq 4   
\end{cases}$$
As a graded vector space, it does not depend on $q$.
\end{prop}
\begin{proof}
    First of all, the Gysin sequence for positive equivariant symplectic cohomology, together with the LES between $H,\mathrm{SH}$ and $\mathrm{SH}_+$ shows that $$\mathrm{SH}^1_{ctr}(B_{p,q})=\mathbb{Q}^{p-1}, \quad \mathrm{SH}^{>1}_{ctr}(B_{p,q})=0$$and 
    $$\mathrm{SH}^{3-\frac{2j}{p}}_{j[g]}(B_{p,q})=\mathbb{Q}^{j}, \quad \mathrm{SH}^{>3-\frac{2j}{p}}_{j[g]}(B_{p,q})=0$$

    Finally, we use the long exact sequence from \cite{CieliebakFrauenfelderOancea2010} together with the $\frac{-2}{p}$-periodicity of $\widecheck{\mathrm{SH}}(B_{p,q})$. The connecting maps between symplectic homology and cohomology factor through $H^\bullet(B_{p,q}, \partial B_{p,q})\rightarrow H^\bullet(B_{p,q})$ all vanish over $\mathbb{Q}$ since $B_{p,q}$ is a rational homology ball. This provides us with short exact sequences \begin{equation}\label{eq:SESBpq}
    0\rightarrow \mathrm{SH}^\bullet_{j[g]}(B_{p,q})\rightarrow \widecheck{\mathrm{SH}}^\bullet_{j[g]}(B_{p,q})\rightarrow (\mathrm{SH}^{3-\bullet}_{-j[g]}(B_{p,q}))^\lor \rightarrow 0 \end{equation}
    For $\bullet=1,2, j=0$ this shows that $\widecheck{\mathrm{SH}}_{ctr}^1(B_{p,q})\simeq \mathbb{Q}^{p-1}\simeq \widecheck{\mathrm{SH}}_{ctr}^0(B_{p,q})$. Using Theorem \ref{thm:monotonecompactifications}, this determines the whole Rabinowitz Floer cohomology via $\frac{-2}{p}$-periodicity as $$\widecheck{\mathrm{SH}}^{k-\frac{2j}{p}}_{j[g]}(B_{p,q})=\mathbb{Q}^{p-1}$$
    Finally, using the short exact sequences \ref{eq:SESBpq}, one obtains that $$\mathrm{SH}^{k-\frac{2j}{p}}_{j[g]}=\mathbb{Q}^{p-1}, k\leq 1\quad \mathrm{SH}^{2-\frac{2j}{p}}_{j[g]}=\mathbb{Q}^{j-1}$$
    This is because when $k\leq 1$, the group $\mathrm{SH}^{3-(k-\frac{2j}{p})}_{-j[g]}=0$, whereas when $k=2$, $$\mathrm{SH}^{3-(2-\frac{2j}{p})}_{-j[g]}=\mathrm{SH}^{3-\frac{2(p-j)}{p})}_{(p-j)[g]}=\mathbb{Q}^{p-j}$$

\end{proof}

\begin{rem}
    There is a presentation of the wrapped Fukaya category of $B_{p,1}$ in \cite{Kar18}. However, extracting the Hochschild invariants of a Chekanov-Eliashberg dg algebra (which would describe the symplectic cohomology from the data of the Fukaya category) is a difficult algebraic problem and it is not described anywhere to the authors knowledge.
\end{rem}

\subsection{An example of J. Wahl}
\label{sec:ssing2}
The $B_{p,q}$'s are rational homology smoothings of hyperquotient singularities which compactify to weighted projective spaces. In fact, there is a similar class of surface singularities like this which appears in \cite[Example 5.9.2]{Wahl1981} and \cite{kollar_is_2008}: the defining equation is
\begin{equation*}
x^ay+y^bz+z^cx=0.
\end{equation*}
This is quasihomogeneous with weights $(w_1,w_2,w_3):=(bc-c+1, ca-a+1,ab-b+1)$ and degree $d=abc+1$. This singularity carries a $\mathbb{Z}/d$ action and the Milnor fiber of its quotient is a rational homology ball. We describe it as follows: consider the family of weighted affine hypersurfaces \begin{equation*} x^ay+y^bz+z^cx+\varepsilon  t=0 \subset \mathbb{P}(w_1,w_2,w_3,d)\setminus \{t=0\}\end{equation*} At $\varepsilon=0$ we get precisely the hyperquotient $\{x^ay+y^bz+z^cx=0\}/\frac{1}{d}(w_1,w_2,w_3)$. Note that $\frac{1}{d}(w_1,w_2,w_3)\equiv \frac{1}{d}(1,-a,ab)$ after multiplying by $w_1$. 

At $\varepsilon\neq 0$ we can eliminate $t$ to get \begin{equation*}\mathbb{P}(w_1,w_2,w_3)\setminus \mathcal{D}_{abc+1}\end{equation*} where $\mathcal{D}_{abc+1}$ is a degree $d$ curve in weighted projective space defined by $x^ay+y^bz+z^cx=0$. When $(w_1,w_2,w_3)$ are pairwise coprime, the curve $\mathcal{D}_{abc+1}$ is of genus $0$ and the complement is a rational homology plane (see \cite[41]{kollar_is_2008}).

\begin{rem}
    When $a=1$, the singularity $xy+y^bz+z^cx=0$ is actually analytically equivalent to the $A_{bc}$ singularity. However, the link has a Seifert structure which is different than the one explored in Example 11.1 We focus on the case $a,b,c>1$.
\end{rem}

Since $H^\bullet(I\mathcal{D})$ is concentrated in even degrees, we immediately see by Theorem \ref{thm:degeneration} that the positive equivariant symplectic cohomology of $\mathbb{P}(w_1,w_2,w_3)\setminus \mathcal{D}$ can be computed from the first page of the Morse-Bott spectral sequence. We will compute enough of it to show the following:
\begin{prop}
    Let $W$ be the $\mathbb{Q}$-homology plane obtained by removing $\mathcal{D}_{abc+1}$ from $\mathbb{P}(w_1,w_2,w_3)$ with $w_1,w_2,w_3$ pairwise coprime as above. Assume that each of the $a,b,c$ are greater than $1$. Then,\begin{equation*}\mathrm{SH}^0_{ctr}(W)=\mathbb{Q}, \quad \mathrm{SH}^1_{ctr}(W)=\mathrm{SH}^{<0}_{ctr}(W)=0.\end{equation*}
    Moreover, \begin{equation*}\mathrm{SH}^\bullet_{S^1,ctr}(W)=H^\bullet_{S^1}(W)\oplus \mathrm{SH}^\bullet_{S^1,ctr,+}(W)\end{equation*}
    as $u$-modules, where $H^\bullet_{S^1}(W)\simeq \mathbb{Q}((u))/ u\mathbb{Q}[[u]]$ and $|u|=2$.
\end{prop}

\begin{proof}
    We focus on the contractible part of symplectic cohomology again: one has $\pi_1(W)=\mathbb{Z}/d$ and $c_1^{TW}\in H^2(W)=\mathbb{Z}/d$ is $d$-torsion, so there is a $\frac{1}{d}\mathbb{Z}$-grading. Moreover, it is not difficult to show that the principal orbit, as well as each of the exceptional orbits, each generate $H_1(W)$.

    The contractible part of $\mathrm{SH}^\bullet_{+,S^1}(W)$ then will have contributions from the principal orbits, as well as the exceptional orbits which asymptote to the cyclic quotient singularities $\frac{1}{w_1}(1,w_3w_2^{-1})$ etc. We will show that these contributions are all in degrees $\geq 3$: this immediately implies the second part of the Proposition and, combining it with the Gysin sequence and the fact that $H^\bullet(W;\mathbb{Q})=\mathbb{Q}[0]$, one obtains the results about $\mathrm{SH}(W)$.

    We proceed to compute $\mathrm{SH}^\bullet_{+,S^1,ctr}$. Firstly, the $Nd$-th cover of a principal orbit will contribute $H^\bullet(S^2;\mathbb{Q})$ shifted upwards by $2N(d-w_1-w_2-w_3)+1$ since $c_1^\mathbb{P}=\frac{w_1+w_2+w_3}{d}\mathcal{D}$. However, \begin{equation*}d-w_1-w_2-w_3=abc-ab-bc-ac+a+b+c-2=(a-1)(b-1)(c-1)-1\geq 1\end{equation*}
    by our assumption that none of the $a,b,c$ is equal to $1$ and moreover the assumption that $w_1,w_2,w_3$ are pairwise coprime, excluding the possibility $(a,b,c)=(2,2,2)$.

    Secondly, the $Nd$-th cover of an exceptional orbit of type $\frac{1}{w_1}(1,w_3w_2^{-1})$ will contribute $H^\bullet(pt)$ shifted up by $2(Na-N-\lfloor \frac{Nw_3}{w_1} \rfloor)+1$. However, \begin{equation*}Na-N-\left\lfloor \frac{Nw_3}{w_1}\right\rfloor\geq Na-N-\frac{Nw_3}{w_1}=N\left(a-1-\frac{ab-b+1}{bc-c+1}\right)> 0\end{equation*}since \begin{equation*}(a-1)(bc-c+1)-(ab-b+1)=(a-1)(b-1)(c-1)-1> 0.\end{equation*}

    Hence, each of these contributions will be in an odd degree $\geq 3$. The situation is entirely symmetric in $w_1,w_2,w_3$ hence we can conclude that the support of $\mathrm{SH}^\bullet_{+,S^1,ctr}(W)$ is inside the odd integers $\geq 3$.  
\end{proof}

\begin{rem}
    One can in fact get much more information about $\mathrm{SH}^\bullet(W)$ from the degeneration of the Morse-Bott spectral sequence for $\mathrm{SH}_{+,S^1}^\bullet(W)$, however the numerology gets quite intricate so we have not pursued this computation any further. For us, the interesting content is that the symplectic cohomology is non-vanishing.
\end{rem}

\section{K3 surfaces and the complement of an ample orbifold rational curve}
\label{sec:K3}
In this section, we consider the complement of an ample star-shaped configuration $D$ of smooth rational curves in a $K3$ surface $S$. This divisor $D$ consists of a central component, together with various legs attached to it, each of which is just a tower of rational curves. The completion of the complement $S\setminus D$ is a smooth Liouville manifold. Each of the non-central components of $D$ can be contracted, resulting in an orbifold K3 surface $\mathcal{S}$ with singularities of type $A_n$. The divisor $D$ has contracted to an orbifold $\mathbb{P}^1$ i.e. $\mathcal{D}=\mathbb{P}^1_{a_1,\dots, a_d}$ and we instead consider $W$ as the complement of $\mathcal{D}$ in $\mathcal{S}$.
\begin{rem}
    For $\mathcal{D}$ to be ample, it is necessary that $0<\mathcal{D}^2=-\mathcal{D} (-K_\mathcal{S}-\mathcal{D})=-\chi^{orb}(\mathcal{D})=-2+\sum \frac{a_i-1}{a_i}$. This requires $d\geq 3$ and is guaranteed when $d>4$.
\end{rem}

A priori, we cannot use \ref{prop:closednessofprincipalorbit} to show that $\check{M}_1$ is closed and apply Theorem \ref{thm:mainthm}. Instead, we find a different way to verify the closedness of the principal orbit. In this example, it is of degree $2$. We first compute $\mathrm{SH}^2$ using positive equivariant $\mathrm{SH}$ and the Gysin sequence. Then, we verify that in our chain model, $\mathrm{CF}^2(H)$ has the same rank as $\mathrm{SH}^2$ (for $H$ of sufficiently high slope), showing that the differential must vanish on it. We then proceed to compute $\mathrm{SH}^\bullet(W)$.
\subsection{Some topology}
\begin{prop}
    Let $S$ be a $K3$ surface and $D$ a star-shaped configuration of rational curves such that $\mathcal{D}$ is ample. Then $H_2(D;\mathbb{Z})\rightarrow H_2(S;\mathbb{Z})$ is injective.
\end{prop}
\begin{proof}
    Suppose some linear combination of classes in $H_2(D;\mathbb{Z})$ is nullhomologous in $H_2(S;\mathbb{Z})$. Write it as $$0=\alpha_0[C_0]+\sum_{j=1}^d\sum^{a_j-1}_{i=1}\alpha_{i,j} [C_{i,j}]$$where $C_0$ is the central curve and $C_{i,j}$ is the $i$-th curve on the $j$-th leg, with $C_{i,1}$ being adjacent to $C_0$. For a fixed $j$, we can intersect this relation with $[C_{j, a_j-1}], [C_{j, a_j-2}],\dots, [C_{j,1}]$. As a consequence, $$\alpha_{i, a_i-2}=2\alpha_{i,a_i-1},\quad \alpha_{i, a_i-3}=3\alpha_{i,a_i-1}, \dots \quad\alpha_{i, 1}=(a_i-1)\alpha_{i,a_i-1}, \quad \alpha_0=a_i\alpha_{i,a_i-1}$$ for each $i$. Set $x_i:=\alpha_{i,a_i-1}$. We have that $$\alpha_0=a_1x_1=a_2x_2=\dots=a_dx_d$$
    Finally, intersect with $[C_0]$. The relation is $$-2\alpha_0+\sum_{i=1}^d (a_i-1)x_i=0$$
    If $\alpha_0=0,$ this forces all $x_i=0$.  So let $\alpha_0\neq 0$. In particular, since $x_i=\frac{\alpha_0}{a_i}$, we get that $$\mathcal{D}^2=-2+\sum \frac{a_i-1}{a_i}=0$$ which contradicts the assumption that $\mathcal{D}$ is ample. 

    Note that $H_2(D;\mathbb{Z})=\mathbb{Z}^{1+\sum_{i=1}^d(a_i-1)}$ and $H_2(S;\mathbb{Z})=\mathbb{Z}^{22}$. Set $a=\sum_{i=1}^d (a_i-1)$. Then, since $K3$ surfaces are simply-connected, $ H^2(S\setminus D; \mathbb{Q})\simeq \mathbb{Q}^{22-a-1}$. 

    Similarly, $H_2(S\setminus D;\mathbb{Z})\simeq H^2(S,D)$ is the kernel of the restriction $H^2(S;\mathbb{Z})\rightarrow H^2(D,\mathbb{Z})$ so it is free of rank $22-(a+1)$.
\end{proof}
\subsection{Positive equivariant $\mathrm{SH}$}
\begin{prop}
\label{star_comp_pos_eq_SH}
    Let $W$ be as above and assume $H_1(W;\mathbb{Z})$ is either $0$ or generated by the image of a simple principal orbit in $H_1(\partial W;\mathbb{Z})$. Let $a=\sum_{i=1}^da_i-d$. Then 
    \begin{equation*}
    \mathrm{SH}^\bullet_{+,S^1}(W)=\begin{cases}
        \mathbb{Q}^{a+1}&\bullet =3\\
        \mathbb{Q}^{a+2}&\bullet >3 \text{ odd}
    \end{cases}
    \end{equation*}
\end{prop}
\begin{proof}
    We use the assumption on $H_1(W)$ to ensure that our formulas from Section \ref{sec:grading_formulae} are applicable. 
    The $N$-th cover of exceptional orbit corresponding to a $\frac{1}{a_s}(1,a_s-1)$ point has Robbin-Salamon index (by \ref{sec:grading_formulae})
    \begin{equation*}
        \mu_{CZ}(\delta_{N/a_s})=2\left(-\frac{\overline{N}+a_s-\overline{N}}{a_s}-\left\lfloor \frac{N}{a_s}\right\rfloor\right)+1=-1-2\left\lfloor \frac{N}{a_s}\right\rfloor.
    \end{equation*}

    In positive equivariant symplectic cohomology, each of these will contribute once in degree $-(\mu_{CZ}-(n-1))+1=2-\mu_{CZ}=3+2\lfloor \frac{N}{a_s}\rfloor$. In particular, the exceptional orbits contribute $\mathbb{Q}^a$ in each odd degree $\geq3$.

    The principal orbits come in a $S^2$ family. The Robbin-Salamon index is $\mu_{CZ}(\gamma_N)=-2N$ and so they contribute $H^\bullet(S^2)$ for each $N$, shifted in degree by $2N+1$. This shows that the spectral sequence degenerates at the first page and the proposition follows.

\end{proof}

\subsection{Constructing $\mathfrak{M}$}
We utilize the criterion from Proposition \ref{prop:closednessusinglocalFloer} to construct a cocycle $\mathfrak{M}$ of degree $2$ so that we can apply Theorem \ref{thm:mainthm}. We do this by computing $\mathrm{SH}^2$ and comparing it to the rank of $\mathrm{CF}^2(H)$ for an admissible Hamiltonian $H$ for symplectic cohomology. We also use the computation of a finite portion of $\mathrm{SH}(W)$, together with the periodicity Theorem \ref{thm:mainthm}, to compute the whole symplectic cohomology.
\begin{prop}
    For $W$ as before, one has $$\mathrm{SH}^0(W)=\mathbb{Q}, \quad \mathrm{SH}^1(W)=0, \quad \mathrm{SH}^2(W)=\mathbb{Q}^{22}$$
\end{prop}
\begin{proof}
By Proposition \ref{star_comp_pos_eq_SH}, along with the Gysin sequence for positive equivariant symplectic cohomology, we see that $\mathrm{SH}^{\bullet\leq 1}_+(W)=0$. From the long exact sequence between positive symplectic cohomology and symplectic cohomology, we see that $\mathbb{Q}=H^0(W)\simeq \mathrm{SH}^0(W), 0=H^1(W)\simeq \mathrm{SH}^1(W)$. The last follows since $H^1(W)=H_3(S,D)$ and $H_2(D)\rightarrow H_2(S)$ is injective.

On the other hand, by the Gysin sequence we have $\mathrm{SH}^2_+(W)\simeq \mathrm{SH}^3_{+,S^1}(W)=\mathbb{Q}^{a+1}$. There is moreover a short exact sequence \begin{equation*}H^2(W;\mathbb{Q})\rightarrow \mathrm{SH}^2(W)\rightarrow \mathrm{SH}^2_+(W)\end{equation*}
This implies that $H^2(W;\mathbb{Q})$ injects into $\mathrm{SH}^2(W)$ and dually, that $\mathrm{SH}_2(W)$ surjects onto $H_2(W;\mathbb{Q})$. However, $H^2(W;\mathbb{Q})\simeq H_2(S,D)\simeq H_2(S;\mathbb{Q})/\mathrm{im} H_2(D;\mathbb{Q})$ which has rank $22-(a+1)$ so we can conclude that $\mathrm{SH}^2(W)$ has rank $22$. 
\end{proof}

\begin{prop}
    For a Hamiltonian $H$ for symplectic cohomology of sufficiently large slope, the first page of the Morse-Bott spectral sequence in degree 2 satisfies $$\bigoplus_{p+q=2}E_1^{p,q}(H)\simeq \mathbb{Q}^{22}\simeq \mathrm{SH}^2(W)$$ so that it already computes the symplectic cohomology in that degree. Hence, the criterion from Proposition \ref{prop:closednessusinglocalFloer} applies.
\end{prop}
\begin{proof}   
    This is straightforward, but we give a detailed proof for completeness and moreover in order to explain where the higher-action corrections come into play.
    
    Take an autonomous Hamiltonian $H^\tau$ of sufficiently high slope so that all Reeb orbits of period $\leq 1$ are accounted for by their Hamiltonian counterparts. Then we apply a Morse-Bott perturbation. One subtlety is that in the construction of our Morse-Bott perturbations of $H^\tau$, as in \cite{bourgeois_morse-bott_2003}, one needs to ensure that each of the orbifold points on the sphere $\partial W/S^1 \simeq \mathbb{P}^1_{a_1,\dots, a_d}$ are accounted for by at least one critical point. It is possible to choose Morse functions which have exactly one critical point for each of these orbifold points, but that necessarily forces the existence of odd-index critical points, i.e. the Morse function cannot be perfect. This is where the higher-action corrections $c$ as in \ref{prop:closednessusinglocalFloer} come in. 
    After perturbation, the orbits contributing in degree $2$ are given by (using the index and grading formulas from before): $$\mathrm{CF}^2(H')=\mathbb{Q}\langle \check{M}_1, \check{\delta}_{1/a_1}, \dots, \check{\delta}_{(a_1-1)/a_1}, \check{\delta}_{1/a_2},\dots, \check{\delta}_{(a_d-1)/a_d} \rangle \oplus C$$ where $C$ is generated by the following types of orbits:
    \begin{itemize}
        \item Other orbits in the principal family which are associated to an index $2$ critical point on $\mathcal{D}$.
        \item Constant orbits corresponding to critical points of index $2$ of a Morse function on $W$.
    \end{itemize}
    After passing to the first page of the Morse-Bott spectral sequence, the first type of orbits disappear, since by definition of the principal orbit, it is the one corresponding to the generator of $H^0(\partial W;\mathbb{Q})$ in the local Floer homology of the family of principal orbits. Moreover, the second type of orbits contribute $H^2(W;\mathbb{Q})\simeq \mathbb{Q}^{22-(a+1)}$. All in all, the contributions add up to exactly $\mathbb{Q}^{22}$ hence we can conclude by quoting the previous proposition.
\end{proof}
\subsection{The full symplectic cohomology}
\begin{prop}
Let $W$ be as above. Then \begin{equation*}
        \mathrm{SH}^\bullet(W)=\begin{cases}
            \mathbb{Q} &\bullet=0
            \\ 0 &\bullet=1
            \\ \mathbb{Q}^{22}&\bullet=2
            \\ \mathbb{Q}^{a}&\bullet=3
            \\\mathbb{Q}^{a+1}&\bullet\geq 4
        \end{cases}
    \end{equation*}
    Moreover, $\widecheck{\mathrm{SH}}$ is of constant rank equal to $a+1=\mathrm{rk}\,H_2(D)$.
\end{prop}
\begin{proof}
    In the previous propositions, we computed $\mathrm{SH}$ up to degree $2$.
    
    For the higher degrees, we use Theorem \ref{thm:mainthm} which shows that $\widecheck{\mathrm{SH}}(W)$ is $2$-periodic, together with the long exact sequence between symplectic cohomology and V-shaped symplectic cohomology. We can read off that 
        \begin{equation*}\mathrm{SH}^\bullet(W)\simeq \widecheck{\mathrm{SH}}^\bullet(W), \bullet \geq 4\end{equation*}
        as well as that there are exact sequences \begin{gather}
            0\rightarrow \mathrm{SH}^3\rightarrow \widecheck{\mathrm{SH}}^3\rightarrow \mathrm{SH}_0\rightarrow 0\\
            0\rightarrow \widecheck{\mathrm{SH}}^1\rightarrow \mathrm{SH}_2\rightarrow \mathrm{SH}^2\rightarrow \widecheck{\mathrm{SH}}^2\rightarrow 0\\
            0\rightarrow \mathrm{SH}^0\rightarrow \widecheck{\mathrm{SH}}^0\rightarrow \mathrm{SH}_3\rightarrow 0.
        \end{gather}
    Since $\mathrm{rk}\,\mathrm{SH}_2=\mathrm{rk}\,\mathrm{SH}^2$, this shows that $\mathrm{rk}\,\widecheck{\mathrm{SH}}^1=\mathrm{rk}\,\widecheck{\mathrm{SH}}^2$. Since V-shaped symplectic cohomology is $2$-periodic, this implies that it must be of constant rank in each degree. The map $\mathrm{SH}_2\rightarrow \mathrm{SH}^2$ factors through $H_2(W)\simeq H^2(W,\partial W)\rightarrow H^2(W)$, which is an isomorphism (over $\mathbb{Q}$). Moreover, we know that $H^2(W)\rightarrow \mathrm{SH}^2(W)$ is injective from before, and dually that $\mathrm{SH}_2(W)\rightarrow H_2(W)$ is surjective. Therefore, \begin{equation*}\widecheck{\mathrm{SH}}^2\simeq \mathrm{SH}^2/\mathrm{im} H^2(W)\simeq \mathrm{SH}^2_+(W)=\mathbb{Q}^{a+1}.\end{equation*}
    Finally, by $2$-periodicity, this also holds for $\widecheck{\mathrm{SH}}^3(W)\simeq \mathrm{SH}_0(W)\oplus \mathrm{SH}^3(W)$, which implies that $\mathrm{SH}^3(W)=\mathbb{Q}^a$, completing the proof.
\end{proof}
\subsection{Examples of orbifold $\mathbb{P}^1$'s in weighted K3 surfaces}
We present here some constructions explicitly describing pairs $(\mathcal{S},\mathcal{D})$.
\begin{eg}[Exceptional unimodal singularities]

    Each of the $14$ exceptional unimodal singularities described by Dolgachev \cite{Dolgachev1974-or} is quasihomogeneous. Their Milnor fiber can be compactified to an orbifold $K3$ surface with type $A$ orbifold points, with the compactifying divisor being an orbifold $\mathbb{P}^1$ with three orbifold points whose isotropies are the triple of Dolgachev numbers associated to the singularity. For example, consider the hypersurface singularity $E_{13}$ defined by $x^3+xy^5+z^2$. The Milnor fiber $W$ can be compactified to the weighted K3 surface $\mathcal{X}_{30}\subset \mathbb{P}(10,4,15,1)$ and the compactifying curve is a hyperplane section isomorphic to $\mathbb{P}^1_{2,4,5}$. We have
        \begin{equation*}
        \mathrm{SH}^\bullet(W)=\begin{cases}
            \mathbb{Q}&\bullet=0,
            \\ 0&\bullet=1,
            \\ \mathbb{Q}^{22}&\bullet=2,
            \\ \mathbb{Q}^{8}&\bullet=3,
            \\\mathbb{Q}^{9}&\bullet\geq 4.
        \end{cases}
    \end{equation*}
\end{eg}
\begin{eg}
    We present also an example with more than three orbifold points.
    Consider the weighted K3 surface $\mathcal{X}_{15}\subset \mathbb{P}(2,3,5,5)$ which contains the rational orbifold curve given as the complete intersection $\mathcal{C}_{3,15}\simeq \mathbb{P}^1_{2,5,5,5}$. An example is given by taking the zero locus of the polynomial $y^5+x^5z+z^3+w^3$ and the boundary divisor is given by the zero locus of $y$. The complement $W$ is given by the affine variety $\{x^5z+z^3+w^3=-1\}\subset \mathbb{A}^3/\mu_3$, which is a free $\mu_3$ quotient of the Milnor fiber of $x^5z+z^3+w^3$. As such, $\pi_1(W)=H_1(W)=\mathbb{Z}/3$ and $W$ satisfies all of our assumptions: in particular, the generic orbit in $\partial W$ satisfies $2[\delta_{1/2}]=[\check{M}_1]\in H_1(\partial W)$ and $5[\delta_{1/5}]=[\check{M}_1]\in H_1(\partial W)$ (for the exceptional orbits of type $\frac{1}{5}$) and hence the same relations hold in $H_1(W)$, showing that the principal orbit generates all orbits in $H_1(W)$.
    \begin{equation*}
        \mathrm{SH}^\bullet(W)=\begin{cases}
            \mathbb{Q}&\bullet=0,
            \\ 0&\bullet=1,
            \\ \mathbb{Q}^{22}&\bullet=2,
            \\ \mathbb{Q}^{13}&\bullet=3,
            \\\mathbb{Q}^{14}&\bullet\geq 4.
        \end{cases}
    \end{equation*}
    
\end{eg}

\section{Three-dimensional log Fano examples}\label{sec:logFanocomputations}

\subsection{A three-dimensional singularity of Thom-Sebastiani type $D_{2m+2}\oplus D_{2m+2}$ }

Consider $$f_m=x^2z+y^2w+z^{2m+1}+w^{2m+1}$$ which defines a canonical three-fold singularity. The Milnor number of $f_m$ is $(2m+2)^2$.
The Milnor fiber $W_{2m+1}$ of $f_m$ compactifies as $\mathcal{X}_{2m+1}\subset\mathbb{P}(m,m,1,1,1)$ and the compactifying divisor $\mathcal{D}_{2m+1}$ is an orbifold surface whose stabilizer locus is the locus $z=w=0$ and forms a whole $\mathbb{P}^1$ with $\frac{1}{m}(0,1,1)$ singularities.

The monotonicity index is $\lambda=2$ since $-K_\mathcal{X}=\mathcal{O}(2)$ and the divisor is a hyperplane section. The Betti numbers of $\mathcal{D}_{2m+1}$ are $(b_0,b_2,b_4)=(1,2m+5,1)$, using Dolgachev's extension \cite{Dolgachev1982} of Griffiths residue formula to weighted projective varieties. The link, which is the boundary of $W_{2m+1}$, is simply-connected and has integral homology $H_2(\partial W_{2m+1};\mathbb{Z})=\mathbb{Z}^{2m+4}$ using \cite{Kollar2006} and hence is diffeomorphic to $\# _{2m+4}(S^2\times S^3)$.

\begin{prop}
    The positive equivariant symplectic cohomology of $W_{2m+1}$ is given by $$\mathrm{SH}_{+,S^1}^\bullet(W_{2m+1})=\begin{cases}
        \mathbb{Q}^m, \bullet =3\\
        \mathbb{Q}^{4m+4}, \bullet =1\\
        \mathbb{Q}^{4m+5}, \bullet <0 \text{ odd }\\
        0, \text{ otherwise}
    \end{cases}$$
\end{prop}
\begin{proof}
    The underlying space of $I\mathcal{D}$ has only even cohomology, so Theorem \ref{thm:degeneration} can be applied to compute the positive equivariant symplectic cohomology. The principal orbits of period $l$ contribute one copy in degree $1-2l$, $b_2=2m+5$ copies in degree $3-2l$ and one copy in degree $5-2l$. The exceptional orbits $\delta_{l/m}$ contribute once in degree $1-2\lfloor \frac{l}{m}\rfloor$ and once in degree $3-2\lfloor \frac{l}{m}\rfloor$.
\end{proof}

\begin{prop}
    The symplectic cohomology of $W_{2m+1}$ is given by $$\mathrm{SH}^\bullet(W_{2m+1})=\begin{cases}
        \mathbb{Q}^{m+(2m+2)^2}, \bullet =3\\
        \mathbb{Q}^{3m+4},\bullet \leq 1\\
        0, \text{ otherwise }
    \end{cases}$$
\end{prop}

\begin{proof}
    First, we observe that $\mathrm{SH}^2_+(W_{2m+1})=0$. This follows by running the usual Morse-Bott spectral sequence: the generators, for a Hamiltonian of high slope, in degree $2$ are exceptional orbits of period $<1$ associated to index 2 critical points on a Morse-Bott manifold which is a three-sphere, together with a principal orbit of period $1$ associated to an index $4$ critical point on a Morse-Bott manifold which is $\partial W$. None of these survive after passing to the first page of the MB spectral sequence, since $H^2(S^3;\mathbb{Q})=0$ and $H^4(\partial W;\mathbb{Q})=0$. Therefore, by the long exact sequence between $\mathrm{SH}$ and $\mathrm{SH}_+$, we see that $\mathrm{SH}^2(W_{2m+1})=0$ and $\mathrm{SH}^3(W_{2m+1})\simeq H^3(W_{2m+1};\mathbb{Q})\oplus \mathrm{SH}^3_+(W_{2m+1})$.

    Moreover, by the Gysin sequence for positive equivariant symplectic cohomology, we can read off that $\mathrm{SH}^1_+(W_{2m+1})=\mathbb{Q}^{3m+4}$ and again, by the LES relating positive to usual symplectic cohomology, we get that $\mathrm{SH}^1(W_{2m+1})=\mathbb{Q}^{3m+4}$ as well as that $\mathrm{SH}^\bullet(W_{2m+1})=0, \bullet \geq 4$.

    Now, by Theorem \ref{thm:monotonecompactifications}, the V-shaped symplectic cohomology is $2$-periodic. Moreover, the LES from \cite{CieliebakFrauenfelderOancea2010} shows that $$\mathrm{SH}^{\leq 0}(W_{2m+1})\simeq \widecheck{\mathrm{SH}}^{\leq 0}(W_{2m+1}),\quad \mathrm{SH}^1(W_{2m+1})\simeq \widecheck{\mathrm{SH}}^1(W_{2m+1})$$ as well as $$\widecheck{\mathrm{SH}}^4(W_{2m+1})\simeq \mathrm{SH}_1(W_{2m+1})\simeq \mathbb{Q}^{3m+4}$$

    In particular, $\widecheck{\mathrm{SH}}^\bullet(W_{2m+1})$ is of constant rank $3m+4$ and the result follows.
\end{proof}

\subsection{The cone on the Johnson-Kollár surface}

\begin{prop}
    Consider the cone on the Johnson-Kollár surface $X_{8k+4}\subset \mathbb{P}(2,2k+1,2k+1,4k+1)$. This is a canonical 3-fold singularity whose Milnor fiber has $\mu=3(12k+9)$. Moreover, the link of this singularity is diffeomorphic to $(S^2\times S^3)\#_7$. The positive equivariant symplectic cohomology is lacunary, and hence an invariant of the contact structure of the boundary (\cite{gutt2017positive}).
    \begin{equation*}
    \mathrm{SH}_{+,S^1}^\bullet\cong\begin{cases}
        \mathbb{Q}^{6k+1}& \bullet=3,
        \\\mathbb{Q}^{12k+9}&\bullet=1,
        \\\mathbb{Q}^{12k+10}&\bullet=-2r-1,\; r\in\mathbb{Z}_{\geq 0}.
    \end{cases}
    \end{equation*}

    Therefore, we obtain an infinite family of inequivalent non-index positive contact structures on $(S^2\times S^3)\#_7$ and their mean Euler characteristic satisfies \begin{equation*}
        \chi_m=\frac{-1}{2} \mathrm{rk}\,H^\bullet_{CR}(X_{8k+4})=-(6k+5)
    \end{equation*} 
\end{prop}
\begin{proof}
    Firstly, the surfaces $X_{8k+4}$ have a single $\frac{1}{4k+1}(1,1)$ point as well as four $\frac{1}{2k+1}(1,k)$ points. They have Picard rank $8$ (\cite{Gugiatti2023FullEC}) hence the Seifert bundle associated to the line bundle $\mathcal{O}(1)$ will be a simply-connected $5$ manifold with $H^2(\partial W)=\mathbb{Z}^7$, with no torsion in its homology (this is because torsion on Seifert manifolds occurs from codimension 1 orbifold loci, as in Kollár's work \cite{Kollar2006}). Therefore, the link of the singularity is diffeomorphic to $(S^2\times S^3)_{\#7}$ by Smale's classification.

    To compute the positive equivariant $SH$, we use the formula from Section 4. Firstly, we note that the singularities in the threefold compactifying the Milnor fiber are of type $\frac{1}{4k+1}(1,2k+1,2k+1)$ and $4\times \frac{1}{2k+1}(1,2,2k)$ with the weight $1$ corresponding to the normal direction of the compactifying divisor. So the contribution from the $i$-th cover of an exceptional orbit corresponding to a $\frac{1}{2k+1}(1,2,2k)$ point is  \begin{equation*}
       -2\left(\frac{i}{2k+1}-\frac{\overline{2i}}{2k+1}-\frac{\overline{2ki}}{2k+1}\right)+1
    \end{equation*} 
    For the $i$-th cover of an exceptional orbit corresponding to an orbifold point of type $\frac{1}{4k+1}$ the contribution is in degree \begin{equation*}
        -2\left(\frac{i}{4k+1}-\frac{\overline{(2k+1)i}}{4k+1}-\frac{\overline{(2k+1)i}}{4k+1}\right)+1
    \end{equation*}
    for an orbit of type $\frac{1}{4k+1}$. Finally, the $i$-th covers of principal orbits contribute to degrees $-2i+1$ resp. $-2i+3$ resp. $-2i+5$. Adding all of these up, which turn out to be all odd integers, we get the result.
\end{proof}

\section{Calabi-Yau threefold compactifications}\label{sec:CY3}

\begin{prop}
    Suppose $W$ is a smooth, simply-connected Liouville domain which compactifies, via a simply-connected divisor $\mathcal{D}$, to a three-dimensional orbifold $\mathcal{X}$ with isolated orbifold points and $c_1^\mathcal{X}=0$. Suppose that:\begin{itemize}
        \item $H^4(W)=0$
        \item $H^3\rightarrow \mathrm{SH}^3(W)$ is injective
        \item every Reeb orbit of $\partial W$ (with its canonical contact structure as a Seifert bundle over $\mathcal{D}$) corresponding to a twisted sector of $\mathcal{X}$ of age one has smaller fractional period than the Reeb orbits corresponding to twisted sectors of age two.
    \end{itemize}  Then: $$\mathrm{SH}^\bullet(W)=\begin{cases}
        0, \bullet <0 \text{ or } \bullet =1\\
        \mathbb{Q}, \bullet =0\\
        \mathbb{Q}^{b_2(W)+1+\alpha_1}, \bullet =2\\
        \mathbb{Q}^{b_3(W)+\alpha_1}, \bullet =3\\
        \mathbb{Q}^{b_2(\mathcal{D})+2\alpha_1}, \bullet =4 \text{ or } \bullet \geq 6\\
        \mathbb{Q}^{b_2(\mathcal{D})+2\alpha_1-1}, \bullet =5
    \end{cases}$$
    where $\alpha_1$ is the number of non-identity twisted sectors on $\mathcal{X}$ with age $1$.
\end{prop}

\begin{proof}
    We first focus on the positive equivariant cohomology. The principal orbits of period $l$ contribute once in degree $1+2l$, $b_2(\mathcal{D})$ many times in degree $3+2l$ and once in degree $5+2l$. The exceptional orbits of period $\frac{l}{w}$ contribute once in degree $1+2(\iota_{age}+\lfloor \frac{l}{w}\rfloor)$. In particular, the positive equivariant symplectic cohomology is concentrated in odd degrees $\geq 3$: $$\mathrm{SH}^1_{+,S^1}(W)=0,\quad \mathrm{SH}^3_{+,S^1}(W)=\mathbb{Q}^{\alpha_1+1}$$
    Thus, by the Gysin sequence, $\mathrm{SH}^0_+(W)=0$ and so $H^0(W)\simeq \mathrm{SH}^0(W)$. Moreover, $\mathrm{SH}^1_+(W)=0=\mathrm{SH}^1(W)$ . Again by the Gysin sequence, $\mathrm{SH}^2_+(W)=\mathbb{Q}^{\alpha_1+1}$, which is generated by the principal orbit, together with $\alpha_1$ many exceptional orbits. By assumption, $H^3(W)\rightarrow \mathrm{SH}^3(W)$ is injective, so there is a short exact sequence $$0\rightarrow H^2(W)\rightarrow \mathrm{SH}^2(W)\rightarrow \mathrm{SH}^2_+(W)\rightarrow 0$$

    In particular, we see that the principal orbit survives to $\mathrm{SH}^2$ and so Theorem \ref{thm:mainthm} can be applied.

    Before we do this, we compute $\mathrm{SH}^3(W)$, which we do by running the Morse-Bott spectral sequence. On the first page, the degree $2$ contributions are given by $H^2(W)$, the minimum of the principal orbit family, as well as the minima of the circle families of age $1$ orbits. In degree $3$, the contributions are from $H^3(W)$, the maxima of the circle families of age $1$ orbits, but no principal orbits, since $H^1(\partial W;\mathbb{Q})=0$. In degree $4$, the contributions are from age $2$ exceptional orbits, together with principal orbits corresponding to $H^2(\partial W;\mathbb{Q})$. In degrees $2$ and $3$, the spectral sequence degenerates at this page, since by assumption, all the orbits in degree $4$ have higher period, hence lower action, than the ones in degree $3$. In particular, $\mathrm{SH}^3(W)=\mathbb{Q}^{b_3(W)}\oplus \mathbb{Q}^{\alpha_1}$.

    Finally, we run the LES from \cite{CieliebakFrauenfelderOancea2010}. This immediately gives $$\widecheck{\mathrm{SH}}^{-1}(W)\simeq \mathrm{SH}_6(W), \quad \widecheck{\mathrm{SH}}^1(W)\simeq \mathrm{SH}_4(W), \quad \mathrm{SH}^4(W)\simeq \widecheck{\mathrm{SH}}^4(W),\quad \mathrm{SH}^{\geq 6}(W)\simeq \widecheck{\mathrm{SH}}^{\geq 6}(W)$$together with a short exact sequence $$0\rightarrow \mathrm{SH}^0(W)\rightarrow \widecheck{\mathrm{SH}}^0(W)\rightarrow \mathrm{SH}_5(W)\rightarrow 0$$
    Since V-shaped symplectic cohomology is $2$-periodic and $\mathrm{SH}_4(W)$ is $\mathbb{Q}$-dual to $\mathrm{SH}^4(W)$, we get that in fact $\widecheck{\mathrm{SH}}^\bullet(W)$ is of constant rank, which we now determine.

    The last piece of the exact sequence is $$0\rightarrow \mathrm{SH}^2(W)\rightarrow \widecheck{\mathrm{SH}}^2(W)\rightarrow \mathrm{SH}_3(W)\rightarrow \mathrm{SH}^3(W)\rightarrow \widecheck{\mathrm{SH}}^3(W)\rightarrow \mathrm{SH}_2(W)\rightarrow 0$$

    The middle connecting map factors through $H_3(W)\rightarrow H^3(W)\simeq H_3(W,\partial W)$. The kernel of this map is identified with the image of $H_3(\partial W)\rightarrow H_3(W)$, which can be read off from $$\underbrace{H_4(W)}_{0}\rightarrow H_4(W,\partial W)\rightarrow H_3(\partial W)\rightarrow H_3(W)$$
    Moreover, $\mathrm{rk}\,H_4(W,\partial W)=\mathrm{rk}\,H^4(W,\partial W)=\mathrm{rk} \,H_2(W)=b_2(W)$.

    All together, the kernel of $\mathrm{SH}_3(W)\rightarrow \mathrm{SH}^3(W)$ is identified, as vector spaces, with $\mathrm{SH}_3^+(W)\oplus \frac{H_3(\partial W)}{\mathrm{im}\,H_4(W,\partial W)}$. Therefore, we conclude that \begin{gather*}     \mathrm{rk}\,\widecheck{\mathrm{SH}}^2(W)=\mathrm{rk}\,\mathrm{SH}^2+\mathrm{rk}\,\mathrm{SH}^3_+(W)+b_3(\partial W)-b_2(W)=\\=b_2(W)+1+\alpha_1+\alpha_1+(b_2(\mathcal{D})-1)-b_2(W)=b_2(\mathcal{D})+2\alpha_1
    \end{gather*} which is, in fact, the rank of $\widecheck{\mathrm{SH}}^\bullet$ in all integer degrees. From this, the rest follows.
\end{proof}



\newpage
\sloppy
\printbibliography
\end{document}